\documentclass[final,5p,times,sort&compress]{elsarticle}
\usepackage{amsmath, amsfonts, amssymb}
\usepackage{graphicx,epsfig}
\usepackage{color}
\usepackage{booktabs}
\usepackage[section]{placeins}
\numberwithin{equation}{section}

\usepackage{array}
\graphicspath{{images2/}}

\usepackage{amsthm}
\theoremstyle{plain}\newtheorem{remark}{Remark}[section]
\theoremstyle{definition}

\journal{Digital Signal Processing}
\usepackage{hyperref}

\begin{document}

\begin{frontmatter}

\title{Linear and synchrosqueezed time-frequency representations revisited.\\ Part II: Resolution, reconstruction and concentration.}
\tnotetext[support]{Work supported by the Engineering and Physical Sciences Research Council, UK}

\author{Dmytro Iatsenko}
\ead{dmytro.iatsenko@gmail.com}

\author{Peter V. E. McClintock\corref{}}
\ead{p.v.e.mcclintock@lancaster.ac.uk}

\author{Aneta Stefanovska\corref{cor1}}
\cortext[cor1]{Corresponding author}
\ead{aneta@lancaster.ac.uk}

\address{Department of Physics, Lancaster University, Lancaster LA1 4YB, UK}

\begin{abstract}

Having reviewed the aspects of the linear and synchrosqueezed time-frequency representations (TFRs) needed for their understanding and correct use in Part I of this review, we now consider three more subtle issues that are nonetheless of crucial importance for effective application of these methods. (i) What effect do the window/wavelet parameters have on the resultant TFR, and how can they most appropriately be chosen? (ii) What are the errors inherent in the two reconstruction methods (direct and ridge) and which of them is the better? (iii) What are the advantages and drawbacks associated with synchrosqueezing? To answer these questions, we perform a detailed numerical and theoretical study of the TFRs under consideration. We consider the relevant estimates in the presence of the complications that arise in practical applications including interference between components, amplitude modulation, frequency modulation, and noise. Taken together, the results provide an in-depth understanding of the issues in question.

\end{abstract}


\end{frontmatter}


\section{Introduction}\label{sec:introd}

Although the main aspects of the windowed Fourier transform (WFT), wavelet transform (WT) and their synchrosqueezed equivalents (SWFT and SWT) were considered in detail in Part I of this review, some important questions still remain to be addressed. Thus, the parameters of the window/wavelet used for the (S)WFT/(S)WT, such as the resolution parameter $f_0$, are often chosen blindly based on established conventions (e.g.\ $f_0=1$ for the Morlet wavelet), without complete understanding of the effects of this and other choices on the resultant TFR and the outcome of the analysis. Next, there are two different methods by which one can reconstruct the parameters of the signal components from its TFR: direct and ridge (see Part I). Both of them are used in the literature, but it has remained unclear which of the two is to be preferred and when. Finally, it is unclear what advantages over the usual WFT/WT are gained by using SWFT/SWT (apart from nicer visual appearance of the latter). It also remains to be established whether or not synchrosqueezing changes the time and/or frequency resolution of the TFR and therefore allows for a more accurate estimation of the components' parameters.

We now study each of these issues, {\it viz.} the effects of the window/wavelet parameters, the relative performance of different reconstruction methods, and the advantages/drawbacks of synchrosqueezing, in four different cases, when the signal is represented as: two interfering tones; an amplitude modulated (AM) component; a frequency modulated (FM) component; and a single tone corrupted by noise. By proceeding in this way, accounting for all possible complications, one can build up quite a complete picture of how the issues in question manifest themselves for an arbitrary signal.

After discussing in Sec.\ \ref{sec:assumptions} the assumptions and conventions used in this work, we introduce in Sec.\ \ref{sec:unifying} a convenient unifying formalism that allows a common set of equations to be used to describe the (S)WFTs and (S)WTs. We lay out the questions to be addressed in Sec.\ \ref{sec:questions}, and in Sec.\ \ref{sec:cases} we consider them in detail in the situation when there are two interfering tones, amplitude modulation, frequency modulation, and noise. In the light of the results obtained, Sec.\ \ref{sec:optres} discusses the choice of optimal window/wavelet parameters, Sec.\ \ref{sec:optrec} considers the performance of different reconstruction methods and their related errors, and Sec.\ \ref{sec:Synchroneed} discusses the extent to which synchrosqueezing is useful, as well as the difference between TFR concentration and resolution. Finally, in Sec.\ \ref{sec:conclusions} we summarise and draw conclusions.

\section{Assumptions and conventions}\label{sec:assumptions}

In the following we assume that the Part I of this work has been read thoroughly, and we use the same notation, terminology and conventions as there (see its Appendix A), e.g.\ the notion of $\epsilon$-supports. However, we now introduce the additional conventions and assumptions that are listed below. Unless otherwise specified, all the following considerations apply only within the assumptions made, e.g.\ most of the discussion is inapplicable if the window function $\hat{g}(\xi)$ is multimodal.

\emph{Assumptions about the form of the window function:} We consider $\hat{g}(\xi)$ to be real, positive (so that $\hat{g}(\xi)=|\hat{g}(\xi)|$) and unimodal, i.e.\ having one dominant peak that decays on both sides of it, with any other peaks being negligible by comparison. This is the most useful and convenient form. We also assume $\hat{g}'(\xi)$ and $\hat{g}''(\xi)$ to be finite for all $\xi$. Note, that we do not assume a finite support for $g(t)$ or $\hat{g}(\xi)$, or their symmetry around the maximum, although window functions symmetric in frequency and time are usually to be preferred.

\emph{Assumptions about the form of the wavelet function:} We consider $\hat{\psi}(\xi)$ to be real, positive and unimodal for $\xi>0$ (the form for $\xi\leq 0$ does not matter since we take only positive frequencies in the WT computation). Therefore, in the following $\hat{\psi}(\xi)=\hat{\psi}^*(\xi)=|\hat{\psi}(\xi)|$. We also assume $\hat{g}'(\xi)$ and $\hat{g}''(\xi)$ to be finite for all $\xi>0$. Note, that we do not assume a finite support for $\psi(t)$ or $\hat{\psi}(\xi)$, or any kind of symmetry around the maximum.

\emph{Calculated TFRs:} For all simulation examples presented in this work we have used Gaussian window for the (S)WFT and Morlet wavelet for the (S)WT. The boundary errors in all TFRs (see Part I) are minimized by padding the signal with exact values (i.e.\ simulating it for a longer period and considering only the central part, with the rest used as padding), but, in all cases considered, predictive padding gives almost the same results. The number of padded values on each side is determined as described in Part I. In view of these issues, we will not discuss reconstruction errors $<0.001$, as they become influenced by boundary effects, and by frequency discretization too in the estimation of frequency from the SWFT/SWT (see below).

\emph{Frequency discretization:} The discretization parameters $\Delta\omega$ and $n_v$, determining the width of the frequency bins for the (S)WFT and (S)WT, respectively, are chosen using the criterion described in Part I. However, when the reconstruction of components from TFR is performed and compared for different window/wavelet parameters, we use much smaller constant $\Delta\omega/2\pi=0.002$ (WFT and SWFT) and $n_v=256$ (WT and SWT): as discussed in Part I, the direct and ridge frequency estimation from synchrosqueezed TFRs (only) can suffer greatly from discretization effects, so that very small bins are needed to remove them from consideration; the chosen values guarantee discretization-related errors to be $\Delta\nu_d/2\pi\leq 0.001$ (SWFT) and $\Delta\nu_d/\nu\leq 0.0015$ (SWT).

\emph{Extraction of the time-frequency support:} In the following simulation examples, we will extract the ridge curve $\omega_p(t)$ by selecting either the highest TFR amplitude peaks at each time -- ``maximum-based'' scheme -- or the peaks that are nearest to the actual component's frequency (which we always know \emph{a priori} for simulated signals) -- the ``frequency-based'' scheme. These approaches are used because of their low computational cost, but note that they work only for the simulated examples, and are inapplicable to real signals (more universal schemes for ridge curve extraction are discussed and compared in \cite{Iatsenko:ridge}). The method will always be specified, and will be selected so as to model the best or most realistic component extraction; when both methods are in principle possible, the effect of each extraction procedure on the results will be discussed.

\emph{Figures:} The yellow regions around the time-averaged TFR amplitudes at each frequency indicate their $\pm\sqrt{2}$ standard deviations in time: in the common case where the time-variation is approximately sinusoidal, this is equivalent to $\pm$oscillations' amplitude.

\section{Unifying formulation}\label{sec:unifying}

The only significant difference between the WFT and WT lies in their linear/logarithmic resolution, so one can generalize all the formulas for them. Therefore, denoting the WFT/WT as $H_s(\omega,t)$, we write
\begin{equation}\label{gformTFR}
\begin{gathered}
H_s(\omega,t)=\int s^{+}(t)h_{u-t}(\omega)dt=\frac{1}{2\pi}\int_0^\infty e^{i\xi t}\hat{s}(\xi)\hat{h}_{\xi}(\omega)d\xi,\\
\begin{aligned}
\mbox{\textbf{WFT:}}\;\;&h_t(\omega)=g(t)e^{-i\omega t},\;\;\hat{h}_\xi(\omega)=\hat{g}(\omega-\xi),\\
\mbox{\textbf{ WT:}}\;\;&h_t(\omega)=\big(\omega/\omega_\psi\big)\psi^*(\omega t/\omega_\psi),\;\;\hat{h}_\xi(\omega)=\hat{\psi}^*(\omega_\psi \xi/\omega).\\
\end{aligned}
\end{gathered}
\end{equation}
The signal's time-domain form can then be reconstructed as (see Part I)
\begin{equation}\label{gformREC}
\begin{gathered}
s^{a}(t)=C_h^{-1}\int H_s(\omega,t)d\mu(\omega),\quad
C_h=\frac{1}{2}\int \hat{h}_\nu(\omega)d\mu(\omega)\\
\begin{aligned}
\mbox{\textbf{WFT:}}\;\;&C_h=C_g,\;\;\mu(\omega)=\omega\in(-\infty,\infty),\\
\mbox{\textbf{ WT:}}\;\;&C_h=C_\psi,\;\;\mu(\omega)=\big(\log\omega\big)\in(-\infty,\infty)\quad \big(\omega\in(0,\infty)\big).\\
\end{aligned}
\end{gathered}
\end{equation}
Additionally, we define
\begin{equation}\label{gformNHQ}
\begin{gathered}
\hat{h}_{\max}\equiv \max_{\xi}\hat{h}_\nu(\xi)=\hat{h}_{\nu}(\nu)=
\left[\begin{array}{l}
\hat{g}(0)\quad\mbox{(WFT)}\\
\hat{\psi}(\omega_\psi)\quad\mbox{(WT)}
\end{array}\right.\\
Q_\nu(\omega)\equiv \frac{C_h^{-1}}{2}\int_{-\infty}^{\mu(\omega)} \hat{h}_{\nu}(\xi)d\mu(\xi)=
\left[\begin{array}{l}
R_g(\omega-\nu)\quad\mbox{(WFT)}\\
1-R_{\psi}(\omega_\psi\nu/\omega)\quad\mbox{(WT)}\\
\end{array}\right.\\
\nu_H(\omega,t)\equiv \partial_t {\rm arg}[H_s(\omega,t)],\;\;\widetilde{Q}_\nu(\omega_1,\omega_2)\equiv Q_\nu(\omega_2)-Q_\nu(\omega_1),\\
\end{gathered}
\end{equation}
where $R_{g,\psi}(\omega)$ are defined in Part I, and we have taken into account that $\hat{g}(\xi)=|g(\xi)|,\hat{\psi}(\xi)=|\hat{\psi}(\xi)|$, according to the assumptions made. Thus, if signal contains the tone $A\cos(\nu t+\varphi)$, then $\widetilde{Q}_\nu(\omega_1,\omega_2)$ represents the relative part of this tone that is contained in the TFR located at frequencies $\omega\in[\omega_1,\omega_2]$.

Finally, in what follows we will also need the expressions for all possible TFR-derived quantities related to the multitone signal, which are:
\begin{equation}\label{gformGEN}
\begin{gathered}
s(t)=\sum_na_n\cos(\nu_n t+\varphi_n)\equiv\sum_na_n\cos\phi_n(t),\\
\Downarrow\\
H_s(\omega,t)=\frac{1}{2}\sum_n a_n\hat{h}_{\nu_n}(\omega)e^{i\phi_n(t)},\\
{\tiny{
|H_s(\omega,t)|^2=\sum_n\frac{a_n^2\hat{h}_{\nu_n}^2(\omega)}{4}+\sum_{n,m>n}\frac{a_na_m \hat{h}_{\nu_n}(\omega)\hat{h}_{\nu_m}(\omega)}{2}\cos \Delta\phi_{nm}(t)
}},\\
{\rm arg}[H_s(\omega,t)]=
\arctan\frac{\sum_n a_n\hat{h}_{\nu_n}(\omega)\sin\phi_n(t)}{\sum_n a_n\hat{h}_{\nu_n}(\omega)\cos\phi_n(t)},\\
\begin{aligned}
\nu_{H}(\omega,t)=&\sum_n \frac{a_n^2\hat{h}_{\nu_n}^2(\omega)}{4|H_s(\omega,t)|^2}\nu_n\\
&+\sum_{n,m>n}\frac{a_na_m \hat{h}_{\nu_n}(\omega)\hat{h}_{\nu_m}(\omega)}{4|H_s(\omega,t)|^2}[\nu_n+\nu_m]\cos\Delta\phi_{nm}(t),
\end{aligned}\\
\end{gathered}
\end{equation}
where $\Delta\phi_{nm}\equiv\phi_n(t)-\phi_m(t)$. Note that, unlike the WFT and WT, which are linear TFRs and therefore satisfy additivity (i.e.\ $H_{s_1(t)+s_2(t)}=H_{s_1(t)}+H_{s_2(t)}$), synchrosqueezing is an inherently nonlinear operation, so that the SWFT/SWT of the sum of signals is not equal to the sum of SWFTs/SWTs for each signal separately. Due to this the expressions for synchrosqueezed TFRs are more complicated and cannot be obtained in a simple form like (\ref{gformGEN}).

\section{The questions to be answered}\label{sec:questions}

In this section, we formulate and discuss the three issues that will be studied below.

\subsection{Choice of the window/wavelet parameters}

\begin{figure*}[t!]
\includegraphics[width=1.0\linewidth]{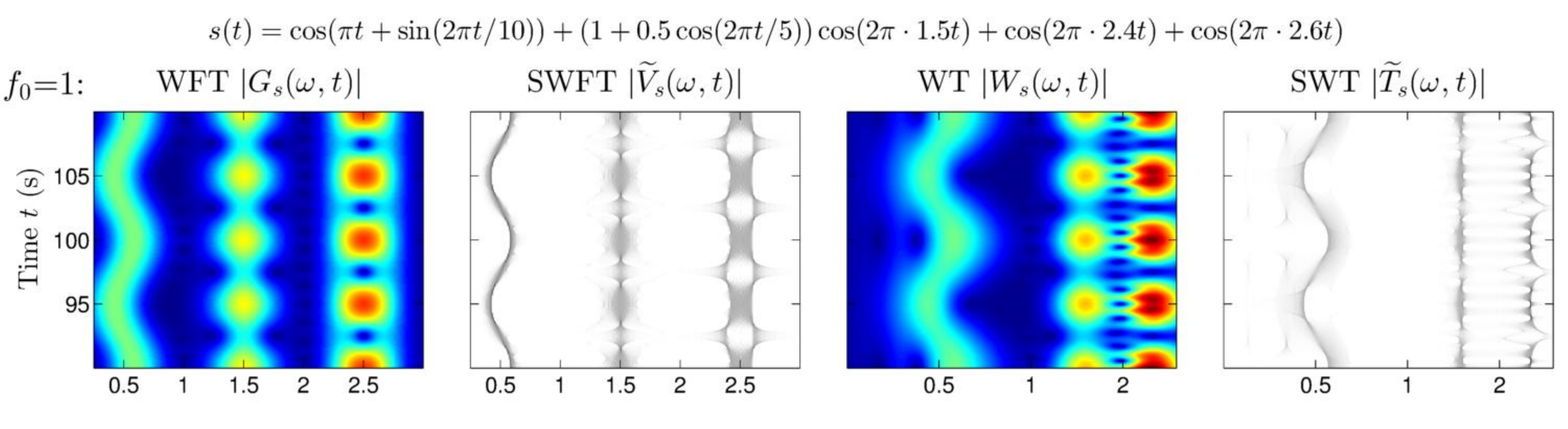}\\
\includegraphics[width=1.0\linewidth]{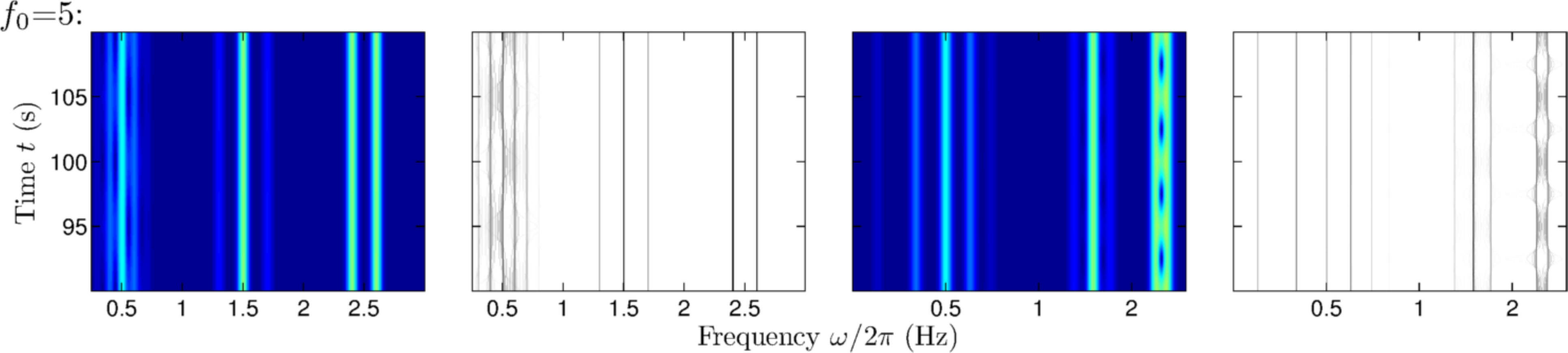}\\
\caption{The WFT, SWFT, WT and SWT calculated at different $f_0$ for a signal consisting of four components $s(t)=s_1(t)+s_2(t)+s_3(t)+s_4(t)$: 1) FM component $s_1(t)=\cos(2\pi\cdot0.5(t+0.25\cos(2\pi\cdot 0.1t)))$; 2) AM component $s_2(t)=(1+0.5\cos(2\pi\cdot0.2t))\cos(2\pi\cdot1.5t)$; 3,4) two tones with close frequencies $s_{3,4}(t)=\cos(2\pi\cdot(2.5\mp0.1)t)$. Note that behavior of the SWFT/SWT is qualitatively the same as for the WFT/WT at each $f_0$. The signal was sampled at 50 Hz for 200 s.}
\label{fig:f0}
\end{figure*}

The parameters of the window/wavelet, e.g.\ the resolution parameter $f_0$, determine the time and frequency resolution of the resultant TFR. Thus, the higher $f_0$ is -- the closer in frequency are the AM/FM components that can be distinguished in a TFR, while at the same time the weaker amplitude and frequency modulations can be represented reliably. The choice of the resolution parameter determines how the TFR treats the AM/FM components present in the signal: either as single entities, or as sums of tones \cite{Putland:00}. This is illustrated in Fig.\ \ref{fig:f0}, where different TFRs are shown for signals consisting of one FM component, one AM component, and two tones. As can be seen, for $f_0=1$, the TFRs treat both AM/FM components as individual entities, while at $f_0=5$ they are represented as sums of tones; at the same time, two tones present in the signal that appear indistinguishable at $f_0=1$ can be resolved at $f_0=5$. Fig.\ \ref{fig:f0} illustrates the importance of selecting appropriate window/wavelet parameters, as well as the inherent difficulties of this task. Thus, for the case considered in the figure, neither $f_0=1$ nor $f_0=5$ is suitable, and one needs to choose the resolution parameter based on a compromise between reliable representation of the AM/FM components and reliable representation of tones.

To choose window/wavelet parameters appropriately, one needs first to understand fully their influence on the TFR properties. For the resolution parameter $f_0$, we wish to answer the questions: how high should $f_0$ be to distinguish between two tones? How low should it be to reliably represent components with particular amplitude and/or frequency modulation, i.e.\ without separating them into multiple TFR lines? Within what range can it be chosen? We investigate these issues in detail by studying the TFR behavior for different signals. Note, that although the simulations are performed for the (S)WFT with Gaussian window and (S)WT with Morlet wavelet, we actually consider the general case of any window/wavelet within the given assumptions, and understand by $f_0$ any set of parameters characterizing its resolution properties.

Because the AM/FM component can be represented as a sum of tones (see Part I), it is useful to consider the WFT/WT of a multitone signal $s(t)=\sum_{m=1}^Ma_m\cos(\nu_m t+\varphi_m)$, and we partition its qualitative behavior into regimes as illustrated in Table \ref{BTtab} below.
\begin{center}
\begin{table}[h]
\caption{Regimes of possible WFT/WT behavior for an $M$-tone signal $s(t)=\sum_{m=1}^Ma_m\cos(\nu_m t+\varphi_m)$.}
\begin{tabular}{| p{1.6cm} | p{6.4cm} |}
\hline
Regime I & All tones are fully resolved, i.e.\ at each time there are $M$ well-separated peaks in the WFT/WT amplitude, and the time-variations of the latter are zero or negligible.\\ \hline
Regime II & Tones partly interfere with each other, so the WFT/WT amplitude varies in time, but there are always $M$ distinct peaks.\\ \hline
Regime III & Tones severely interfere, so that the nearest ones sometimes merge, i.e.\ there exist moments when there are fewer than $M$ (non-negligible) peaks in the WFT/WT amplitude.\\ \hline
Regime IV & Tones are completely merged so that, although the WFT/WT amplitude strongly varies in time, it always has only one dominant peak.\\ \hline
\end{tabular}
\label{BTtab}
\end{table}
\end{center}
Evidently, for a few unrelated tones, one should aim at Regime I,  while for AM/FM components -- at Regime IV. The differences between regimes will become clearer in the forthcoming sections.

To distinguish between different types of WFT/WT behavior, we introduce the mean number of peaks $\langle N_p\rangle$, defined simply as the time-averaged number of peaks in the WFT/WT amplitude
\begin{equation}\label{NP}
\langle N_p\rangle=\langle \#\omega_p(t):\;|\partial_{\omega}H_s(\omega_p(t),t)|=0,\;|\partial_{\omega}^2H_s(\omega_p(t),t)|<0 \rangle
\end{equation}
Thus, $\langle N_p\rangle=M$ indicates Regime I or II, $1<\langle N_p\rangle<M$ implies Regime III, and $\langle N_p\rangle\approx 1$ is characteristic of the IV type of behavior. Note that there might be many small spurious peaks (e.g.\ due to round-off errors), so we discard peaks smaller than $10^{-6}$ of the total summed amplitude of all peaks at each time. In addition, to quantify the interference between tones, we will also introduce the interference measure $\eta$, but it will be defined for each particular case separately.

Regarding behavior of the SWFT/SWT, its classification is more problematic due to its high complexity. For example, even in the noiseless case the synchrosqueezed TFRs often contain much more TFSs than there are tones in the signal. However, as will be seen, the SWFT/SWT and WFT/WT from which it is constructed both behave in qualitatively the same way, so it is enough to classify only behavior of the latter.

\subsection{Choice of the reconstruction method}

As discussed in Part I, there are two possible methods -- direct and ridge -- by which a component can be reconstructed from its support in a TFR. However, to the best of our knowledge, there are no works comparing their performance, so it is not clear in what cases which method should be used. This issue will be thoroughly investigated in the following sections. We will quantify the relative error of reconstruction $\varepsilon_{a,\phi,f}$ of amplitude, phase and frequency of the AM/FM component as
\begin{equation}\label{apfrel}
\begin{aligned}
&\varepsilon_a\equiv\frac{\sqrt{\langle [A_{rec}(t)-A_{true}(t)]^2 \rangle}}{\sqrt{\langle [A_{true}(t)]^2 \rangle}}\\
&\varepsilon_\phi\equiv \sqrt{1-|\langle e^{i(\phi_{rec}(t)-\phi_{true}(t))} \rangle|^2}\\
&\varepsilon_f\equiv
\left[\begin{array}{l}
\frac{\sqrt{\langle [\nu_{rec}(t)-\nu_{true}(t)]^2\rangle}}{2\pi}\mbox{ for (S)WFT}\\
\frac{\sqrt{\langle [\nu_{rec}(t)-\nu_{true}(t)]^2\rangle}}{\langle \nu_{true}(t)\rangle}\mbox{ for (S)WT}\\
\end{array}\right.
\end{aligned}
\end{equation}
where $A_{rec,true}(t),\phi_{rec,true}(t),\nu_{rec,true}(t)$ denote the reconstructed and true parameters. The difference between the definitions of $\varepsilon_f$ for the (S)WFT and (S)WT accounts for the linear and logarithmic frequency resolutions of these transforms. Note that, generally, the reconstruction of phase is much more precise than that of amplitude or frequency, so the form of $\varepsilon_\phi$ was chosen so as to magnify the corresponding error and make it comparable to the others. The analytic expressions for the errors of each reconstruction method will be given at the end, in Sec.\ \ref{sec:optrec}.

\subsection{Advantages/drawbacks of synchrosqueezing}

As can be seen from Fig.\ \ref{fig:f0}, if the components are not reliably represented in the WFT/WT, they will be not well reflected in the SWFT/SWT either (as already noted in Part I). Therefore, it is not immediately obvious what are the advantages of SWFT/SWT over the usual WFT/WT, apart from being more visually appealing. Does synchrosqueezing improve the time or frequency resolution of the transform? Or does it allow more accurate reconstruction of the components? And generally, does the concentration of the TFR represent the primary characteristic of its performance? To answer these questions, we compare not only the performances of different reconstruction methods, but also the accuracy of estimates obtained by these methods from the usual and synchrosqueezed TFRs. For example, to understand whether synchrosqueezing improves frequency resolution it is sufficient to study the case of two interfering components: the resolution can be regarded as being increased only if one is able to extract the parameters of these components more accurately from the SWFT/SWT than from the underlying WFT/WT.

\section{Simulation study}\label{sec:cases}

\subsection{Resolution of two tones}\label{sec:OV}

Consider a two-tone signal
\begin{equation}\label{OVs}
s(t)=A[\cos(\nu_1t+\varphi_1)+r\cos(\nu_2 t+\varphi_2)],\;\;\nu_1<\nu_2,\\
\end{equation}
and for further convenience define
\begin{equation}\label{OVnt}
\Delta\nu\equiv \nu_2-\nu_1;\mbox{  }\Delta\phi(t)\equiv \Delta\nu t+(\varphi_2-\varphi_1);\mbox{  }\bar{\nu}=\frac{\nu_1+\nu_2}{2}.
\end{equation}

Then the corresponding WFT/WT and related quantities will be (\ref{gformGEN}):
\begin{equation}\label{OVamp}
\begin{aligned}
H_s(\omega,t)=&\frac{Ae^{i(\nu_1 t+\varphi_1)}}{2}\left[\hat{h}_{\nu_1}(\omega)+\hat{h}_{\nu_2}(\omega)e^{i\Delta\phi(t)}\right],\\
|H_s(\omega,t)|^2=&\frac{A^2}{4}\left[\hat{h}_{\nu_1}^2(\omega)+\hat{h}_{\nu_2}^2(\omega)+2r\hat{h}_{\nu_1}(\omega)\hat{h}_{\nu_2}(\omega)\cos\Delta\phi(t)\right],\\
\nu_H(\omega,t)=&\frac{\nu_1\hat{h}_{\nu_1}^2(\omega)+\nu_2\hat{h}_{\nu_2}^2(\omega)+2r\bar{\nu}\hat{h}_{\nu_1}(\omega)\hat{h}_{\nu_2}(\omega)\cos\Delta\phi(t)}
{\hat{h}_{\nu_1}^2(\omega)+\hat{h}_{\nu_2}^2(\omega)+2r\hat{h}_{\nu_1}(\omega)\hat{h}_{\nu_2}(\omega)\cos\Delta\phi(t)}.\\
\end{aligned}
\end{equation}
As can be seen, the main supports of the different tones in $H_s(\omega,t)$ might overlap with each other, so that each term is non-negligible for some $\omega$. In this situation one says that the two tones interfere in the TFR, which is reflected in the appearance of terms $\sim \hat{h}_{\nu_1}(\omega)\hat{h}_{\nu_2}(\omega)\cos\Delta\phi(t)$, causing harmonic oscillations in the squared TFR amplitude $|H_s(\omega,t)|^2$ and an instantaneous frequency $\nu_H(\omega,t)$. These terms will be called \emph{interference terms}, and they are mainly responsible for the different types of TFR behavior, considered below. Note that, as seen from (\ref{OVamp}), changing the phase lag $\varphi_{1,2}$ in (\ref{OVs}) does not lead to any qualitative or quantitative changes in TFR behavior, but only shifts it in time.

\subsubsection{Representation}\label{sec:OVA}

\begin{figure*}[t!]
\includegraphics[width=1.0\linewidth]{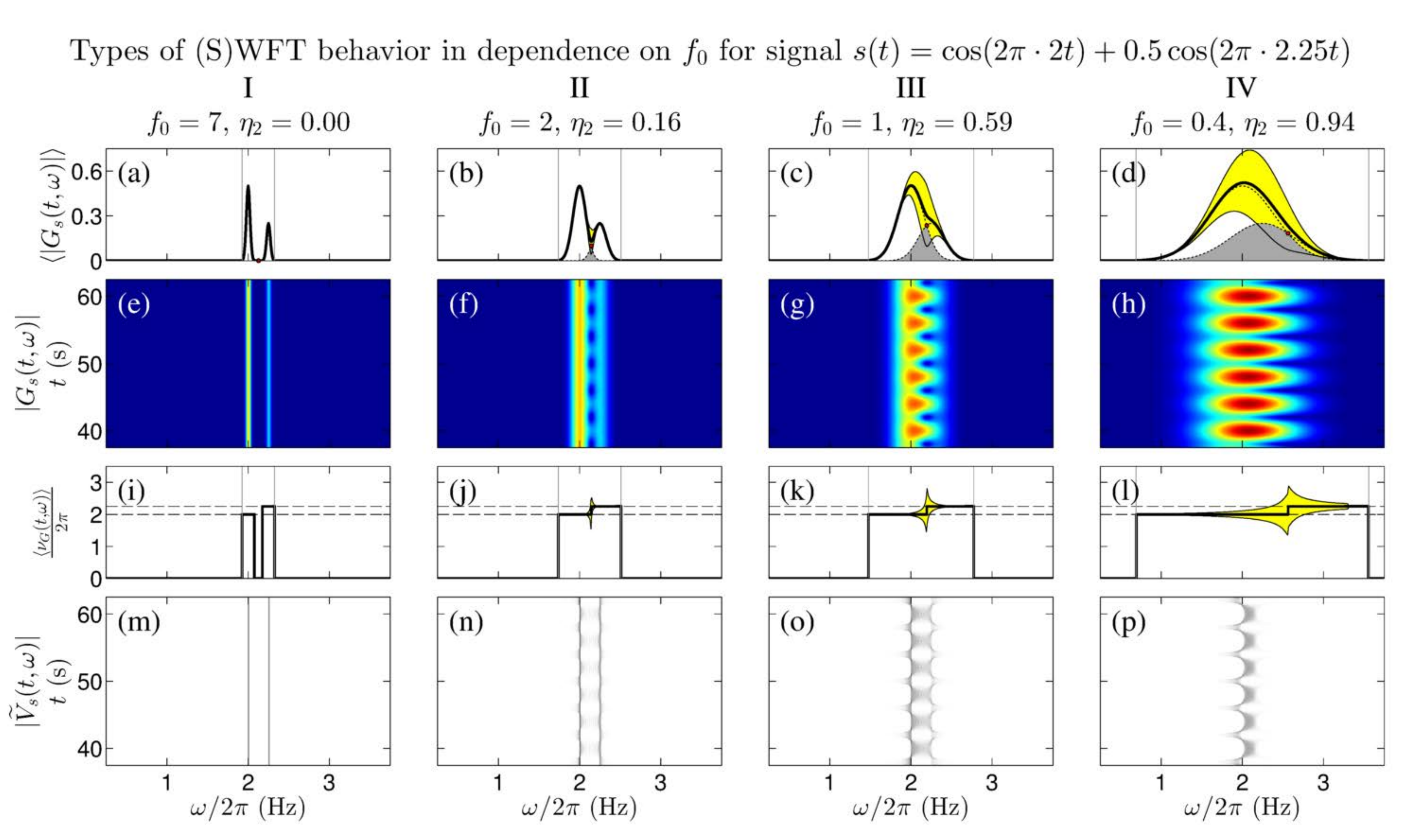}\\
\caption{Behavior of the Gaussian window WFT in dependence on $f_0$ for a two-tone signal $s(t)=\cos(2\pi\cdot2 t)+0.5\cos(2\pi\cdot2.25 t)$, sampled at $20$ Hz for $500$ s. For illustrational purposes, Gaussian window $\hat{g}(\xi)$ is ``cutted'' to compact frequency support $[\xi_1(0.001),\xi_2(0.001)]$, and the boundaries of the joint support of both peaks $[\nu_1+\xi_1(0.001),\nu_2+\xi_2(0.001)]$ are shown by gray lines in (a-d) and (i-l). (a-d): Time-averaged WFT amplitudes; dotted lines show $\frac{1}{2}\hat{g}(\omega-\nu_1)$ and $\frac{1}{2}\hat{g}(\omega-\nu_2)$, corresponding to the WFTs of each tone separately, with red dot indicating point of their intersections and gray region showing the area shared by both of them. (e-h): WFT amplitudes in time-frequency domain. (i-l): Time-averaged WFT frequency $\nu_{G}(\omega,t)$, with dashed lines showing the frequencies of each tone $\nu_1,\nu_2$. (m-p): SWFT amplitudes in time-frequency domain. Values of $\eta_2$ indicate relative overlap for the second component (ratio of gray-shaded area to all area below the right dotted peak in (a-d), see text and (\ref{OVeta})), which in our case of $r=0.5<1$ is the maximum among two $\eta_2=r^{-1}\eta_1=\max[\eta_1,\eta_2]$. Note, that (d,h,l,p) correspond to Regime IV only if we take $\epsilon>0.06$ in its condition (see (\ref{OVt4}) below), for lower precision smaller $f_0$ are needed.}
\label{fig:OVbt}
\end{figure*}

Different types of TFR behavior (see Table \ref{BTtab}) for the two-tone signal (\ref{OVs}) are illustrated in Fig.\ \ref{fig:OVbt} for the Gaussian window (S)WFT; for other windows, as well as for the (S)WT, all is qualitatively the same. As can be seen, for sufficiently high $f_0$, the WFT amplitude has two well-separated peaks at all times (Regime I). For lower $f_0$, although there are still two distinct peaks in WFT at all times, ``bridges'' begin to appear between them at certain times, reflecting interference between components and causing localised corruption of the WFT (Regime II). Regime III behaviour appears when we further decrease $f_0$, so the window/wavelet frequency resolution becomes insufficient to resolve the two tones, and they become mixed, i.e.\ sometimes there are two peaks in WFT amplitude and sometimes only one. Finally, for a very low value of the resolution parameter, the frequency resolution becomes so poor (although time resolution is extremely sharp) that two tones are completely merged, appearing as a single AM/FM component in the TFR (Regime IV).

Comparing (a-d) and (m-p) in Fig.\ \ref{fig:OVbt}, it can be seen that, in agreement with what was said above, synchrosqueezing does not change the qualitative behavior of the TFR. Thus, when interference is present, it affects both the TFR amplitude and instantaneous frequency (\ref{OVamp}) (as shown in the first and third row in Figure \ref{fig:OVbt}), both of which are used in the SWFT/SWT construction. Note that, for an SWT based on the wavelet with compact frequency support, the resolution of two tones was also considered in \cite{Wu:11}, but here we will give a more general and detailed treatment for both the (S)WFT and (S)WT.

In Fig.\ \ref{fig:OVbt}(a-d), the dotted lines indicate the ``interference-free'' WFT amplitudes, i.e.\ as they would be if the signal consisted only of a single tone, with the gray-shaded area showing the region shared by each of the ``non-interfering'' peaks. These can be associated with the two terms in $H_s(\omega,t)$ (\ref{OVamp}): at each frequency $\omega$, the upper of the two dotted lines represents the amplitude of the currently dominant component, while the lower one reflects the contribution of the other tone and equals the amplitude of the oscillations at its frequency. It seems logical, therefore, to measure the severity of interference in terms of the relative overlaps $\eta_{1,2}$, defined as the ratio of the interference area (gray-shaded in Fig.\ \ref{fig:OVbt}) to the total area under each peak. Using (\ref{gformNHQ}), one can write it as
\begin{equation}\label{OVeta}
\begin{aligned}
\eta_1=r\eta_2
=&\frac{\int \min[\hat{h}_{\nu_1}(\omega),\hat{h}_{\nu_2}(\omega)]d\mu(\omega)}{\int \hat{h}_{\nu_1}(\omega)d\mu(\omega)}\\
=&\frac{\int_{-\infty}^{\mu(\omega_\times)} \hat{h}_{\nu_1}(\omega)d\mu(\omega)+\int_{\mu(\omega_\times)}^\infty \hat{h}_{\nu_1}(\omega)d\mu(\omega)}
{\int \hat{h}_{\nu_1}(\omega)d\mu(\omega)}\\
=&Q_{\nu_1}(\omega_\times)+r[1-Q_{\nu_2}(\omega_\times)],
\end{aligned}
\end{equation}
where $Q_{\nu}(\omega)$ is defined in (\ref{gformNHQ}), and the ``intersection frequency'' $\omega_\times$ is determined as the frequency where the WFT/WT amplitudes of separate (non-interfering) tones intersect each other (red point in Figure \ref{fig:OVbt}), i.e.\ \begin{equation}\label{omisect}
\begin{gathered}
\omega_\times\in[\nu_1,\nu_2]:\;
\hat{h}_{\nu_1}(\omega_\times)=r\hat{h}_{\nu_2}(\omega_{\times})\\
\begin{aligned}
&\Rightarrow \omega_{\times}=\bar{\nu}-\frac{\log r}{f_0^2\Delta\nu}
\quad\mbox{(Gaussian window WFT)},\\
&\Rightarrow \omega_{\times}=\sqrt{\nu_1\nu_2}\exp\Big[-\frac{\log r}{(2\pi f_0)^2\log\frac{\nu_2}{\nu_1}}\Big]
\quad\mbox{(lognormal wavelet WT)},
\end{aligned}
\end{gathered}
\end{equation}
The relative overlaps $\eta_{1,2}$ provide rough measures of the error that one can expect while extracting and reconstructing by a direct method the first and second components from the TFR of a summed signal (\ref{OVs}); a more rigorous relation will be considered in the reconstruction subsection below. Interestingly, the jump in $\langle \nu_G(\omega,t)\rangle$ seen in Fig.\ \ref{fig:OVbt} (j-l) occurs exactly at $\omega=\omega_{\times}$, and the same situation is observed for the WT. Note also that, in the case of $r=1$ (tones of equal amplitude), for WFT with symmetric $\hat{g}(\xi)$ one always has $\omega_\times=\bar{\nu}$, while for the WT with $\hat{\psi}(\xi)$ symmetric over $\log\xi$ it is always $\omega_\times=\sqrt{\nu_1\nu_2}$.

\begin{center}
\begin{table*}[t!]
  \begin{subequations}
  \begin{tabular}{| >{\centering\arraybackslash}m{1.5cm} | m{16cm} |}
  \hline
Regime & \multicolumn{1}{c|}{Condition}
\\ \hline
\textbf{I} &
\begin{equation}\label{OVt1}
\max[\eta_1,\eta_2]\leq\epsilon
\end{equation}
\begin{equation}\label{OVt1a}
\overset{approx.}{\Rightarrow}
\left[\begin{array}{rl}
\mbox{(S)WFT:} & \Delta \nu>\xi_2(\epsilon)-\xi_1(\epsilon),\\
\mbox{(S)WT:} & \nu_2/\nu_1>\xi_2(\epsilon)/\xi_1(\epsilon).\\
\end{array}\right.
\end{equation}
\\ \hline
\textbf{II} &
\begin{equation}\label{OVt2}
\left\{\begin{array}{l}
\max[\eta_1,\eta_2]>\epsilon\\
\langle N_p\rangle=2\Rightarrow
[\hat{h}_{\nu_1}(\omega)+r \hat{h}_{\nu_2}(\omega)]\mbox{ has minimum in }\omega\in[\nu_1,\nu_2]\\
\end{array}\right.
\end{equation}
\\ \hline
\textbf{III} &
\begin{equation}\label{OVt3}
1<\langle N_p\rangle<2\;\Rightarrow\;
\left\{\begin{array}{l}
\max[\eta_1,\eta_2]< 1-\epsilon\\
{[}\hat{h}_{\nu_1}(\omega)+r \hat{h}_{\nu_2}(\omega)]\mbox{ has no minimum in }\omega\in[\nu_1,\nu_2]\\
\end{array}\right.
\end{equation}
\\ \hline
\textbf{IV} &
\begin{equation}\label{OVt4}
\max[\eta_1,\eta_2]\geq 1-\epsilon
\end{equation}
\begin{equation}\label{OVt4a}
\overset{approx.}{\Rightarrow}
\left[\begin{array}{rl}
\mbox{(S)WFT:}&
\left[\begin{array}{l}
\omega_{\times}\geq\nu_2+\xi_2(2\epsilon)\\
\omega_{\times}\leq\nu_1+\xi_1(2\epsilon)
\end{array}\right.
\Rightarrow
r^{\mp 1}\geq\exp\left\{\frac{f_0\Delta\nu}{2}[2n_G(2\epsilon)+f_0\Delta\nu]\right\}
\mbox{ for Gaussian window}\\
 & \\
\mbox{(S)WT:}&
\left[\begin{array}{l}
\omega_{\times}\geq\omega_\psi\nu_2/\xi_1(2\epsilon)\\
\omega_{\times}\leq\omega_\psi\nu_1/\xi_2(2\epsilon)
\end{array}\right.
\Rightarrow
r^{\mp 1}\geq\exp\left\{\frac{(2\pi f_0)\log\frac{\nu_2}{\nu_1}}{2}\Big[2n_G(2\epsilon)+(2\pi f_0)\log\frac{\nu_2}{\nu_1}\Big]\right\}
\mbox{ for lognormal wavelet}\\
\end{array}\right.
\end{equation}
\\ \hline
  \end{tabular}
  \end{subequations}
\caption{Conditions for each type of behavior (illustrated in Fig.\ \ref{fig:OVbt}) for the two-tone signal (\ref{OVs}), where we have used notation (\ref{OVnt}). Value of $\epsilon$ is some predefined accuracy (we use $\epsilon=0.001$) that determines how high (low) the interference should be to regard the tones as fully merged (separated) in the TFR.}
\label{tab:OVc}
\end{table*}
\end{center}

Now, using (\ref{OVamp}) and (\ref{OVeta}), we can formulate the conditions for each type of behavior, which are summarized in Table \ref{tab:OVc}. The I type of behavior requires the interference to be very small which, rewritten in terms of relative overlaps, requires that $\eta_{1,2}$ be smaller than $\epsilon$ (\ref{OVt1}). This condition can be rewritten approximately in terms of the $\epsilon$-supports: it can be shown that the expressions (\ref{OVt1}) and \ref{OVt1a} are fully equivalent if either: (a) $r=1$, i.e.\ tones have equal amplitude, and the window (wavelet) FT is symmetric over $\omega$ ($\log \frac{\omega}{\omega_\psi}$); (b) $\epsilon=0$ and the window/wavelet has a compact frequency support $|\mu(\xi_{1,2}(0))|<\infty$. The latter case is especially simple, since tones with any $r$ are obviously separated if $\hat{h}_{\nu_{1,2}}(\omega)$ has non-overlapping finite frequency supports.

However, in other cases (\ref{OVt1a}) is only approximate, depending on the amplitude ratio $r$ and the (time-dependent) positions of the minima between the two TFR amplitude peaks. Nevertheless, it can be shown that (\ref{OVt1a}) implies $\eta_1\leq(1+r)\epsilon/2,\;\eta_2\leq(1+r^{-1})\epsilon/2$, so that if tones are non-negligible compared to each another, one will always have $\eta_{1,2}\leq O(\epsilon)$. Therefore, the considered approximation appears to be very accurate (especially for the Gaussian window WFT, see Fig.\ \ref{OVbdWFT} below).

The conditions for the II and III Regimes, (\ref{OVt2}) and (\ref{OVt3}), are devised by investigating the minima of $|H_s(\omega,t)|^2$ (\ref{OVamp}). Thus, since the window/wavelet FT is assumed to be positive and unimodal, it can be shown that, in terms of the number of peaks, the tones are most merged and most resolved in the TFR when in (\ref{OVamp}) $\Delta\phi(t)=0$ and $\Delta\phi(t)=\pi$, respectively; the TFR amplitudes (\ref{OVamp}) for these limiting cases are:
\begin{equation}\label{OVms}
\left.|H_s(\omega,t)|^2\right|_{\Delta\phi(t)=0,\pi}=\frac{A^2}{4}\Big[\hat{h}_{\nu_1}(\omega)\pm r\hat{h}_{\nu_2}(\omega)\Big]^2
\end{equation}
Obviously, if for the most-merged case $\Delta\phi(t)=0$ one still has two peaks in the $|H_s(\omega,t)|^2$, then there will always be two peaks in the TFR amplitude. Similarly, if only one peak appears in the most-resolved case ($\Delta\phi(t)=\pi$), then the TFR amplitude will always have a single peak. Taken with (\ref{OVamp}), this logic gives (\ref{OVt2}) and (\ref{OVt3}).
It should be clarified that, in (\ref{OVt3}), this is the condition $\max[\eta_1,\eta_2]<1-\epsilon$ which establishes that the two tones are not \emph{always} merged into one TFR peak (i.e.\ $\langle N_p\rangle>1$). Thus, comparing the equation for $\omega_\times$ (\ref{omisect}) with (\ref{OVamp}), one can see that $|H_s(\omega,t)|^2$ for $\Delta\phi(t)=\pi$ drops to zero (i.e.\ has a minimum, implying multiple peaks) at $\omega=\omega_\times$. Therefore, the necessary condition for $\langle N_p\rangle=1$ is non-existence of $\omega_\times$ which, according to (\ref{OVeta}), is equivalent to $\max[\eta_1,\eta_2]=1$, being violated by $\max[\eta_1,\eta_2]<1-\epsilon$ in (\ref{OVt3}).

Finally, Regime IV appears when almost all of the area under the TFR amplitude corresponding to one tone is included into the area of the other, as reflected in (\ref{OVt4}). When this happens, the interference becomes extremely strong, and the tones behave as a single component in the TFR nearly all the time. In fact, one may require tones to be \emph{always} merged and set $\epsilon=0$ in (\ref{OVt4}); this is a stricter condition and is equivalent to the non-existence of $\omega_\times$ (which is necessary for $\langle N_p\rangle=1$, as discussed previously). However, for many windows/wavelets such situation is in principle impossible, e.g.\ if $\hat{h}_{\nu_1}(\omega)$ and $\hat{h}_{\nu_2}(\omega)$ have finite supports, there will always be an intersection between the two (since the support of the latter will end after the support of the former). The definition through $\epsilon$ is therefore much more general and convenient. With decreasing $\epsilon$, the time intervals for which more than one peak exists are ``squeezed'' around the times when $\Delta\phi(t)=\pi$, becoming very short in Regime IV as defined by small $\epsilon$ (in fact, the Regime IV condition is equivalent to $\langle N_p\rangle\leq 1+O(\epsilon)$). The approximate conditions (\ref{OVt4a}) imply that $\max[\eta_1,\eta_2]\geq1-\alpha\epsilon,\alpha\leq1$, thus being slightly stricter and implying (but not being implied by) (\ref{OVt4}).

By performing numerical simulation, one can locate the regions of $\{r,\Delta\nu,f_0\}$ space corresponding to each type of behavior. They are shown for a Gaussian window WFT in Fig.\ \ref{OVbdWFT} together with $\langle N_p\rangle$ (\ref{NP}) and maximum relative overlap $\max[\eta_1,\eta_2]$ (\ref{OVeta}). Both qualitatively and quantitatively, the WFT behavior, apart from $r$, depends only on the product $f_0\Delta\nu$, so that the pictures in Fig.\ \ref{OVbdWFT} remain the same for all $f_0$. Additionally, for symmetric $\hat{g}(\xi)$ (such as Gaussian), all the pictures are also symmetric with respect to $r\rightarrow1/r$, which in this case is equivalent to exchanging $\eta_1\leftrightarrow\eta_2$. Note that the border of Regime I predicted by (\ref{OVt1}) (blue line in Fig.\ \ref{OVbdWFT} (a)) is almost independant of $r$ and lies very close to the predictions of (\ref{OVt1a}), shown by the dashed light-blue line.

\begin{figure*}[t!]
\includegraphics[width=1.0\linewidth]{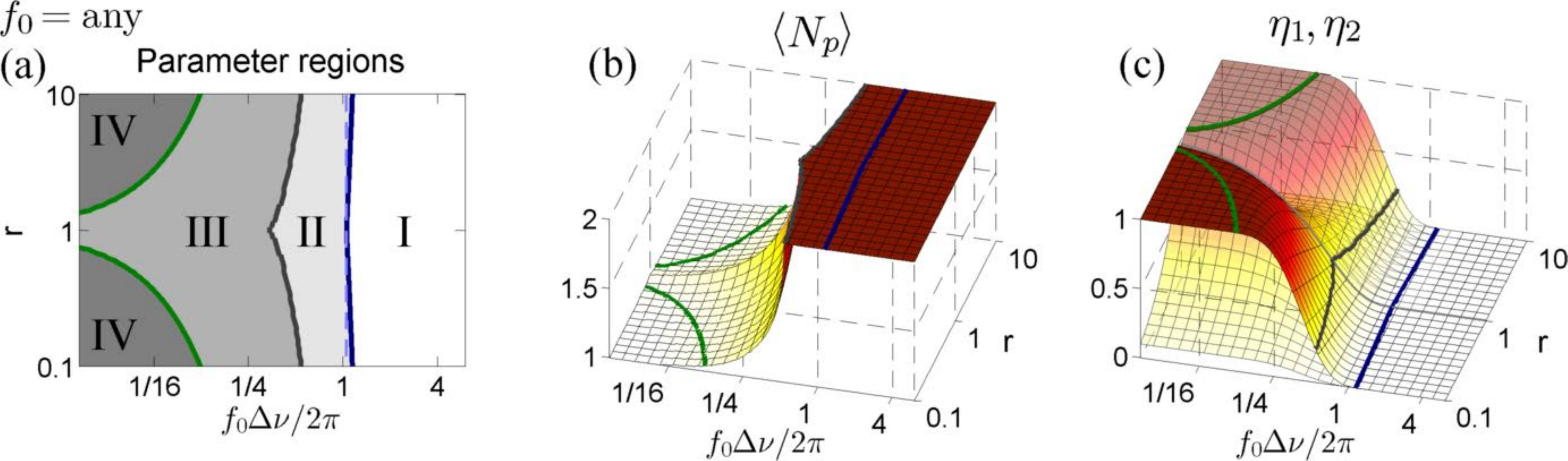}\\
\caption{Dependence of the WFT behavior on signal (\ref{OVs}) parameters $r,\nu_1,\nu_2\equiv\nu_1+\Delta\nu$ and Gaussian window resolution parameter $f_0$. (a): Regions of parameter space corresponding to each type of behavior, according to (\ref{OVt1}),(\ref{OVt2}),(\ref{OVt3}),(\ref{OVt4}); dashed light-blue line shows boundary of the I-type behavior as predicted by approximate (\ref{OVt1a}). (b): Mean number of peaks $\langle N_p\rangle$ (\ref{NP}). (c): Relative overlaps $\eta_1$ (transparent) and $\eta_2$ (opaque), as calculated from (\ref{OVeta}); light-grey line shows their intersection. For determining Regimes I and IV we used $\epsilon=0.001$ in (\ref{OVt1}) and (\ref{OVt4}).}
\label{OVbdWFT}
\end{figure*}

\begin{figure*}[t!]
\includegraphics[width=1.0\linewidth]{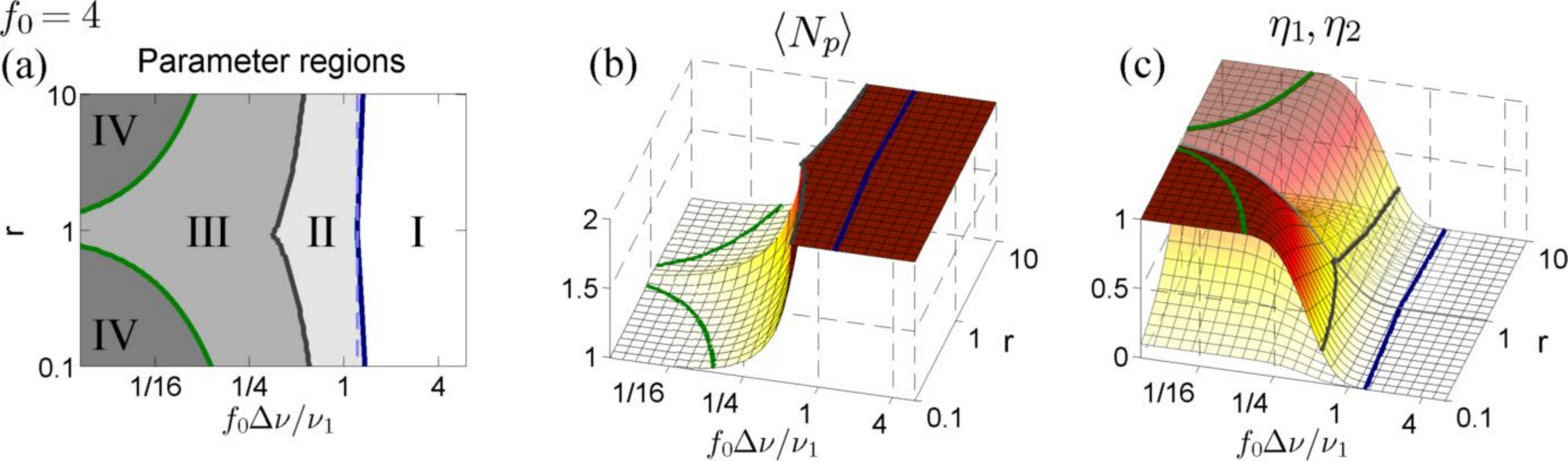}\\
\includegraphics[width=1.0\linewidth]{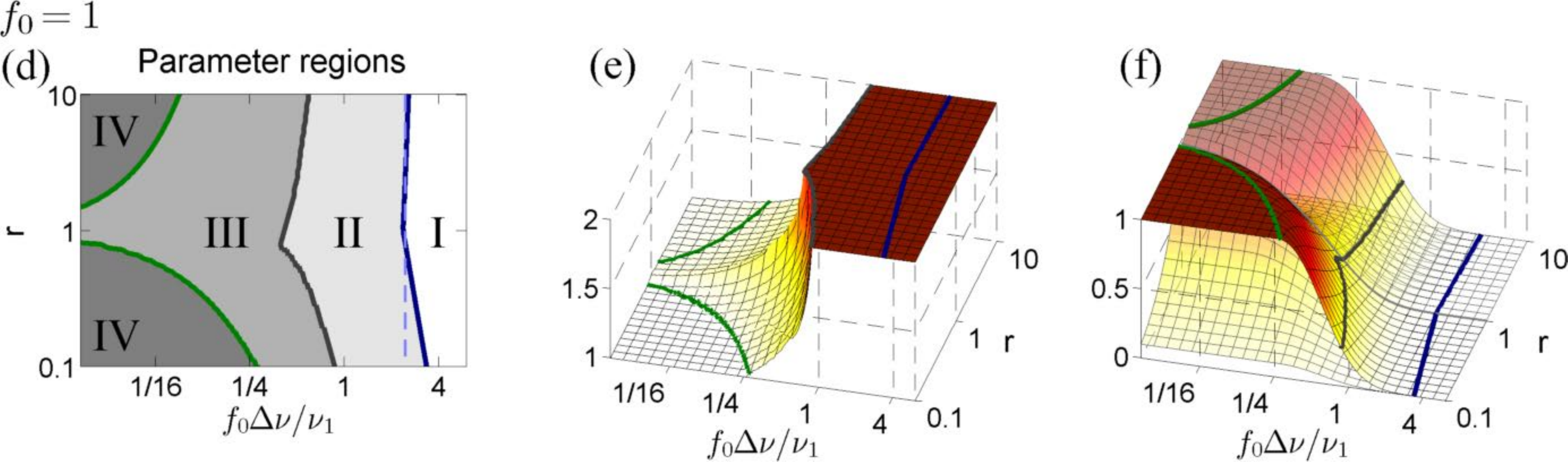}\\
\caption{Same as Figure \ref{OVbdWFT}, but for the Morlet wavelet WT.}
\label{OVbdWT}
\end{figure*}

\begin{figure*}[t!]
\includegraphics[width=1.0\linewidth]{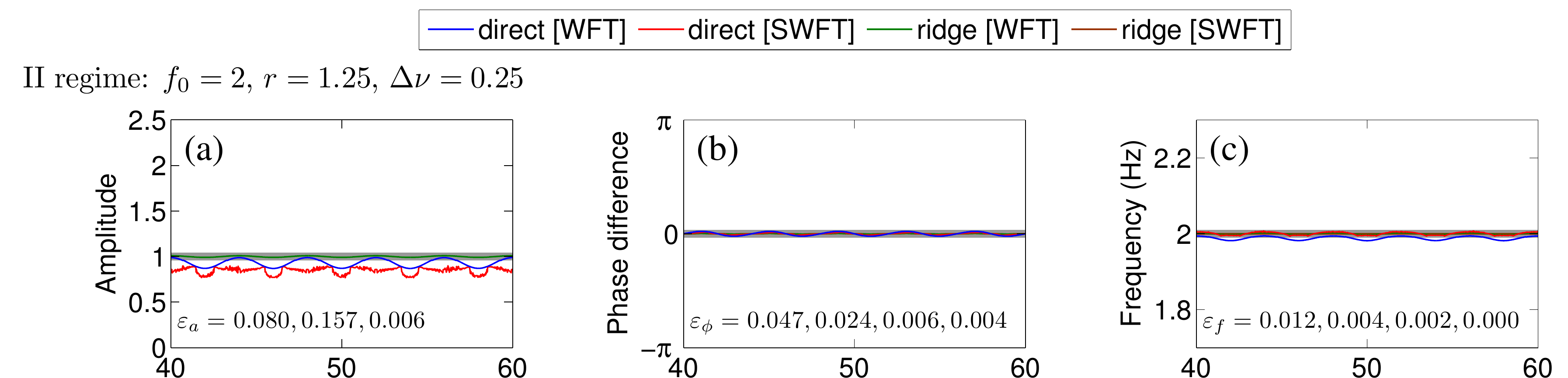}\\
\includegraphics[width=1.0\linewidth]{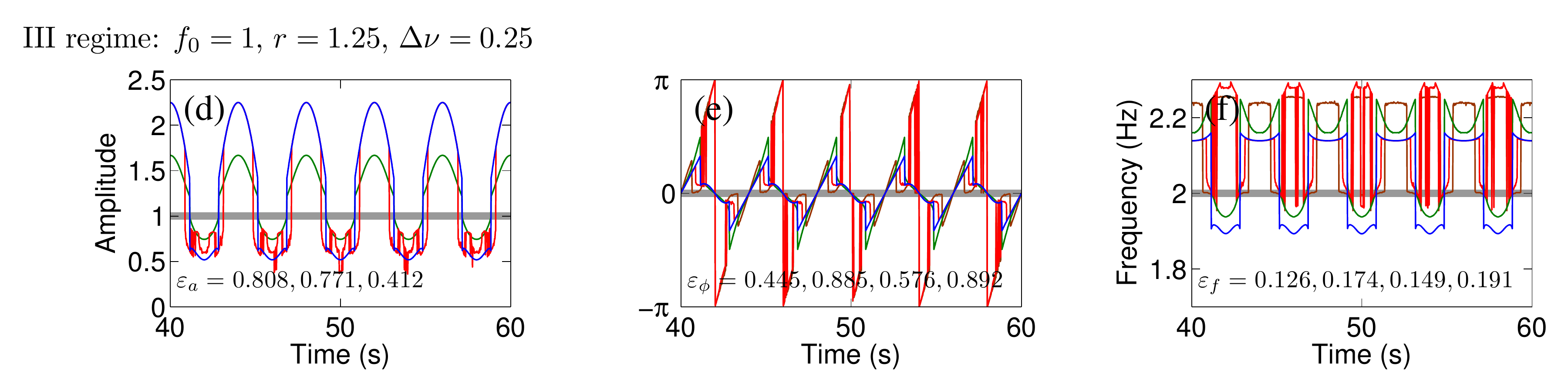}\\
\caption{Amplitude (a,d), phase (b,e) and frequency (c,f) of a two-tone signal (\ref{OVs}) as reconstructed from its WFT and SWFT (colored lines), compared to the true values (thick gray lines). Values of $\varepsilon_{a,\phi,f}$ shown are in the same order as lines in legend, corresponding to ${\rm direct [WFT]}$ (blue), ${\rm direct [SWFT]}$ (red), ${\rm ridge [WFT]}$ (green), ${\rm ridge [SWFT]}$ (brown). In (a,d), ridge reconstruction from the SWFT is not shown as it is not appropriate for amplitude (see Part I). In (b,e), the difference between the reconstructed and true phase is shown. The signal (\ref{OVs}) was sampled at $50$ Hz for $100$ s, and it was simulated with $\varphi_1=\varphi_2=0$ and $\nu_1/2\pi=2$; all other parameters are indicated on the figure.}
\label{OVex}
\end{figure*}

Fig.\ \ref{OVbdWT} is the analog of Fig.\ \ref{OVbdWFT} for the Morlet wavelet WT. Due to its logarithmic frequency resolution one might expect that, in analogy to the WFT, the behavior of the WT is determined only by $f_0\log\nu_2/\nu_1\approx f_0\Delta\nu/\nu_1+O(\Delta\nu^2/\nu_1^2)$ and $r$. However, although this is indeed the case for different $\nu_1$, it does not hold true for different $f_0$. Thus, as seen from a comparison of Fig.\ \ref{OVbdWT} (a-c) and (d-f), the behavior of the WT, apart from $r$ and $f_0\Delta\nu/\nu_1$ (or any other simple combination of $f_0,\omega_\psi,\nu_{1,2}$), is characterized also by $f_0$ separately. This arises because the form of the Morlet wavelet changes qualitatively with $f_0$ (mainly due to the admissibility term $e^{-(2\pi f_0)^2/2}$); in the case of the lognormal wavelet, for which $\hat{\psi}(\omega_\psi/\omega)$ does not change its form with $f_0$, the WT's behavior is characterized only by $r$ and $f_0\log\nu_2/\nu_1$, as expected (not shown). For the same reasons, $\langle N_p\rangle$ and $\max(\eta_1,\eta_2)$, although symmetric under $r\rightarrow r^{-1}$ for Gaussian window WFT (Fig.\ \ref{OVbdWFT}) and lognormal wavelet WT (not shown), are no longer symmetric for the Morlet wavelet. Thus, as seen from Fig.\ \ref{OVbdWT}, it is easier to separate tones perfectly (Regime I) when the lower-frequency tone has a higher amplitude than the higher-frequency one ($r<1$). As a result of this asymmetry, the exact and approximate conditions for Regime I type behavior, (\ref{OVt1}) and (\ref{OVt1a}), exhibit poorer agreement in the case of WT as compared to the WFT (and better agreement for $r>1$ than for for $r<1$). Note, however, that with increase of $f_0$ in the Morlet wavelet the WT behavior becomes more and more symmetric under $r\rightarrow r^{-1}$ (compare Fig.\ 1(d-f) and (a-c)), as well as less dependent on $f_0$ separately from $f_0\Delta\nu/\nu_1$.

\subsubsection{Reconstruction}

To relate the relative overlaps $\eta_{1,2}$ (\ref{OVeta}) to the errors of parameter estimation, consider the reconstruction of each tone by the direct method (see Part I). We denote the point in the WFT/WT at which the supports of each tone are separated as $\omega_s(t)$, so that at each time the first tone $s_1(t)=A\cos(\nu_1 t+\varphi_1)$ is reconstructed from the frequency region $[\omega_-^{(1)}(t),\omega_+^{(1)}(t)]=(-\infty,\omega_s(t)]$, while the second one $s_2(t)=rA\cos(\nu_2 t+\varphi_2)$ -- from $[\omega_s(t),\infty)$. The reconstructed signals $s_{1,2}^{rec}(t)$ are then
\begin{equation}\label{OVrec}
\begin{aligned}
s_1^{rec}(t)-A\cos(\nu_1 t+\varphi_1)=&A\Big[-\Big(1-Q_{\nu_1}\big(\omega_s(t)\big)\Big)\cos(\nu_1 t+\varphi_1)\\
&+rQ_{\nu_2}\big(\omega_s(t)\big)\cos(\nu_2 t+\varphi_2)\Big]\\
s_2^{rec}(t)-rA\cos(\nu_1 t+\varphi_1)=&A\Big[\Big(1-Q_{\nu_1}\big(\omega_s(t)\big)\Big)\cos(\nu_1 t+\varphi_1)\\
&-rQ_{\nu_2}\big(\omega_s(t)\big)\cos(\nu_2 t+\varphi_2)\Big]\\
\end{aligned}
\end{equation}
According to the definition of the TFS (see Part I), the frequency $\omega_s(t)$ corresponds to a minimum of $|H_s(\omega,t)|^2$ at each time, which will obviously lie between the tone frequencies $\omega_s(t)\in(\nu_1,\nu_2)$. As discussed in Part I, for $r=1$ and symmetric $\hat{g}(\xi)$ (WFT) or logarithmically symmetric $\hat{\psi}(\xi)$, one has $\omega_s(t)=[(\nu_1+\nu_2)/2\mbox{ (WFT) or }\sqrt{\nu_1\nu_2}\mbox{ (WT)}]=\omega_\times$, where the intersection frequency $\omega_\times$ is defined in (\ref{omisect}). In other cases it will be not so, but one can still expect $\langle \omega_s(t)\rangle\approx\omega_\times$. Therefore, setting $\omega_s(t)\approx\omega_\times$ in (\ref{OVrec}), one obtains
\begin{equation}\label{maxerreta}
\begin{gathered}
\max\big(A^{-1}|s_1^{rec}(t)-s_1(t)|\big)\approx \eta_1,\\
\max\big(r^{-1}A^{-1}|s_2^{rec}(t)-s_2(t)|\big)\approx \eta_2,
\end{gathered}
\end{equation}
which shows that the relative overlaps $\eta_{1,2}$ (\ref{OVeta}) indeed have a direct relationship with the quality of the tones' representation in the TFR.

We now investigate numerically the performance of the different reconstruction methods (direct and ridge-based) for the two-tone signal (\ref{OVs}). To extract the time-frequency supports of the two tones, at each time we find the two most dominant peaks in the TFR amplitude, pick the one nearest to the actual tone frequency, and select the corresponding TFS around it (``frequency-based'' scheme). This is done for each tone separately so, when they merge into a single peak, they will have the same extracted support. Such procedure give the most appropriate TFS for reliable study of the current case. We then apply the direct and ridge reconstruction methods to obtain amplitude, phase and frequency and calculate the respective errors $\varepsilon_{a,\phi,f}$ (\ref{apfrel}).

Figure \ref{OVex} shows examples of amplitudes, phases and frequencies reconstructed from the (S)WFT in two cases, corresponding to Regimes II and III (for (S)WT all remains qualitatively similar). In Regime II (a-c), when interference is present but not very strong, ridge methods by far outperform direct ones. The performance of ridge reconstruction from the WFT and SWFT is almost the same (except for the amplitudes), which cannot be said about the direct methods. Thus, for amplitude reconstruction, ${\rm direct [WFT]}$ performs better than ${\rm direct [SWFT]}$, while for phase/frequency we observe the opposite situation, i.e.\ direct reconstruction from the SWFT gives slightly better results than that from the WFT (although in all cases direct estimates are less accurate than ridge-based). All changes when the WFT enters Regime III (Fig.\ \ref{OVex}(d-f)). In this case, where the two tones are often merged into a single peak, all methods fail.

\begin{figure*}[t!]
\includegraphics[width=1.0\linewidth]{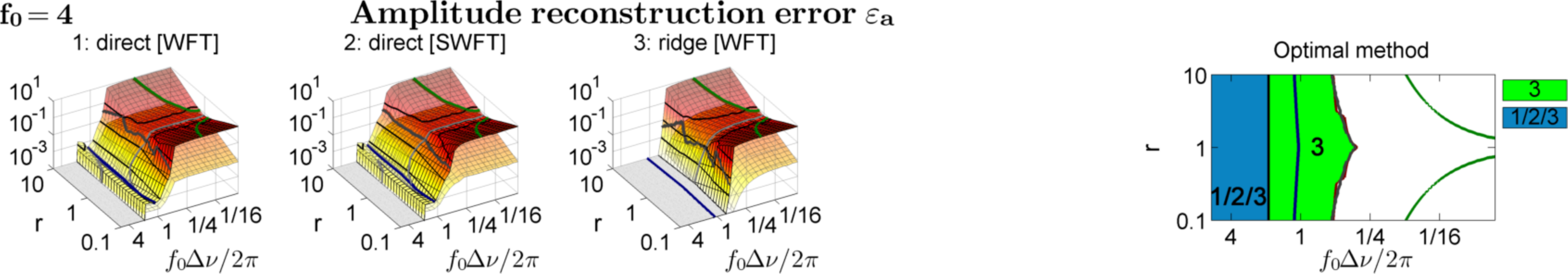}\\
\includegraphics[width=1.0\linewidth]{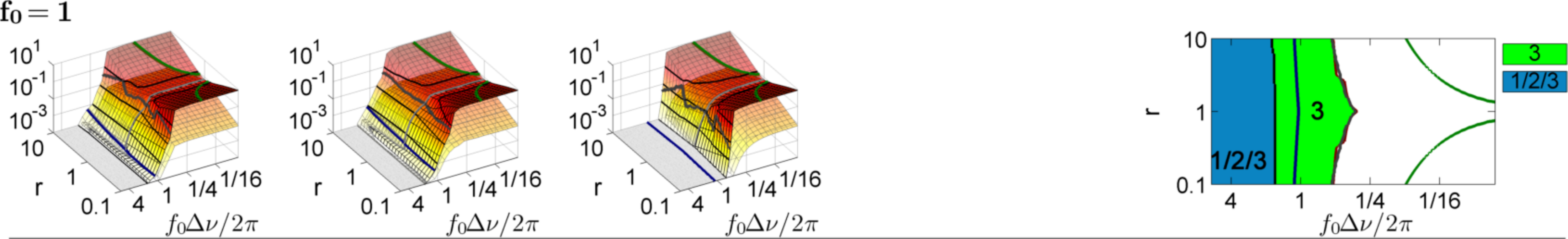}\\
\includegraphics[width=1.0\linewidth]{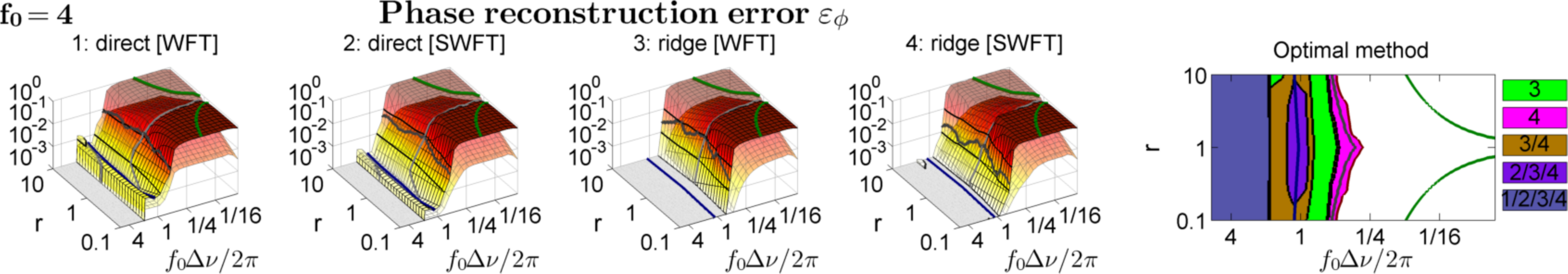}\\
\includegraphics[width=1.0\linewidth]{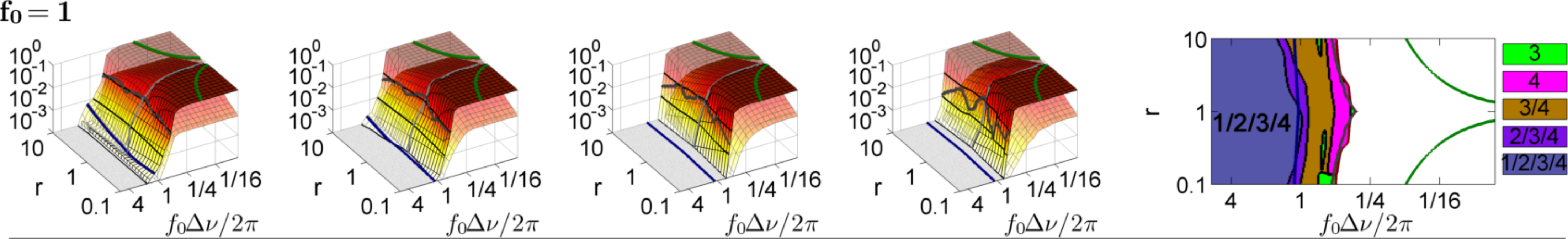}\\
\includegraphics[width=1.0\linewidth]{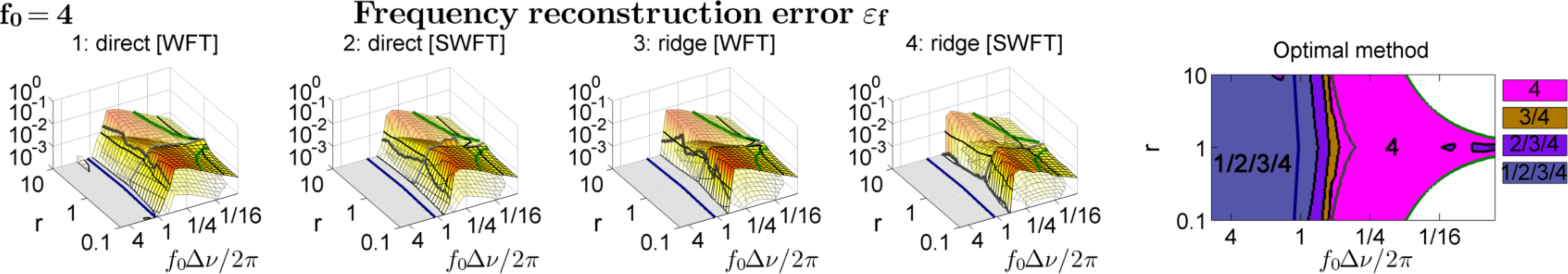}\\
\includegraphics[width=1.0\linewidth]{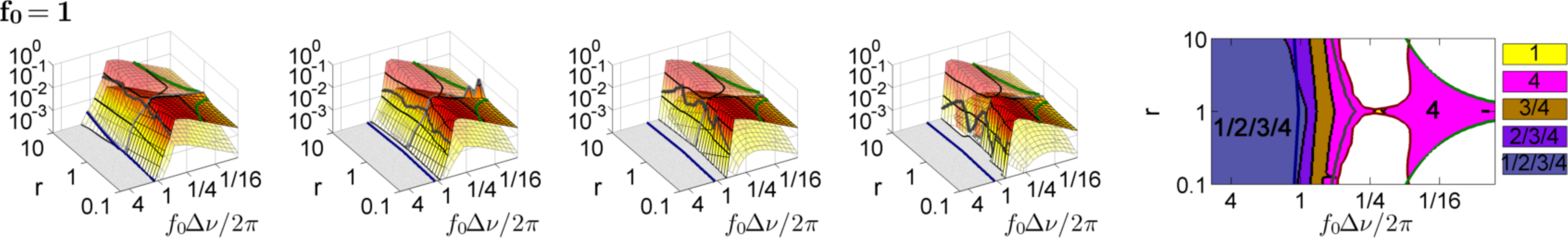}\\
\caption{Left four panels in each row show the (S)WFT-based amplitude (1-2 rows), phase (3-4 rows) and frequency (5-6 rows) reconstruction errors for the two interfering tones $s(t)=\cos \nu_1 t+r\cos(\nu_1+\Delta\nu)t$ in dependence on parameters $r,\Delta\nu$ and Gaussian window resolution parameter $f_0$. Errors for the reconstruction of the first component (of frequency $\nu_1$) are shown transparent, and for the second one (of frequency $\nu_1+\Delta\nu$) are opaque. We used $\nu_1/2\pi=1$, but results do not depend on its value. Estimation of amplitude from the SWFT ridges is not appropriate, and so the corresponding results are discarded. Thin black lines indicate the levels of $0.001,0.01,0.1,1$, and the behavior of errors $<0.001$ for simplicity is not shown; thick green, gray and blue lines are the same as in Fig. \ref{OVbdWFT}, showing the borders of regions corresponding to behavior in Regimes I - IV.
The right panels in each row show regions in parameter space where each method is optimal, i.e.\ gives the smallest estimation error; if the resultant errors of two methods differ on less than $0.001$, they are regarded as having similar performance (the corresponding regions are denoted as `'method1/method2`'). Comparison does not make much sense when WFT behavior is of IV type or if all reconstruction errors are $>0.1$, so that all methods fail; the corresponding regions are white-filled in the right-hand panels. The signal was sampled at 50 Hz for 100 s.}
\label{OVrecWFT}
\end{figure*}

\begin{figure*}[t!]
\includegraphics[width=1.0\linewidth]{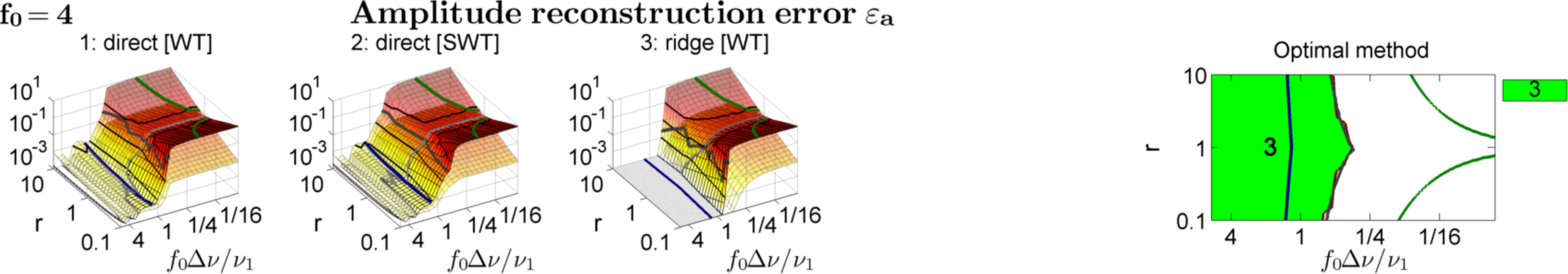}\\
\includegraphics[width=1.0\linewidth]{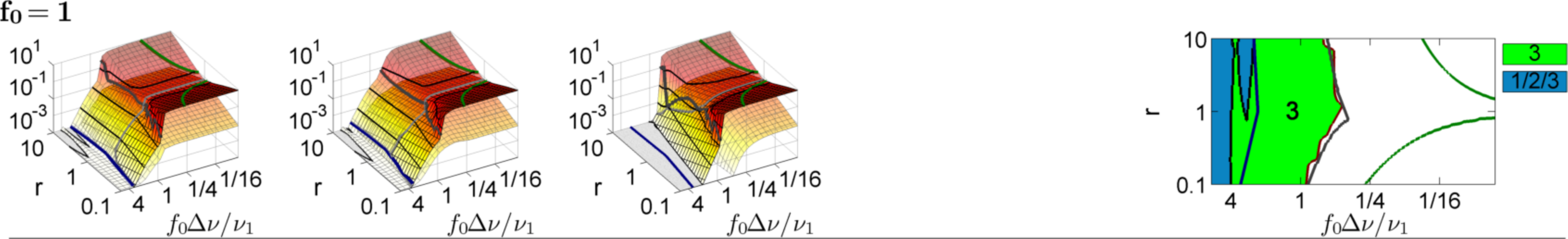}\\
\includegraphics[width=1.0\linewidth]{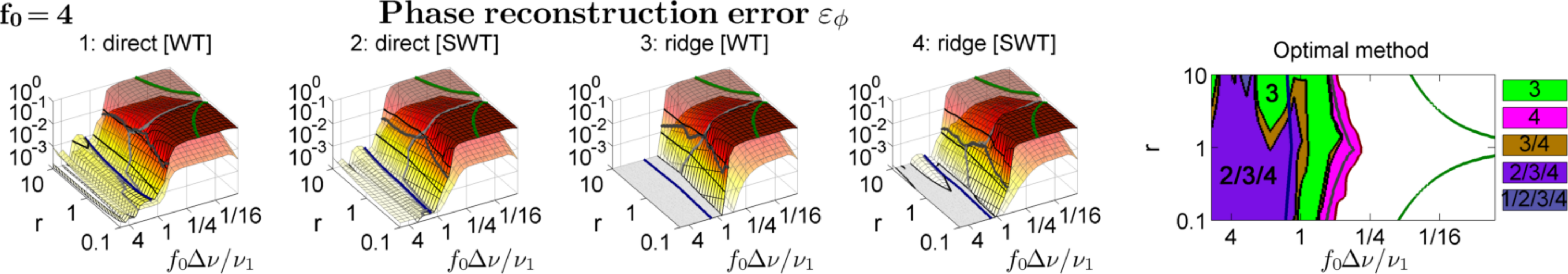}\\
\includegraphics[width=1.0\linewidth]{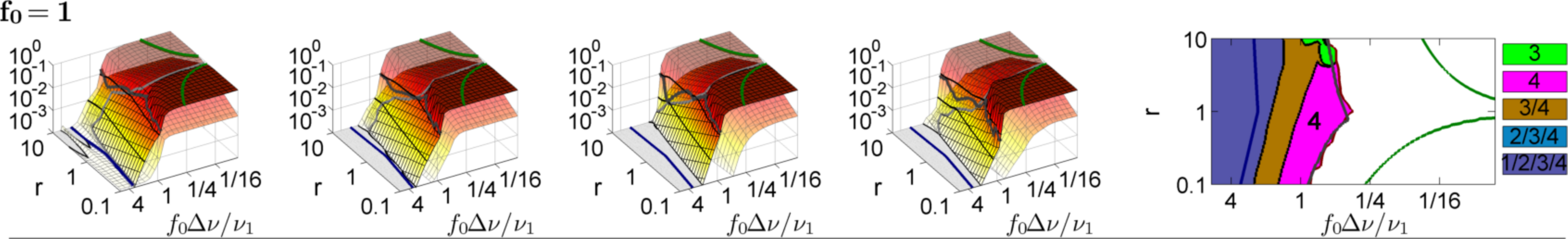}\\
\includegraphics[width=1.0\linewidth]{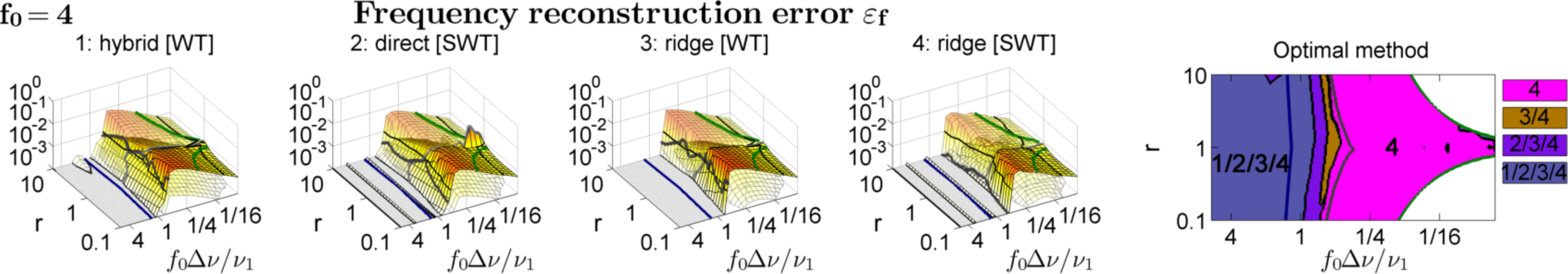}\\
\includegraphics[width=1.0\linewidth]{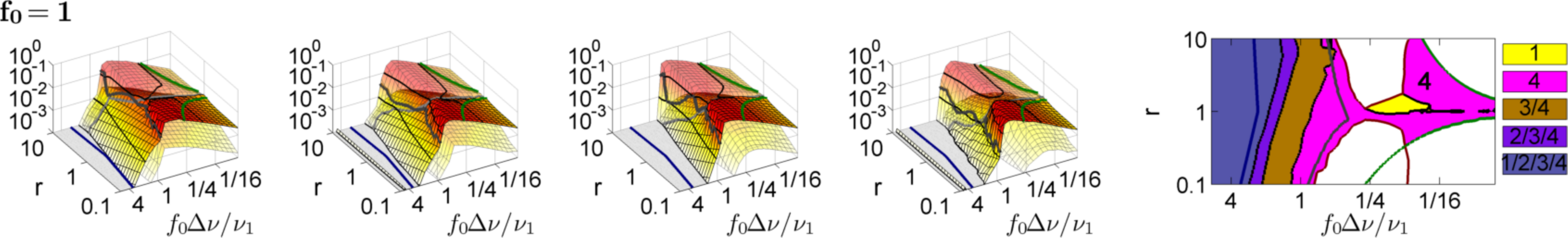}\\
\caption{Same as Fig.\ \ref{OVrecWFT}, but all parameters are reconstructed from the signal's WT/SWT. The direct frequency reconstruction for Morlet wavelet is not possible, so we use hybrid reconstruction (see Part I). Note, that the frequency reconstruction error for the WT (\ref{apfrel}) is defined in relation to the component's mean frequency, so $\varepsilon_f$ here is divided on $\nu_1/2\pi$ ($(\nu_1+\Delta\nu)/2\pi$) for the first (second) tone in comparison to $\varepsilon_f$ in Fig.\ \ref{OVrecWFT}; this makes the depicted dependences consistent for different $\nu_1$.}
\label{OVrecWT}
\end{figure*}

Figure \ref{OVrecWFT} shows the full dependence of all errors on the $f_0$ and signal parameters for (S)WFT-based reconstruction by different methods, while Figure \ref{OVrecWT} shows the same information for the WT. One can see, that reconstruction errors are well correlated with relative overlaps $\eta_{1,2}$ shown previously, and (not surprisingly) that the accuracy of (S)WFT-based reconstruction depends on a combination $f_0\Delta\nu$, while the performance of (S)WT-based methods depends on $f_0$ and $\Delta\nu/\nu_1$ separately, at least for the Morlet wavelet (for lognormal wavelet all will depend only on $r$ and $f_0\log\frac{\nu_2}{\nu_1}$). It is clear that the errors are to a large extent determined by the TFR behavior, being for all methods negligible in Regime I, appearing in II and becoming very large in III-IV (e.g.\ the maximum error among the two, for each tone, crosses level of 10\% almost exactly at the border of III region, indicating that 10\% relative error is an appropriate threshold). Note that, for the direct methods, even for some parameters corresponding to Regime I, there exists a non-negligible error: this is entirely on account of boundary effects (see Part I) and it decreases with increasing signal time-length. Although we pad the signal with original values to suppress boundary errors, the direct method estimates are much more susceptible to boundary effects than the ridge-based ones, thus requiring more padded values to achieve the specified precision than is given by the corresponding formula in Part I.

Both Figs.\ \ref{OVrecWFT} and \ref{OVrecWT} confirm what was already seen from Fig.\ \ref{OVex}. Thus, for the present case of interfering tones, ridge-based reconstruction outperforms the direct methods in all cases. This was to be expected, because a TFR at the amplitude peak should be less corrupted by interference with nearby components than a TFR contained in the wider component support. Thus, while in direct reconstruction one integrates over the whole TFS, including regions close to the other component (i.e.\ with considerable interference), ridge reconstruction accounts for interference only at the more distal peak. This is not so clear, however, for the case of phase/frequency reconstruction from SWFT/SWT ridges, but results indicate that there is a similar situation. Comparing estimates obtained from the WFT/WT and SWFT/SWT, one can see that the synchrosqueezing does not provide significant advantages in terms of the components' reconstruction: the accuracy of SWFT/SWT-based estimates is usually comparable or lower than that of WFT/WT-based ones, although in some parameter regions they are slightly better for the direct reconstruction of phase/frequency. As already discussed, the performance of each method depends strongly on the type of TFR behavior, so that an appropriate choice of window/wavelet parameters is therefore essential for accurate reconstruction.

\newpage

\newpage

\subsection{Amplitude modulation}\label{sec:AM}

We now consider the AM component with sinusoidal amplitude modulation:
\begin{equation}\label{AMs}
\begin{aligned}
s(t)=&A[1+r_a\cos(\nu_a t+\varphi_a)]\cos(\nu t+\varphi)\\
=&A\Big[\cos(\nu t+\varphi)+\frac{r_a}{2}\cos[(\nu-\nu_a)t+(\varphi-\varphi_a)]\\
&+\frac{r_a}{2}\cos[(\nu+\nu_a)t+(\varphi+\varphi_a)]\Big]\\
\end{aligned}
\end{equation}
where the amplitude definition implues $r_a\leq1,\;0\leq\nu_a<\nu$. For convenience, we also denote
\begin{equation}\label{AMnt}
\phi_a(t)\equiv\nu_a t+\varphi_a
\end{equation}

It is evident that the AM component (\ref{AMs}) can be represented as a sum of three tones: a main tone of frequency $\nu$ and two side tones at $\nu\pm\nu_a$ of equal amplitude; the latter appear as a result of the amplitude modulation and therefore will be called ``AM-induced''. Therefore, all the formulas and classification for multitone signals also apply for the AM component (\ref{AMs}). Using (\ref{gformGEN}), one obtains for the signal (\ref{AMs}):
\begin{equation}\label{AMamp}
\begin{aligned}
H_s(\omega,t)=&\frac{Ae^{i(\nu t+\varphi)}}{2}\left[\hat{h}_{\nu}(\omega)+\frac{r_a}{2}\hat{h}_{\nu+\nu_a}(\omega)e^{i\phi_a(t)}
+\frac{r_a}{2}\hat{h}_{\nu-\nu_a}(\omega)e^{-i\phi_a(t)}\right],\\
|H_s(\omega,t)|^2&=\frac{A^2}{4}{\Big[}\hat{h}_{\nu_1}^2(\omega)+\frac{r_a^2}{4}{\Big(}\hat{h}_{\nu+\nu_a}(\omega)+\hat{h}_{\nu-\nu_a}(\omega){\Big)}^2\\
&+r_a\hat{h}_{\nu}(\omega){\Big(}\hat{h}_{\nu+\nu_a}(\omega)+\hat{h}_{\nu-\nu_a}(\omega){\Big)}\cos \phi_a(t)\\
&-r_a^2\hat{h}_{\nu+\nu_a}(\omega)\hat{h}_{\nu-\nu_a}(\omega)\sin^2\phi_a(t){\Big]},\\
\nu_H(\omega,t)&=\nu+\frac{r_a\nu_a}{2}\frac{A^2}{4|H_s(\omega,t)|^2}
\left(\hat{h}_{\nu+\nu_a}(\omega)-\hat{h}_{\nu-\nu_a}(\omega)\right)\\
&\times\left( \frac{r_a}{2}[\hat{h}_{\nu+\nu_a}(\omega)+\hat{h}_{\nu-\nu_a}(\omega)]+\hat{g}_{\nu}(\omega)\cos\phi_a(t) \right).\\
\end{aligned}
\end{equation}
Comparing with the two-tone case (\ref{OVamp}), we now have two interference terms, $\sim\cos\phi_a(t)$ and $\sim\sin^2\phi_a(t)$. They are responsible for interference between the main tone and the AM-induced ones, and between the two AM-induced ones in their own rights, respectively. Note that all time behavior depends on $\phi_a(t)$, so that the phase shifts $\varphi,\varphi_a$ in (\ref{AMs}) do not change anything qualitatively.

\subsubsection{Representation}

\begin{figure*}[t!]
\includegraphics[width=1.0\linewidth]{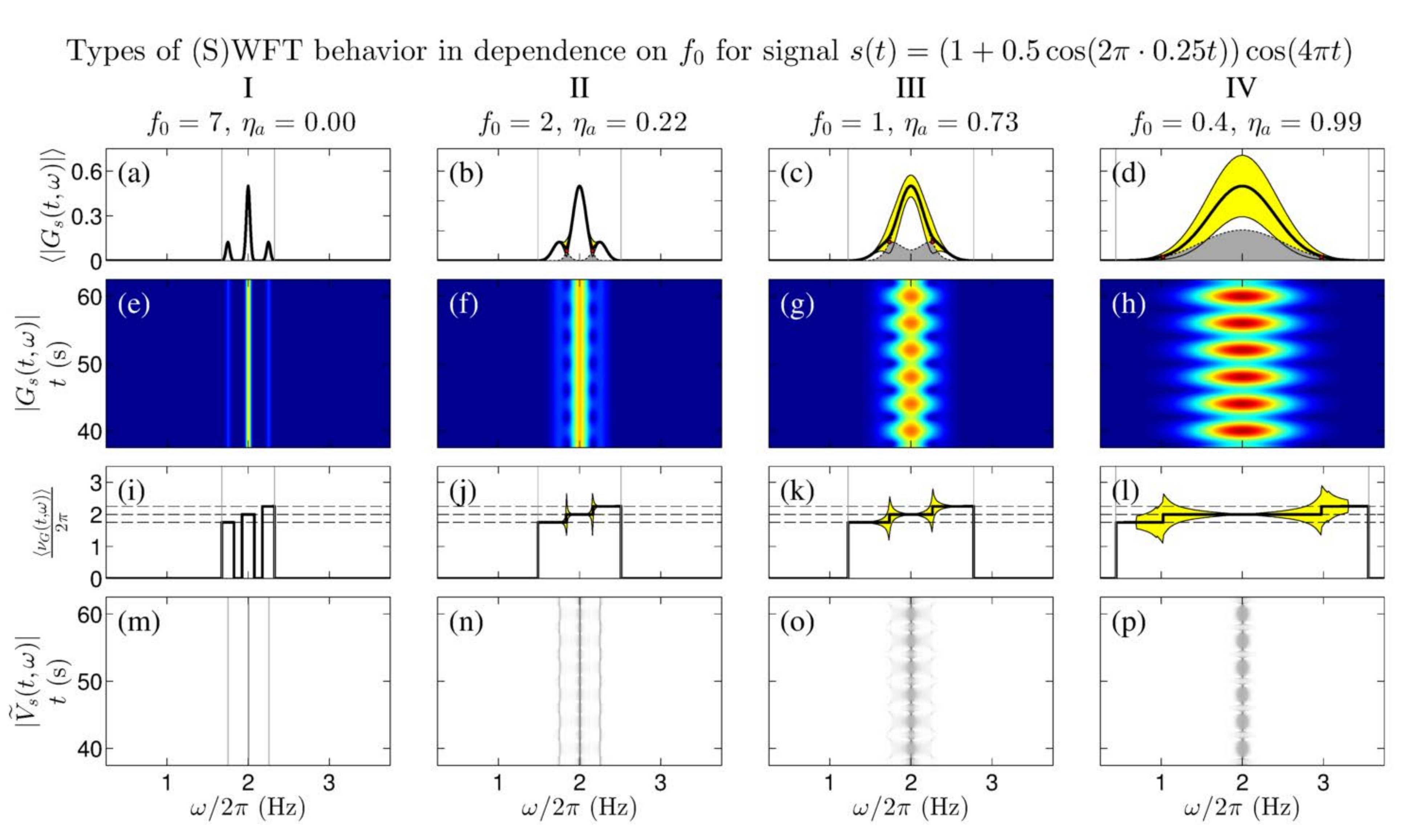}\\
\caption{Behavior of the Gaussian window WFT in dependence on $f_0$ for amplitude-modulated signal $s(t)=(1+0.5\cos(2\pi t/4))\cos(4\pi t)$, sampled at $20$ Hz for $500$ s. For illustrational purposes, Gaussian window $\hat{g}(\xi)$ is ``cutted'' to compact frequency support $[\xi_1(0.001),\xi_2(0.001)]$, and the boundaries of the joint support of all tones $[\nu-\nu_a+\xi_1(0.001),\nu+\nu_a+\xi_2(0.001)]$ are shown by gray lines in (a-d) and (i-l). (a-d): Time-averaged WFT amplitudes; dotted lines show $\frac{1}{2}\hat{g}(\omega-\nu)$ and $\frac{1}{2}\frac{r_a}{2}[\hat{g}(\omega-\nu-\nu_a)+\hat{g}(\omega-\nu+\nu_a)]$, with red dots indicating points of intersection between these functions and gray region showing the area shared by both of them. (e-h): WFT amplitudes in time-frequency domain. (i-l): Time-averaged WFT frequency $\nu_{G}(\omega,t)$, with dashed lines showing the frequencies of each tone $\nu,\nu-\nu_a,\nu+\nu_a$. (m-p): SWFT amplitudes in time-frequency domain. Values of $\eta_a$ indicate relative overlap defined in (\ref{AMeta}) (ratio of gray-shaded area to all area below the lower dotted line in (a-d), see text and (\ref{AMeta})). The rapid ``cuts'' in yellow regions in (l) occur at places where the support of the main peak ends (due to restricting $\hat{g}(\xi)$ to a finite support).}
\label{AMbt}
\end{figure*}

From (\ref{AMs}) it is clear that one can look onto AM component from different perspectives: either as an oscillation with amplitude modulation, or as the superposition of three tones with particular amplitude, phase and frequency relationships. The way in which the component (\ref{AMs}) is represented in the TFR is determined by the window/wavelet parameters: if the time resolution is large, it will be treated as a single component; if the frequency resolution is large, it will be treated as three independent tones. Different types of TFR behavior (see Table \ref{BTtab}) for the AM component (\ref{AMs}) are illustrated in Fig.\ \ref{AMbt} for the Gaussian window (S)WFT; for other windows, as well as for the (S)WT, all is qualitatively the same. Note that, in contrast to the previously considered case of two interfering tones, where the desired behavior was that of the Regime I type, we now want the AM component to be represented consistently as a single entity in the TFR, which corresponds to the Regime IV.

Similarly to (\ref{OVeta}), to quantify TFR behavior in the present case we can use the interference measure $\eta_a$, based on (\ref{OVeta}) and the analogy of the current expression for $|H_s(\omega,t)|^2$ (\ref{AMamp}) with that obtained previously for two tones (\ref{OVamp}). Thus, we define $\eta_a$ as the ratio of the area shared by $\frac{A}{2}\hat{h}_{\nu}$ (the TFR amplitude of the main tone separately) and $\frac{A}{2}\frac{r_a}{2}[\hat{h}_{\nu+\nu_a}(\omega)+\hat{h}(\nu-\nu_a)]$ (the sum of the TFR amplitudes of each of the AM-induced tones), and which is shown as the gray-shaded area in Fig.\ \ref{AMbt} (b-d), to the total area under the latter:
\begin{equation}\label{AMeta}
\begin{aligned}
\eta_a
=&\frac{\int {\min}[\hat{h}_{\nu}(\omega),\frac{r_a}{2}(\hat{h}_{\nu+\nu_a}(\omega)+\hat{h}_{\nu-\nu_a}(\omega))]d\mu(\omega)}
{\frac{r_a}{2}\int [\hat{h}_{\nu+\nu_a}(\omega)+\hat{h}_{\nu-\nu_a}(\omega)]d\mu(\omega)}\\
=&\frac{1}{\frac{r_a}{2}\int [\hat{h}_{\nu+\nu_a}(\omega)+\hat{h}_{\nu-\nu_a}(\omega)]d\mu(\omega)}
\Bigg(
\int_{-\infty}^{\mu(\omega_\times^{(1)})}\hat{h}_\nu(\omega)d\mu(\omega)\\
&+\frac{r_a}{2}\int_{\mu(\omega_\times^{(1)})}^{\mu(\omega_\times^{(2)})}[\hat{h}_{\nu+\nu_a}(\omega)+\hat{h}_{\nu-\nu_a}(\omega)]d\mu(\omega)
+\int_{\mu(\omega_\times^{(1)})}^{\infty}\hat{h}_\nu(\omega)d\mu(\omega)
\Bigg)\\
=&\frac{1}{2}\widetilde{Q}_{\nu+\nu_a}(\omega_\times^{(1)},\omega_\times^{(2)})
+\frac{1}{2}\widetilde{Q}_{\nu+\nu_a}(\omega_\times^{(1)},\omega_\times^{(2)})
+\frac{1}{r_a}[1-\widetilde{Q}_{\nu}(\omega_\times^{(1)},\omega_\times^{(2)})]\\
\end{aligned}
\end{equation}
where the intersection frequencies $\omega_\times^{(1,2)}$ are determined as
\begin{equation}\label{AMcross}
\begin{gathered}
\omega_\times^{(1)}<\nu,\omega_{\times}^{(2)}>\nu:\;
\hat{h}_\nu(\omega_\times^{(1,2)})=\frac{r_a}{2}[\hat{h}_{\nu-\nu_a}(\omega_\times^{(1,2)})+\hat{h}_{\nu+\nu_a}(\omega_\times^{(1,2)})],\\
\begin{aligned}
\Rightarrow\omega_\times^{(1,2)}=&\nu\pm\frac{\log\bigg[r_a^{-1}e^{f_0^2\nu_a^2/2}+\sqrt{r_a^{-2}e^{f_0^2\nu_a^2}-1}\bigg]}{f_0^2\nu_a}\\
&\mbox{ for the Gaussian window WFT}
\end{aligned}
\end{gathered}
\end{equation}
If there are no solutions $\omega_\times^{(1)}<\nu$ ($\omega_\times^{(2)}>\nu$) of (\ref{AMcross}), we take $\mu(\omega_{\times}^{(1)})=-\infty$ ($\mu(\omega_{\times}^{(2)})=\infty$), while if there are few solutions we take the ones closest to $\nu$ from each side. Note that the jumps in $\langle\nu_G(\omega,t)\rangle$ seen in Fig.\ \ref{AMbt} (j-l) occur almost exactly at the intersection frequencies (as was noticed previously for the two-tone signal); this is observed for the WT as well.

\begin{center}
\begin{table*}[t!]
  \begin{subequations}
  \begin{tabular}{| >{\centering\arraybackslash}m{1.5cm} | m{16cm} |}
  \hline
Regime & \multicolumn{1}{c|}{Condition}
\\ \hline
\textbf{I} &
\begin{equation}\label{AMt1}
\eta_a\leq\epsilon
\end{equation}
\begin{equation}\label{AMt1a}
\overset{approx.}{\Rightarrow}
\left[\begin{array}{rl}
\mbox{(S)WFT:} & \nu_a>\xi_2(\epsilon)-\xi_1(\epsilon),\\
\mbox{(S)WT:} & 1+\nu_a/\nu>\xi_2(\epsilon)/\xi_1(\epsilon)\\
\end{array}\right.
\end{equation}
\\ \hline
\textbf{II} &
\begin{equation}\label{AMt2}
\left\{\begin{array}{l}
\eta_a>\epsilon\\
\langle N_p\rangle=3\overset{approx.}{\Rightarrow}
[\hat{h}_{\nu_1}(\omega)+\frac{r_a}{2} (\hat{h}_{\nu+\nu_a}(\omega)+\hat{h}_{\nu-\nu_a}(\omega))]\mbox{ has two minimums in }\omega\in[\nu-\nu_a,\nu+\nu_a]\\
\end{array}\right.
\end{equation}
\\ \hline
\textbf{III} &
\begin{equation}\label{AMt3}
1<\langle N_p\rangle<3\;\overset{approx.}{\Rightarrow}\;
\left\{\begin{array}{l}
\eta_a< 1-\epsilon\\
{[}\hat{h}_{\nu_1}(\omega)+\frac{r_a}{2} (\hat{h}_{\nu+\nu_a}(\omega)+\hat{h}_{\nu-\nu_a}(\omega))]
\mbox{ has less than two minimums in }\omega\in[\nu-\nu_a,\nu+\nu_a]\\
\end{array}\right.
\end{equation}
\\ \hline
\textbf{IV} &
\begin{equation}\label{AMt4}
\eta_a\geq 1-\epsilon
\end{equation}
\begin{equation}\label{AMt4a}
\overset{approx.}{\Rightarrow}
\left[\begin{array}{rl}
\mbox{(S)WFT:}&
\left\{\begin{array}{l}
\omega_{\times}^{(1)}\leq\nu-\nu_a+\xi_1(2\epsilon)\\
\omega_\times^{(2)}\geq\nu+\nu_a+\xi_2(2\epsilon)\\
\end{array}\right.
\Rightarrow
r_a\leq\frac{e^{f_0^2\nu_a^2/2}}{\cosh[f_0\nu_a(n_G(2\epsilon)+f_0\nu_a)]}
\mbox{ for Gaussian window}\\
 & \\
\mbox{(S)WT:}&
\left\{\begin{array}{l}
\omega_{\times}^{(1)}\leq(\nu-\nu_a)\frac{\omega_\psi}{\xi_2(2\epsilon)}\\
\omega_{\times}^{(2)}\geq(\nu+\nu_a)\frac{\omega_\psi}{\xi_1(2\epsilon)}
\end{array}\right.
\\
\end{array}\right.
\end{equation}
\\ \hline
  \end{tabular}
  \end{subequations}
\caption{Conditions for each type of behavior (illustrated in Fig.\ \ref{AMbt}) for the AM component (\ref{AMs}), where we have used notations (\ref{AMnt}). Value of $\epsilon$ is some predefined accuracy (we use $\epsilon=0.001$) that determines how high (low) the interference should be to regard the tones as fully merged (separated) in the TFR.}
\label{tab:AMc}
\end{table*}
\end{center}

The conditions for each type of TFR behavior in the case of the AM component (\ref{AMs}) can be derived in a similar way as was done for a two-tone signal; they are summarized in Table \ref{tab:AMc}. The expressions determining the borders of Regimes I and IV, (\ref{AMt1}) and (\ref{AMt4}), are closely similar to those devised previously for the two-tone signal, (\ref{OVt1}) and (\ref{OVt4}), and the same considerations apply. Thus, e.g.\ (\ref{AMt4a}) are stricter than (\ref{AMt4}), implying the latter but being not implied by it. The conditions for the Regime II and III behavior types, (\ref{AMt2}) and (\ref{AMt3}), are also derived in a similar manner to that done for the two-tone case. Thus, for the signal (\ref{AMs}) the three tones appear to be most-merged and most-separated in the TFR when $\phi_a(t)=0$ and $\phi_a(t)=\pi$, respectively. In these cases TFR amplitude (\ref{AMamp}) simplifies to
\begin{equation}\label{AMms}
\left.|H_s(\omega,t)|^2\right|_{\phi_a(t)=0,\pi}=\frac{A^2}{4}\Big[\hat{h}_{\nu}(\omega)\pm \frac{r_a}{2}(\hat{h}_{\nu+\nu_a}(\omega)+\hat{h}_{\nu-\nu_a}(\omega))\Big]^2
\end{equation}
which then leads to (\ref{AMt2}) and (\ref{AMt3}). Note, however, that now the statement that there is minimum (maximum) number of peaks for $\phi_a(t)=0$ ($\phi_a=\pi$) in the general case is not exact, but only an approximation (due to the appearance of an additional interference term $\sim \sin^2\phi_a$ in $|H_s(\omega,t)|^2$ (\ref{AMamp})), although the quality of this approximation is usually excellent.

\begin{figure*}[t!]
\includegraphics[width=1.0\linewidth]{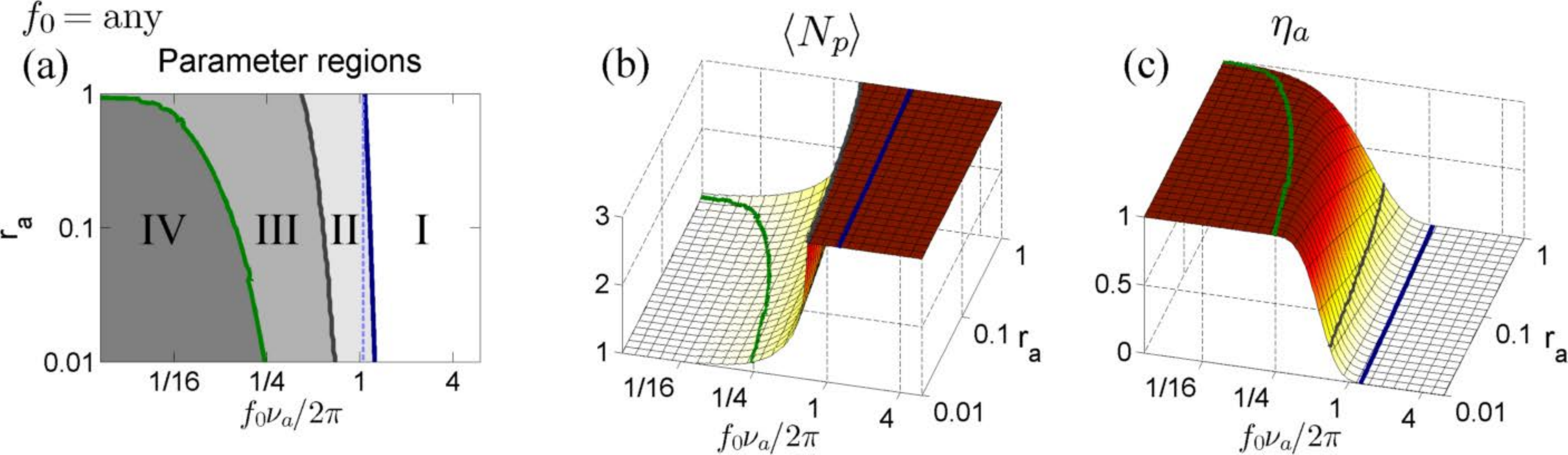}\\
\caption{Dependence of the WFT behavior on signal (\ref{AMs}) parameters $r_a,\nu_a$ and Gaussian window resolution parameter $f_0$. (a): Regions of parameter space corresponding to each type of behavior, according to (\ref{AMt1}),(\ref{AMt2}),(\ref{AMt3}),(\ref{AMt4}); dashed light-blue line shows boundary of the I-type behavior as predicted by approximate (\ref{AMt1a}). (b): Mean number of peaks $\langle N_p\rangle$ (\ref{NP}). (c): Relative overlap $\eta_a$ (\ref{AMeta}). For determining Regimes I and IV we used $\epsilon=0.001$ in (\ref{AMt1}) and (\ref{AMt4}); for all $f_0$ we assume that $\nu_a<\nu$.}
\label{AMbdWFT}
\end{figure*}

For the AM component (\ref{AMs}), the behavior of the Gaussian window WFT, as characterized by mean number of peaks $\langle N_p\rangle$ (\ref{NP}) and $\eta_a$ (\ref{AMeta}), is illustrated in Fig.\ \ref{AMbdWFT} in terms of its dependence on signal parameters $r,\nu,\nu_a$ and the window resolution parameter $f_0$. Similarly to the case of a two-tone signal, all depends only on $r_a$ and $f_0\nu_a$, and not on $f_0$ or $\nu_a$ separately. Note, that the agreement between (\ref{AMt1}) and (\ref{AMt1a}) (solid blue and dashed light-blue lines in Fig.\ \ref{AMbdWFT} (a)) is now not so good as it was for the two-tone signal, although still satisfactory.

\begin{figure*}[t!]
\includegraphics[width=1.0\linewidth]{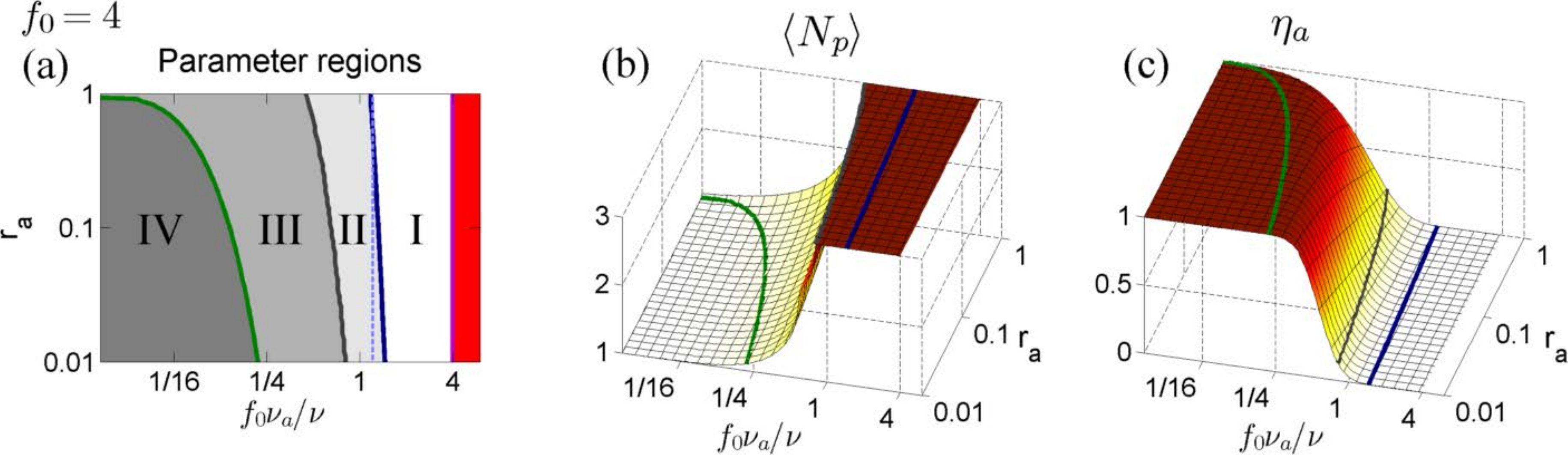}\\
\includegraphics[width=1.0\linewidth]{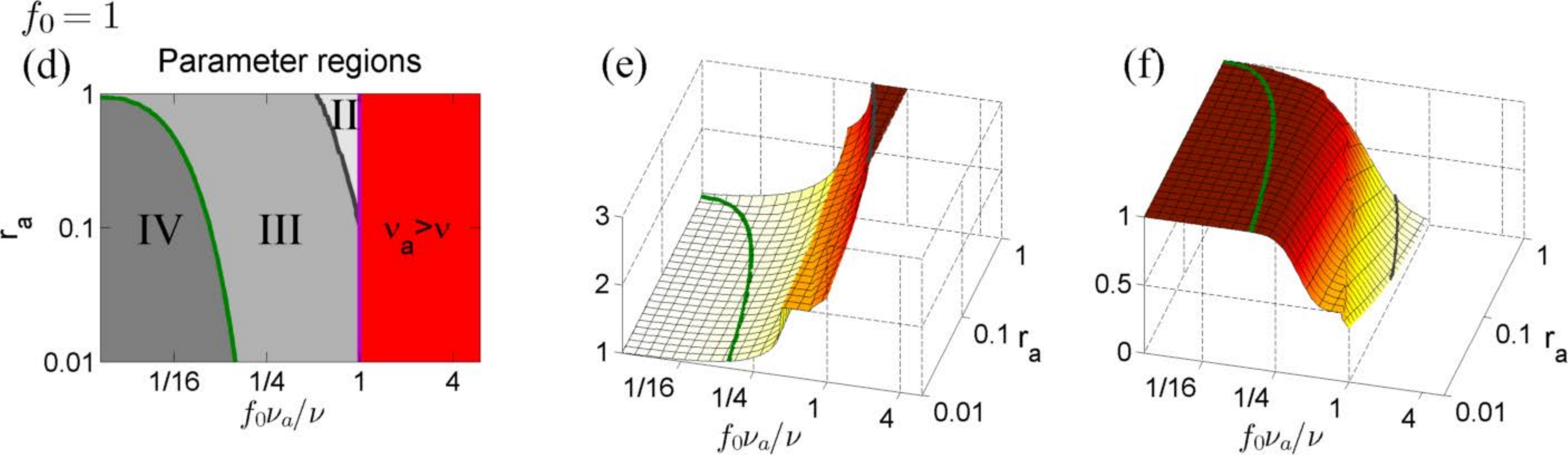}\\
\caption{Same as Figure \ref{AMbdWFT}, but for the Morlet wavelet WT. In (a,d), red filling indicates an inappropriate regions where $\nu_a>\nu$; in (b) and (c) such regions are omitted (as $\eta_a$ (\ref{AMeta}) is not well-defined for them).}
\label{AMbdWT}
\end{figure*}

The corresponding behavior characteristics of the Morlet wavelet WT are shown in Fig.\ \ref{AMbdWT}. In contrast to WFT, where by choosing a high enough $\nu$ we can in principle investigate any $f_0\nu_a$ subject to the condition $\nu_a<\nu$, for the WT the meaningful parameter is the frequency ratio $\nu_a/\nu$, so that one is restricted to $f_0\nu_a/\nu\leq f_0$, establishing an $f_0$-dependent cut-off in $f_0\nu_a/\nu$; the ``prohibited'' region $\nu_a>\nu$, where the amplitude becomes ill-defined due to its varying faster that the main oscillation (and as a result analytic approximation becomes inaccurate, see Part I) is shown as red in Fig.\ \ref{AMbdWT}(a). Additionally, as was mentioned before (see discussion of Fig.\ \ref{OVbdWT}), the WT behavior, apart from $r_a$ and $f_0\nu_a/\nu$, also depends on $f_0$ separately. However, in the present case this will be so even for the lognormal wavelet, whose behavior will also depend on three parameters: $r_a$, $f_0\log(1+\frac{\nu_a}{\nu})$ and $f_0\log(1-\frac{\nu_a}{\nu})$. This is because the AM-induced tones are located around the main one symmetrically on a linear, but not logarithmic frequency scale. Note that, for both the WFT and WT, Regime IV occupies a much wider parameter space compared to what was seen for the two-tone signal, meaning that due to the specific relationships between amplitudes, phases and frequencies of the tones in (\ref{AMs}), which characterize amplitude modulation, it is much easier for them to be treated by the TFR as a single component. Note also, that for $r_a\rightarrow 1$, corresponding to the maximum allowed value, Regime IV is generally not possible, as seen from Fig.\ \ref{AMbdWFT} and Fig.\ \ref{AMbdWT}.

\subsubsection{Reconstruction}

We now investigate the performance of different reconstruction methods for the amplitude-modulated signal (\ref{AMs}). In all TFRs, we use ``maximum-based'' curve extraction, selecting the TFS around the maximum peak in the TFR amplitude at each time. Reconstruction methods are then applied, the resultant signals are compared with their true values, and the errors $\varepsilon_{a,\phi,f}$ (\ref{apfrel}) are calculated.

\begin{figure*}[t!]
\includegraphics[width=1.0\linewidth]{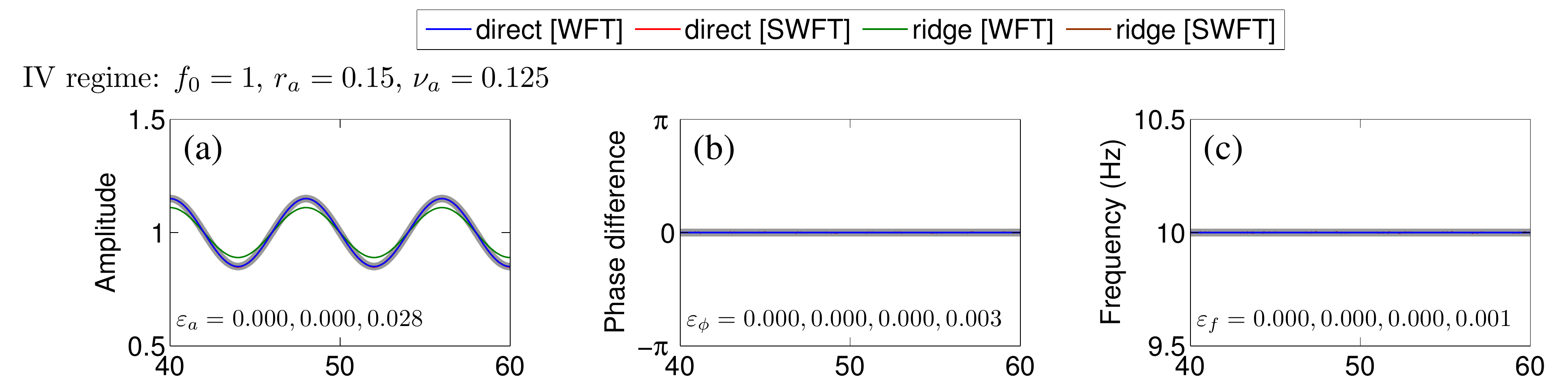}\\
\includegraphics[width=1.0\linewidth]{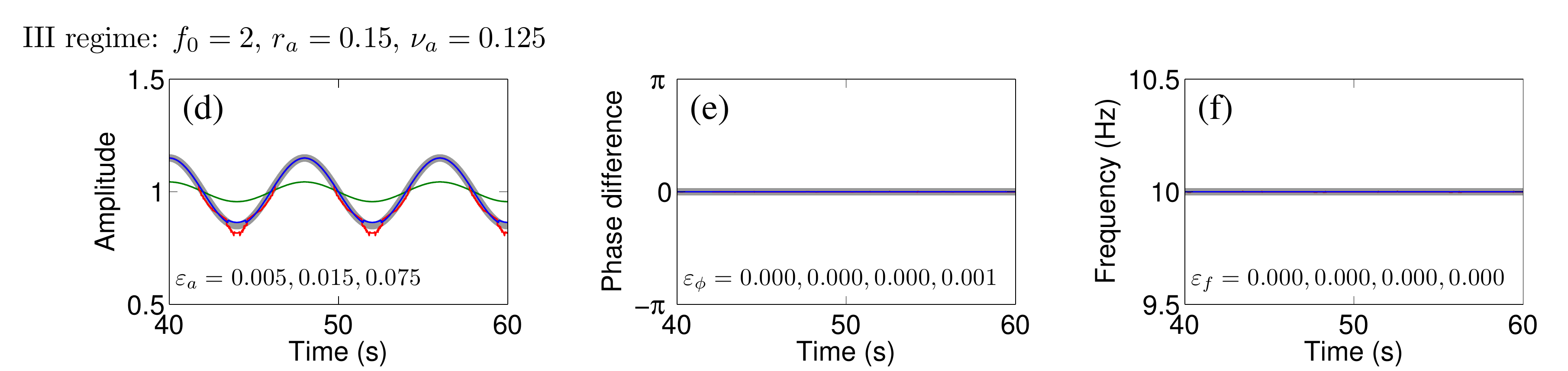}\\
\caption{Amplitude (a,d), phase (b,e) and frequency (c,f) of the AM component (\ref{AMs}) as reconstructed from its WFT and SWFT (colored lines), compared to the true values (thick gray lines). Values of $\varepsilon_{a,\phi,f}$ shown are in the same order as lines in legend, corresponding to ${\rm direct [WFT]}$ (blue), ${\rm direct [SWFT]}$ (red), ${\rm ridge [WFT]}$ (green), ${\rm ridge [SWFT]}$ (brown). In (a,d), ridge reconstruction from the SWFT is not shown as it is not appropriate for amplitude (see Part I). In (b,e), the difference between the reconstructed and true phase is shown. The signal (\ref{AMs}) was sampled at $50$ Hz for $100$ s, and it was simulated with $\varphi=\varphi_a=0$ and $\nu/2\pi=10$; all other parameters are indicated on the figure.}
\label{AMex}
\end{figure*}

Fig.\ \ref{AMex} shows examples of amplitudes, phases and frequencies reconstructed from the (S)WFT in two cases, corresponding to Regimes IV and III; for the (S)WT all remains qualitatively similar. First of all, in both regimes, phase/frequency are reconstructed perfectly by all methods, so amplitude modulation (without frequency modulation) does not influence their reconstruction, which is true for all Regimes I-IV. However, this is typical only for the WFT with symmetric windows $\hat{g}(\xi)$: in this case the peak in the WFT amplitude occurs exactly at the main tone frequency $\omega_p(t)=\nu$ and the WFT is symmetric around it, so one can recover the exact phase/frequency by direct methods; furthermore, one has $\nu_G(\omega_p(t),t)=\nu_G(\nu,t)=\nu$, as follows from (\ref{AMamp}), so that the ridge frequency estimate will be also exact, and the same can be shown for the phase. On the other hand, for the WFT with frequency-asymmetric windows, as well as for the WT in general (since it is inherently asymmetric for any wavelet, although can be symmetric on a logarithmic frequency scale), the amplitude peak will be not located at $\nu$, and the TFR will be also be asymmetric around the peak. As a result, amplitude modulation will introduce errors in both the ridge and direct estimates of the phase and frequency (although direct estimates remain exact for Regime IV). Therefore, (S)WFTs based on windows symmetric in frequency generally offer more accurate phase/frequency estimation than other TFRs.

Next, we see that for the amplitude reconstruction direct methods greatly outperform ridge ones, giving almost perfect estimates both in the regions of Regimes IV and III. This is to be expected, since by definition direct methods should give exact estimates in the case when extracted TFS contain all component. Thus, the quality of direct estimates can be worsened only by the separation of amplitude modulation into few peaks, so that the extracted TFS will contain only part of the component. Indeed, in panel (d) one can see (especially for ${\rm ridge [SWFT]}$), that direct methods become most inaccurate when amplitude is minimal ($\phi_a(t)=\pi$), which correspond to maximum separation between tones (and so occurrence of additional peaks in Regime III, although for the presented in Fig.\ \ref{AMex} case they are small). At the same time, reconstruction from ridges gives approximate estimates, which depend explicitly on the law of amplitude modulation (see Sec.\ \ref{sec:optrec} below): the more pronounced and fast it is (in respect to window/wavelet time resolution), the more ridge reconstruction is inaccurate. This concerns phase and frequency estimates as well -- except in the case of WFTs with symmetric $\hat{g}(\xi)$, for which they are exact.

\begin{figure*}[t!]
\includegraphics[width=1.0\linewidth]{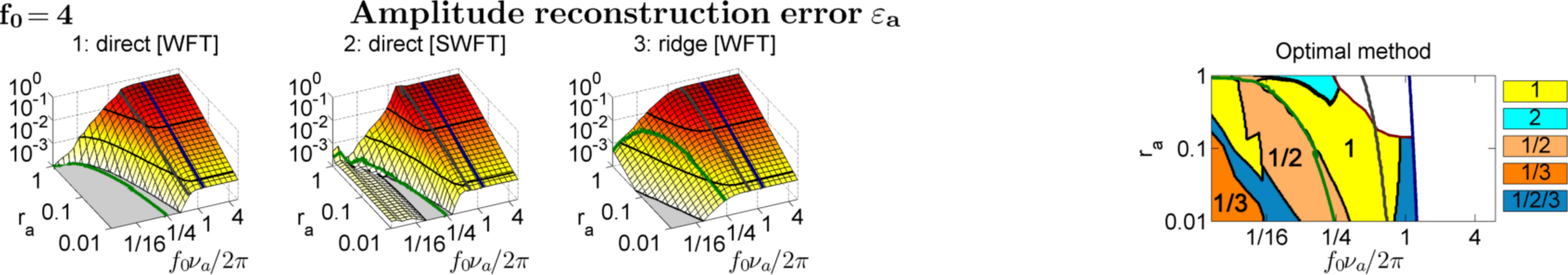}\\
\includegraphics[width=1.0\linewidth]{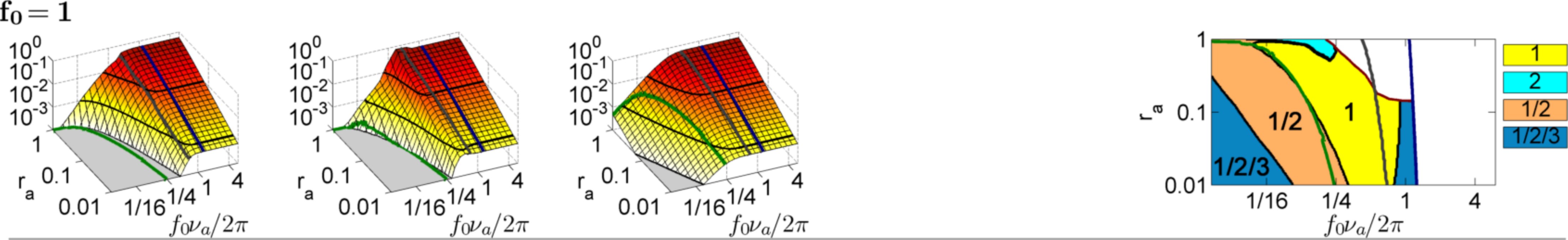}\\
\includegraphics[width=1.0\linewidth]{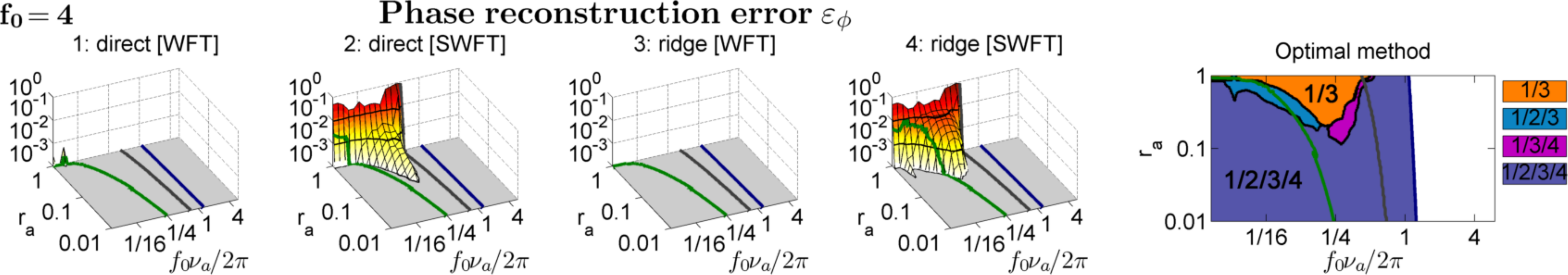}\\
\includegraphics[width=1.0\linewidth]{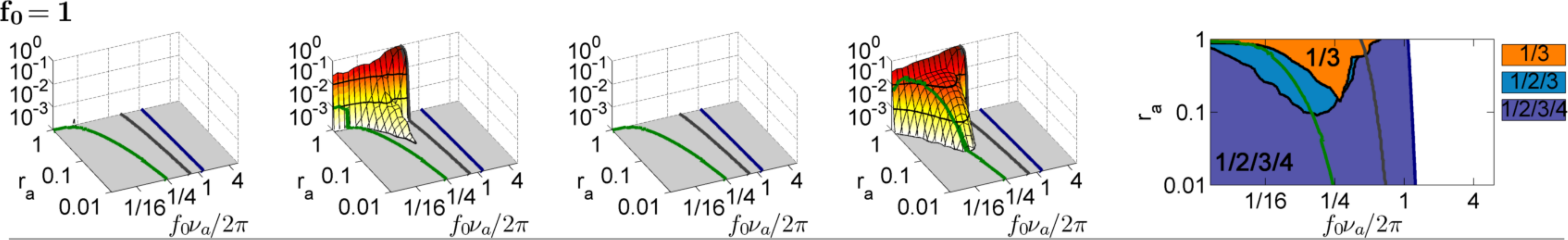}\\
\includegraphics[width=1.0\linewidth]{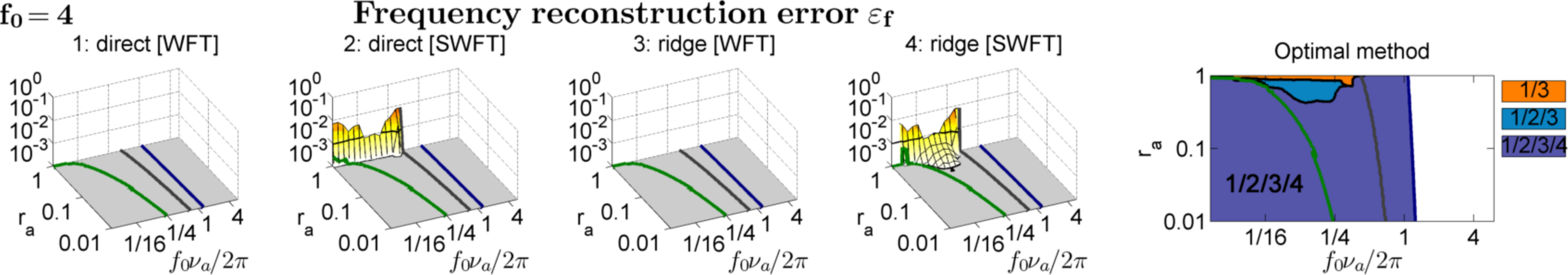}\\
\includegraphics[width=1.0\linewidth]{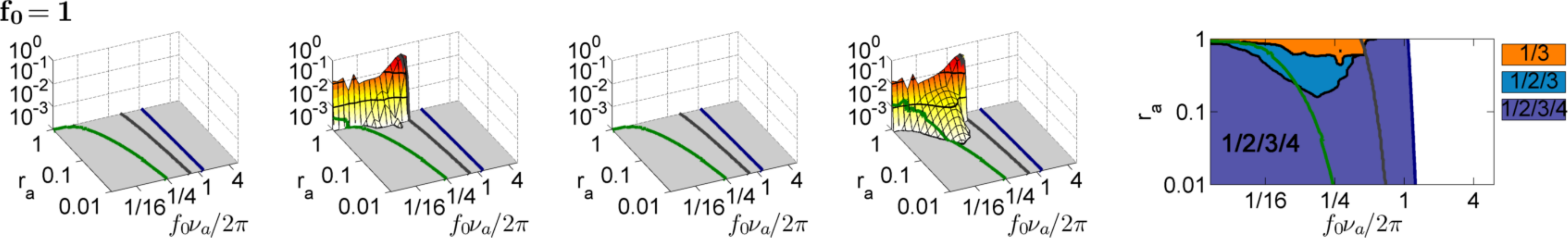}\\
\caption{Left four panels in each row show the (S)WFT-based amplitude (1-2 rows), phase (3-4 rows) and frequency (5-6 rows) reconstruction errors for the AM component $s(t)=(1+r_a\cos(\nu_a t))\cos(\nu t)$ in dependence on parameters $r_a,\nu_a$ and Gaussian window resolution parameter $f_0$. We used $\nu/2\pi=10$, but results do not depend on its value, at least assuming $\nu_a<\nu$ (so the analytic signal approximation remains valid). Estimation of amplitude from SWFT ridges is not appropriate, and so the corresponding results are discarded. Thin black lines indicate the levels of $0.001,0.01,0.1,1$, and the behavior of errors $<0.001$ for simplicity is not shown; thick green, gray and blue lines are the same as in Fig. \ref{AMbdWFT}, showing the borders of regions corresponding to behavior in Regimes I-IV. The right panels in each row show regions in parameter space where each method is optimal, i.e.\ gives the smallest estimation error; if the resultant errors of two methods differ on less than $0.001$, they are regarded as having similar performance (the corresponding regions are denoted as `'method1/method2`'). Comparison does not make much sense when WFT behavior is of I type or if all reconstruction errors are $>0.1$, i.e.\ all methods fail; the corresponding regions are white-filled in the right-hand panels. The signal was sampled at 50 Hz for 100 s.}
\label{AMrecWFT}
\end{figure*}

\begin{figure*}[t!]
\includegraphics[width=1.0\linewidth]{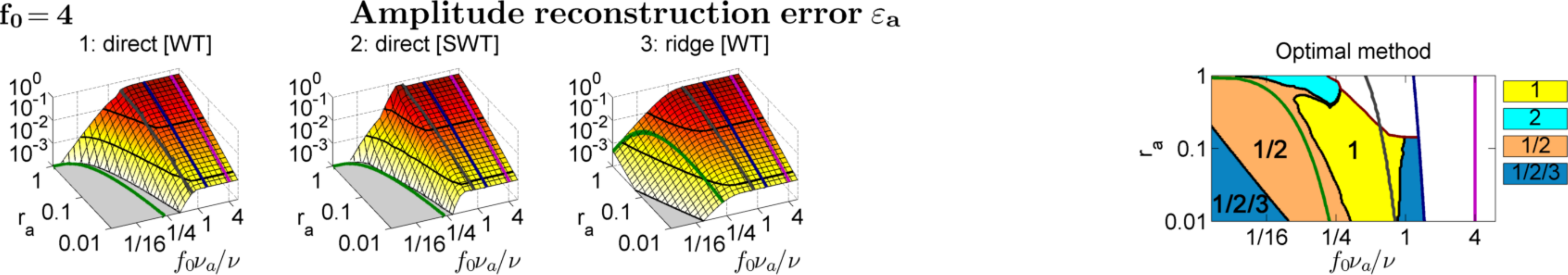}\\
\includegraphics[width=1.0\linewidth]{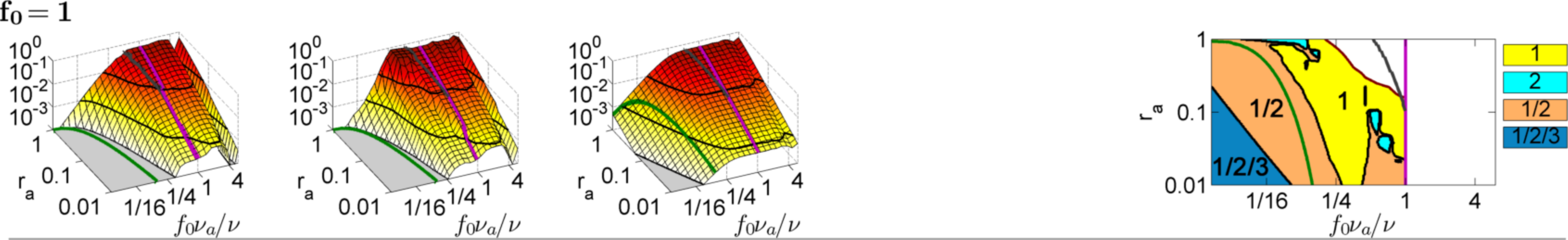}\\
\includegraphics[width=1.0\linewidth]{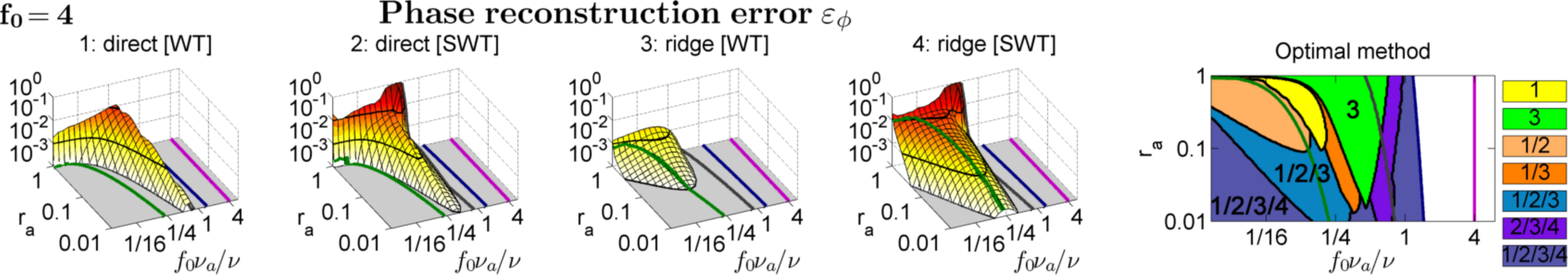}\\
\includegraphics[width=1.0\linewidth]{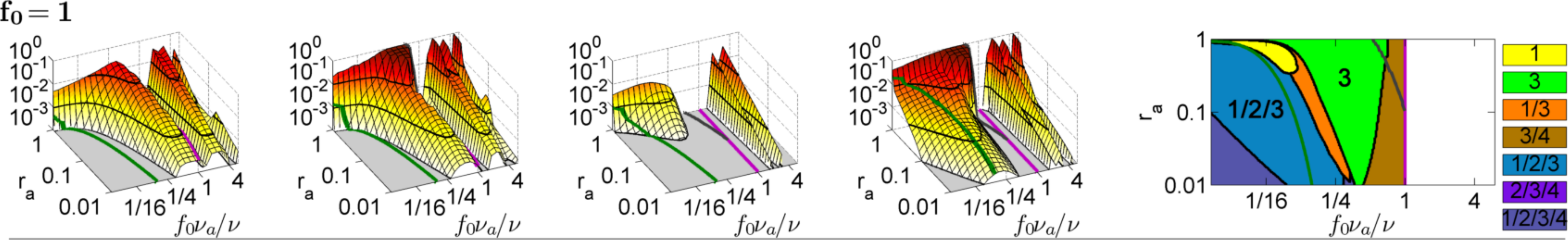}\\
\includegraphics[width=1.0\linewidth]{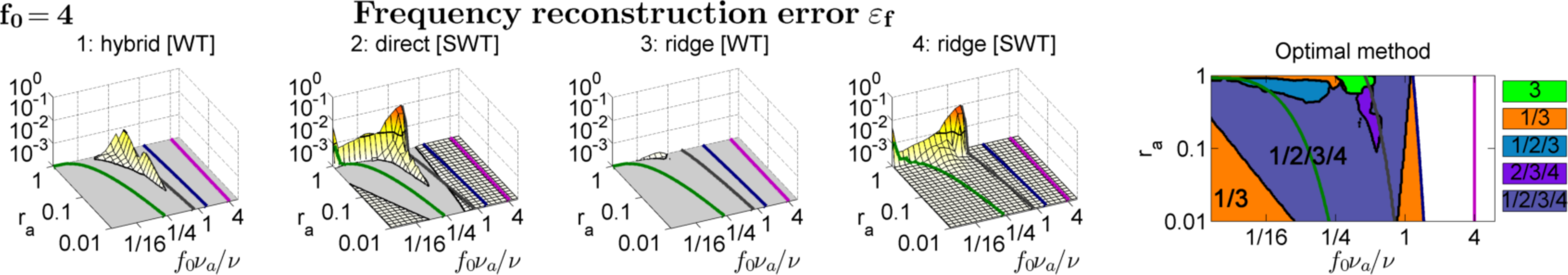}\\
\includegraphics[width=1.0\linewidth]{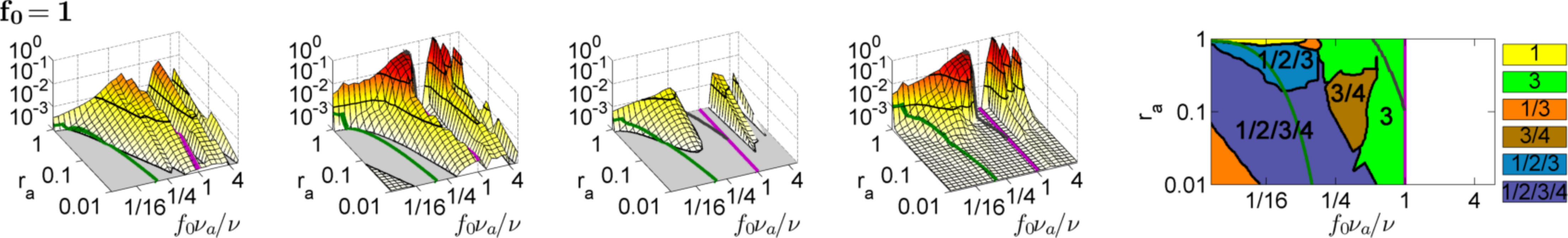}\\
\caption{Same as Fig.\ \ref{AMrecWFT}, but all parameters are reconstructed from the signal's WT/SWT. The direct frequency reconstruction for Morlet wavelet is not possible, so we use hybrid reconstruction (see Part I). Thick magenta lines separate the parameter regions where $\nu_a>\nu$, so that the analytic signal approximation is not valid there (implying additional theoretical error which cannot be reduced and thus making any comparisons not appropriate). Note that the frequency reconstruction error for the WT (\ref{apfrel}) is defined in relation to the component's mean frequency, so $\varepsilon_f$ here is divided on $\nu/2\pi=10$ in comparison to $\varepsilon_f$ in Fig.\ \ref{AMrecWFT}; this makes depicted dependencies consistent for different frequencies $\nu$ in (\ref{AMs}).}
\label{AMrecWT}
\end{figure*}

Fig.\ \ref{AMrecWFT} shows how the reconstruction errors for the (S)WFT depend on the signal and window/wavelet parameters; Fig.\ \ref{AMrecWT} provides the corresponding information for the (S)WT. In the former case, the error depends mainly on $f_0\Delta\nu$, whereas for the (S)WT there is an additional, separate, dependence on $f_0$. As discussed above, the phase/frequency estimates are exact for the Gaussian window WFT (though not the Morlet wavelet WT), so logically they should be exact for SWFT as well. However, errors appear in SWFT-based estimates for high $r_a$ (see Fig.\ \ref{AMrecWFT}), related to the more complex behavior of the synchrosqueezed TFRs. This is because, for considerable amplitude modulation, the SWFT power is redistributed to the side frequencies at certain times, leaving only a small amount of power around the main frequency $\nu$. So a few TFSs appear, none of which contains the full component. The same applies to the SWT. Note, that in the present case one can alternatively use the frequency-based scheme for TFS extraction (instead of maximum-based), but it will not change anything for either (S)WFT- or (S)WT-based reconstruction: the only effect is to reduce (though not eliminate) errors in the SWFT/SWT-based phase/frequency estimation for large $r_a$.

Clearly, direct methods perform better than ridge-based ones for the amplitude-modulated signal (\ref{AMs}), and are exact in Regime IV, in contrast to the ridge methods. Synchrosqueezing does not give any advantages (only drawbacks) in terms of reconstruction, so that SWFT/SWT-based estimates are always of the same or worse quality than WFT/WT-based ones, by any method. The performance of all methods depends largely on the TFR behavior regime, so the choice of suitable window/wavelet parameters is of crucial importance.

\newpage

\newpage

\subsection{Frequency modulation}\label{sec:FM}

We now consider the FM component with sinusoidal frequency modulation, which can be represented as
\begin{equation}\label{FMs}
\begin{aligned}
s(t)=&A\cos(\nu t+\varphi+r_b\sin(\nu_b t+\varphi_b))\\
=&A\,{\rm Re}\sum_{n=-\infty}^\infty J_n(r_b)e^{i[(\nu+n\nu_b) t+(\varphi+n\varphi_b)]}\\
\end{aligned}
\end{equation}
where we have used the expansion $e^{ia\sin\phi}=\sum_{n=-\infty}^\infty J_n(a)e^{in\phi}$, with $J_n(a)=(-1)^nJ_{-n}(a)$ denoting the $n$th order Bessel functions of the first kind. For convenience we also denote
\begin{equation}\label{FMnt}
\phi(t)\equiv \nu t+\varphi;\;\phi_b(t)\equiv\nu_b t+\varphi_b
\end{equation}
Due to the definition of the phase, the parameters in (\ref{FMs}) should obey $r_b\nu_b<\nu$. However, for the analytic signal approximation to hold (see Part I), which is needed for a meaningful time-frequency analysis, one needs all non-negligible terms in the expansion (\ref{FMs}) to correspond to positive frequencies $\nu+n\nu_b>0$. This can be formulated as $n_J(r_b)\nu_b<\nu$, where $n_J(r_b)$ is determined as the maximal order of $J_n(r_b)$ for which it is non-negligible (for $r_b\in[0,1]$ one can safely take $n_J(r_b)=2$):
\begin{equation}\label{FMnj}
n_J(r_b)=\max[n:\;J_n(r_b)>\epsilon_J]
\end{equation}
where $\epsilon_J$ stands for the specified precision under which to regard a value as negligible (here we use $\epsilon_J=0.02$). Thus, one can effectively reduce the sum in (\ref{FMs}), as well as in all the following formulas, to just the terms with $|n|\leq n_J(r_b)$.

Clearly, the FM component (\ref{FMs}) can be represented in the form of a multitone signal: the main tone at frequency $\nu$, and many side tones of pairwise-equal amplitudes at $\nu\pm n\nu_b$; the latter appear due to the frequency modulation and will be called ``FM-induced''. Therefore, the classification and formulas for the multitone signals apply here as well. Using (\ref{gformGEN}), one obtains for the signal (\ref{FMs}):
\begin{equation}\label{FMtfr}
\begin{aligned}
Z_n^{(\pm)}(\omega)\equiv& \hat{h}_{\nu+n\nu_b}(\omega)\pm (-1)^n\hat{h}_{\nu-n\nu_b}(\omega)\\
H_s(\omega,t)=&\frac{A}{2}e^{i\phi(t)}
\bigg[J_0(r_b)+\sum_{n=1}^\infty J_n(r_b)\\
&\times\Big( \hat{h}_{\nu+n\nu_b}(\omega)e^{in\phi_b(t)}+(-1)^n\hat{h}_{\nu-n\nu_b}(\omega)e^{-in\phi_b(t)} \Big)
\bigg]\\
|H_s(\omega,t)|^2=&\frac{A^2}{4}{\bigg[}J_0^2(r_b)\hat{h}_\nu^2(\omega)
+2J_0(r_b)\sum_{n=1}^\infty J_n(r_b)Z_n^{(+)}(\omega)\cos n\phi_b(t)\\
&+\sum_{n,m=1}^\infty J_n(r_b)J_m(r_b)Z_n^{(+)}(\omega)Z_m^{(+)}(\omega)\cos[(n-m)\phi_b(t)]\\
&-4\sum_{n,m=1}^\infty J_n(r_b)J_m(r_b)\hat{h}_{\nu+n\nu_b}(\omega)\hat{h}_{\nu-m\nu_b}(\omega)\\
&\times \sin n\phi_b(t)\sin m\phi_b(t){\bigg]}\\
\nu_H(\omega,t)=&\nu+\nu_b\frac{A^2}{4|H_s(\omega,t)|^2}\sum_{n=1}^\infty nJ_n(r_b)Z_n^{(-)}(\omega)\\
&\times\sum_{m=0}^\infty J_m(r_b)Z_m^{(-)}(\omega)\cos[(n-m)\phi_b(t)]
\end{aligned}
\end{equation}
As can be seen, all is quite complicated, so even restricting the summations to $|n|,|m|\leq n_J(r_b)$, there still remains many interference terms. The only thing which can immediately be seen from (\ref{FMtfr}) is that the phase lags $\varphi,\varphi_b$ in (\ref{FMs}) do not influence anything and can be omitted from the analysis.

\subsubsection{Representation}

Similarly to the previously considered case of AM component, the FM component (\ref{FMs}) can be perceived from different viewpoints: either as an oscillation with frequency modulation; or as a multitone signal with particular relationships between the parameters of the tones. The way it is represented in the TFR is determined by the window/wavelet parameters: for large time resolution one will have a single component; and for large frequency resolution, implying small time resolution, it will be perceived as a sum of tones. Fig.\ \ref{AMbt} illustrates different types of (S)WFT behavior, according to Table \ref{BTtab}, for the FM component (\ref{FMs}); for (S)WFTs based on other windows, as well as for (S)WTs, all is qualitatively the same. Obviously, one usually aims at Regime IV while choosing the window/wavelet parameters in this case.

\begin{figure*}[t!]
\includegraphics[width=1.0\linewidth]{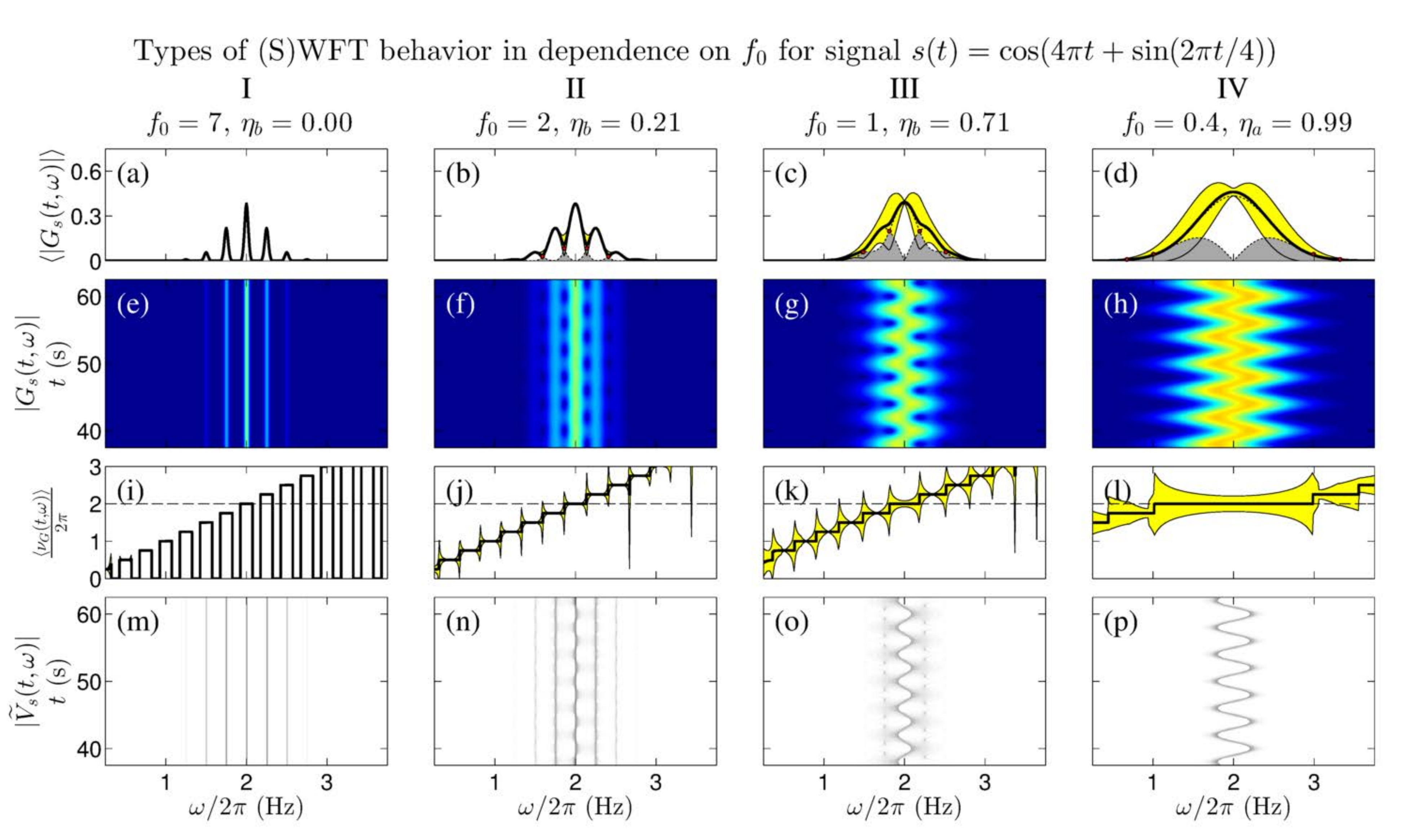}\\
\caption{Behavior of the Gaussian window WFT in dependence on $f_0$ for frequency-modulated signal $s(t)=\cos(4\pi t+\sin(2\pi t/4))$, sampled at $20$ Hz for $500$ s. For illustrational purposes, Gaussian window $\hat{g}(\xi)$ is ``cutted'' to compact frequency support $[\xi_1(0.001),\xi_2(0.001)]$. (a-d): Time-averaged WFT amplitudes; dotted lines show $\frac{1}{2}\hat{h}^{(+)}(\omega)$ and $\frac{1}{2}|\hat{h}^{(-)}(\omega)|$, defined in (\ref{FMhpm}), with red dots indicating first four points of intersection between these functions and gray region showing the area shared by both of them. (e-h): WFT amplitudes in time-frequency domain. (i-l): Time-averaged WFT frequency $\nu_{G}(\omega,t)$, with dashed line showing the frequency of the main tone $\nu$. (m-p): SWFT amplitudes in time-frequency domain. Values of $\eta_b$ indicate relative overlap defined in (\ref{FMeta}) (ratio of gray-shaded area to all area below the lower dotted line in (a-d), see text and (\ref{FMeta})).}
\label{FMbt}
\end{figure*}

To approach the relatively complicated case (\ref{FMtfr}), we introduce two measures:
\begin{equation}\label{FMhpm}
\begin{aligned}
&\hat{h}^{(+)}(\omega)=J_0(r_b)\hat{h}_{\nu}(\omega)+\sum_{n=1}^\infty J_{2n}(r_b)\Big[\hat{h}_{\nu+2n\nu_a}(\omega)+\hat{h}_{\nu-2n\nu_a}(\omega)\Big],\\
&\hat{h}^{(-)}(\omega)=\sum_{n=1}J_{2n-1}(\omega)\Big[\hat{h}_{\nu+(2n-1)\nu_a}(\omega)-\hat{h}_{\nu-(2n-1)\nu_a}(\omega)\Big].
\end{aligned}
\end{equation}
There are some good reasons for considering such quantities. First, it appears that the initial in $\langle\nu_H(\omega,t)\rangle$ (see Fig.\ \ref{FMbt}(j-l)) occur almost exactly at the intersection between $\hat{h}^{(+)}(\omega)$ and $|\hat{h}^{(-)}(\omega)|$. Secondly, consider the limit of vanishing frequency modulation $\nu_a\rightarrow0$, in which one should recover the TFR of the main tone. In this limit $\hat{h}^{(+)}(\omega)\rightarrow \hat{h}_\nu(\omega),\;\hat{h}^{(-)}(\omega)\rightarrow0$ for any $r$, i.e.\ all energy will be contained in $\hat{h}^{(+)}(\omega)$. Furthermore, for the WFT with windows symmetric in frequency, one has $\hat{h}^{(-)}(\nu)=0$, so all of the amplitude at the main tone frequency is concentrated in $\hat{h}^{(+)}(\nu)$.

Based on such considerations, we determine the interference measure $\eta_b$ for the FM signal (\ref{FMs}) as the area shared by $\hat{h}^{(+)}(\omega)$ and $|\hat{h}^{(-)}(\omega)|$ compared to the total area under $|\hat{h}^{(-)}(\omega)|$:
\begin{equation}\label{FMeta}
\eta_b=\frac{\int \min(\hat{h}^{(+)}(\omega),|\hat{h}^{(-)}(\omega)|)d\mu(\omega)}{\int|\hat{h}^{(-)}(\omega)|d\mu(\omega)}
\end{equation}
The value of $\eta_b$ is the intuitive measure of interference in our case: the larger it is -- the more the FM-induced tones interfere and behave as a single entity in the TFR. A simple manifestation of this is that for $\nu_b\rightarrow0$ one has $\eta_b\rightarrow1$, i.e.\ there is single component in the TFR.

\begin{center}
\begin{table*}[t!]
  \begin{subequations}
  \begin{tabular}{| >{\centering\arraybackslash}m{1.5cm} | m{16cm} |}
  \hline
Regime & \multicolumn{1}{c|}{Condition}
\\ \hline
\textbf{I} &
\begin{equation}\label{FMt1}
\eta_b\leq\epsilon
\end{equation}
\begin{equation}\label{FMt1a}
\overset{approx.}{\Rightarrow}
\left[\begin{array}{rl}
\mbox{(S)WFT:} & \nu_b>\xi_2(\epsilon)-\xi_1(\epsilon),\\
\mbox{(S)WT:} & 1+[n_J(r_b)+\nu/\nu_b]^{-1}>\xi_2(\epsilon)/\xi_1(\epsilon)\\
\end{array}\right.
\end{equation}
\\ \hline
\textbf{II} &
\begin{equation}\label{FMt2}
\left\{\begin{array}{l}
\eta_a>\epsilon\\
\langle N_p\rangle=2n_J(r_b)+1\overset{approx.}{\Rightarrow}
\left\{\begin{array}{l}
|\hat{h}^{(+)}(\omega)+\hat{h}^{(-)}(\omega)|\mbox{ has }2n_J(r_b)\mbox{ minimums}\\
|\hat{h}^{(+)}(\omega)-\hat{h}^{(-)}(\omega)|\mbox{ has }2n_J(r_b)\mbox{ minimums}\\
\end{array}\right.\\
\end{array}\right.
\end{equation}
\\ \hline
\textbf{III} &
\begin{equation}\label{FMt3}
1<\langle N_p\rangle<2n_J(r_b)+1\;\overset{approx.}{\Rightarrow}\;
\left\{\begin{array}{l}
\eta_b< 1-\epsilon\\
\left[\begin{array}{l}
|\hat{h}^{(+)}(\omega)+\hat{h}^{(-)}(\omega)|\mbox{ has less than }2n_J(r_b)\mbox{ minimums}\\
|\hat{h}^{(+)}(\omega)-\hat{h}^{(-)}(\omega)|\mbox{ has less than }2n_J(r_b)\mbox{ minimums}\\
\end{array}\right.\\
\end{array}\right.
\end{equation}
\\ \hline
\textbf{IV} &
\begin{equation}\label{FMt4}
\eta_b\geq 1-\epsilon
\end{equation}
\\ \hline
  \end{tabular}
  \end{subequations}
\caption{Conditions for each type of behavior (illustrated in Fig.\ \ref{FMbt}) for the FM component (\ref{AMs}), where we have used notations (\ref{FMnt}), (\ref{FMnj}) and (\ref{FMhpm}). Value of $\epsilon$ is some predefined accuracy (we use $\epsilon=0.001$) that determines how high (low) the interference should be to regard the tones as fully merged (separated) in the TFR.}
\label{tab:FMc}
\end{table*}
\end{center}

Restricting ourselves to terms with $|n|\leq n_J(r_b)$ in (\ref{FMs}), so that there are $2n_J(r_b)+1$ tones in total, one can derive conditions for each type of TFR behavior in the case of FM component (\ref{FMs}), that are summarized in Table \ref{tab:FMc}. All conditions are devised in a similar way to that used previously for two tones and for the AM component. Thus, in the case of the FM component (\ref{FMs}), although hard to prove rigorously, it appears that the tones are most merged and most separated in the TFR when $\phi_b(t)=0,\pi$ and $\phi_b(t)=\pi/2,3\pi/2$, respectively. For these cases the TFR amplitudes (\ref{FMtfr}) become
\begin{equation}\label{FMms}
\begin{aligned}
\left.|H_s(\omega,t)|^2\right|_{\phi_b(t)=0,\pi}&=\frac{A^2}{4}\Big[\hat{h}^{(+)}\pm \hat{h}^{(-)}\Big]^2,\\
\left.|H_s(\omega,t)|^2\right|_{\phi_b(t)=\pi/2,3\pi/2}&=\frac{A^2}{4}\Big[|\hat{y}^{(+)}(\omega)|^2+|\hat{y}^{(-)}(\omega)|^2\Big],\\
\end{aligned}
\end{equation}
where
\begin{equation}\label{FMhtilde}
\begin{aligned}
\hat{y}^{(+)}(\omega)&\equiv J_0(r_b)\hat{h}_\nu(\omega)+\sum_{n=1}^\infty (-1)^nJ_{2n}(\omega)\Big[\hat{h}_{\nu+2n\nu_b}(\omega)+\hat{h}_{\nu-2n\nu_b}(\omega)\Big],\\
\hat{y}^{(-)}(\omega)&\equiv \sum_{n=1}^\infty (-1)^nJ_{2n-1}(\omega)\Big[\hat{h}_{\nu+(2n-1)\nu_b}(\omega)+\hat{h}_{\nu-(2n-1)\nu_b}(\omega)\Big].
\end{aligned}
\end{equation}
Applying to (\ref{FMms}) the same logic as used previously (see Sec.\ \ref{sec:OVA}) gives the conditions for Regimes II and III, while conditions for the Regimes I and IV are defined in the usual way. Note that, in contrast to previously considered signals, for the FM component (\ref{FMs}) it is not rigorously proven that for Regime III the requirement $\eta_b< 1-\epsilon$ assures $\langle N_p\rangle>1$, but this seems to be the case in practice.

In what follows, we consider $r_b\in[0,1]$, so that one can restrict the consideration to only five tones $\sim J_{0,\pm1,\pm2}(r_b)$ in (\ref{FMs}), since $J_{n\geq3}(r_b\in[0,1])$ is negligible. In this case all is quite straightforward, i.e.\ there is one dominant tone in (\ref{FMs}), corresponding to the main frequency $\nu$, and few FM-induced tones of smaller amplitudes. For higher values, e.g.\ $r_b\gtrsim 1.45$, $J_1(r_b)$ becomes higher than $J_0(r_b)$ so, instead of one main tone, one has two dominant side tones, as well as additional FM-induced tones to be considered (as $J_{n>2}(r_b)$ is no longer negligible). With further increase in $r_b$, the FM-induced tones with frequencies more distant from $\nu$ become dominant, and additional terms become non-negligible. However, $\eta_b$ as defined in (\ref{FMeta}) allows for appropriate discrimination between Regimes for any $r_b$, and all considerations apply in the general case. We consider below a simple five-tone case with $r_b\in[0,1]$.

\begin{figure*}[t!]
\includegraphics[width=1.0\linewidth]{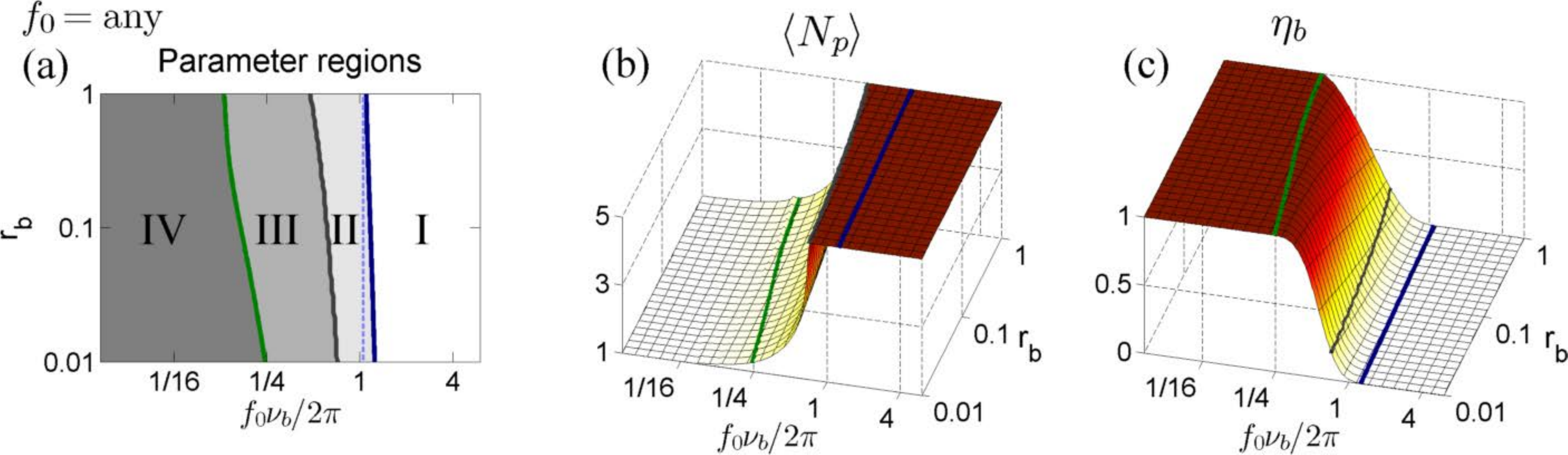}\\
\caption{Dependence of WFT behavior for the signal $s(t)={\rm Re}\sum_{n=-2}^2 J_n(r_b)e^{i(\nu+n\nu_b) t}\approx \cos(\nu t+r_b\sin \nu_b t)$ on parameters $r_b,\nu_b$ and Gaussian window resolution parameter $f_0$. (a): Regions of parameter space corresponding to each type of behavior, according to (\ref{FMt1}),(\ref{FMt2}),(\ref{FMt3}),(\ref{FMt4}); dashed light-blue line shows boundary of the I-type behavior as predicted by approximate (\ref{FMt1a}). (b): Mean number of peaks $\langle N_p\rangle$ (\ref{NP}). (c): Relative overlap $\eta_b$ (\ref{FMeta}). For determining I and IV Regimes we used $\epsilon=0.001$ in (\ref{FMt1}) and (\ref{FMt4}); for all $f_0$ we assume that $2\nu_b<\nu$.}
\label{FMbdWFT}
\end{figure*}

For the Gaussian-window WFT of a signal (\ref{FMs}), where only terms $\sim J_{0,\pm1,\pm2}(r_b)$ are taken into account, the parameter regions of each behavior and the corresponding values of $\langle N_p\rangle$ (\ref{NP}) and $\eta_b$ (\ref{FMeta}) are shown in Fig.\ \ref{FMbdWFT}. As usual, all measures depend only on $r_a$ and $f_0\nu_b$. In general, everything is quite similar to the case of the amplitude-modulated signal (Fig.\ \ref{AMbdWFT}), although here $r_b=1$ is not the maximum allowed value, and the region of Regime IV behavior does not squeeze down to zero anywhere.

\begin{figure*}[t!]
\includegraphics[width=1.0\linewidth]{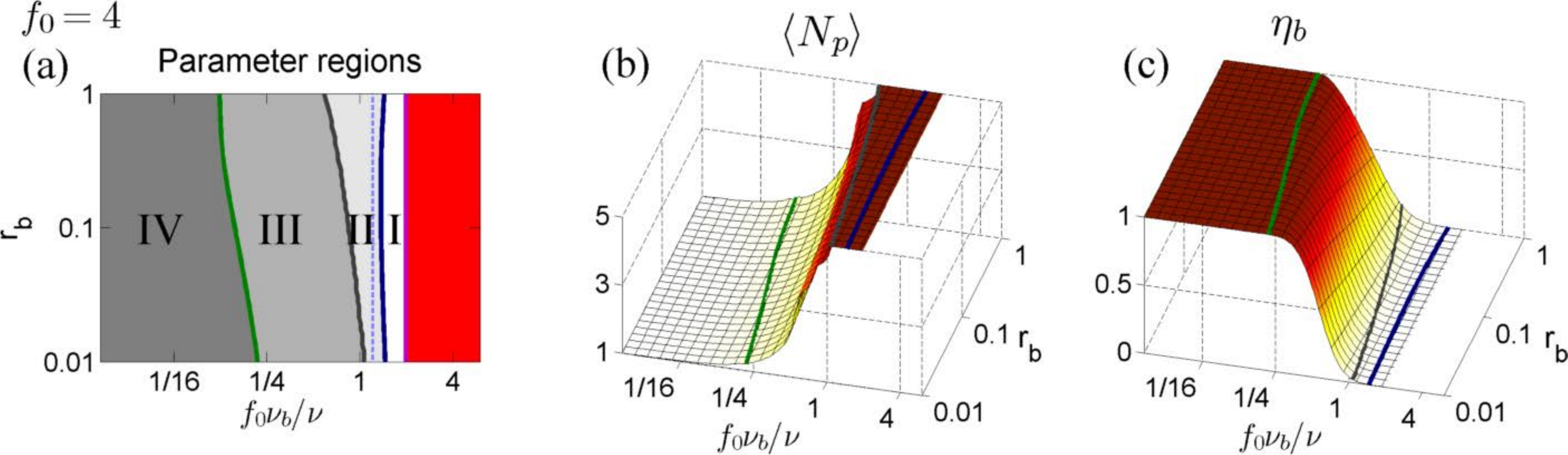}\\
\includegraphics[width=1.0\linewidth]{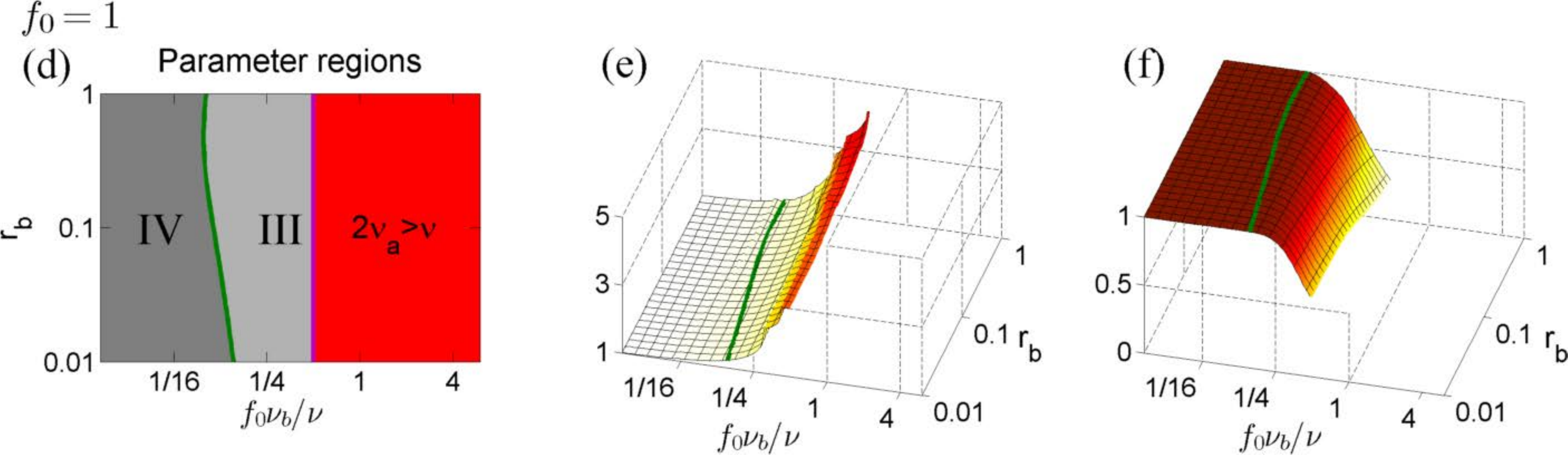}\\
\caption{Same as Figure \ref{FMbdWFT}, but for the Morlet wavelet WT. In (a,d), red filling indicate an inappropriate regions where $2\nu_b>\nu$; in (b) and (c) such regions are omitted (as $\eta_b$ (\ref{FMeta}) is not well-defined for them).}
\label{FMbdWT}
\end{figure*}

Fig.\ \ref{FMbdWT} represents an analog of Fig.\ \ref{FMbdWFT}, but for the WT with the Morlet wavelet. Similarly to the case of amplitude modulation (see Fig.\ \ref{AMbdWT}), in $f_0\nu_b/\nu$ we cannot exceed some threshold, in our case determined as $f_0/n_J(r_b)=f_0/2$: for appropriate analysis we need all non-negligible tones in (\ref{FMs}) be located on the positive frequency axis ($\nu>n_J(r_b)\nu_a=2\nu_a$), so that the analytic signal amplitude and phase approximation are valid (see Part I). The disallowed parameter region is shown in red in Fig.\ \ref{FMbdWT}. The usual contrast with the WFT is that the dependence of the WT behavior and the characterization of its parameters in terms of $r_b$ and $f_0\nu_b/\nu$ differs for different $f_0$ and, in the case considered, this feature appears to be a general property characteristic of any wavelet.

\subsubsection{Reconstruction}

We now investigate the performance of different reconstruction methods for the frequency-modulated signal (\ref{FMs}). In all TFRs, we use ``maximum-based'' curve extraction, selecting TFS around the highest peak in the TFR amplitude at each time. Reconstruction methods are then applied, the resultant signals are compared with their true values, and the errors $\varepsilon_{a,\phi,f}$ (\ref{apfrel}) are calculated. It should be noted, that ``frequency-based'' TFS extraction (see Sec.\ \ref{sec:assumptions}) is not appropriate here, because it will introduce bias and lead to the extracted TFS being quite different from what one encounters in real life (e.g.\ selecting different tones when the TFR behavior is within Regime)I; for this reason, it will be not discussed further.

\begin{figure*}[t!]
\includegraphics[width=1.0\linewidth]{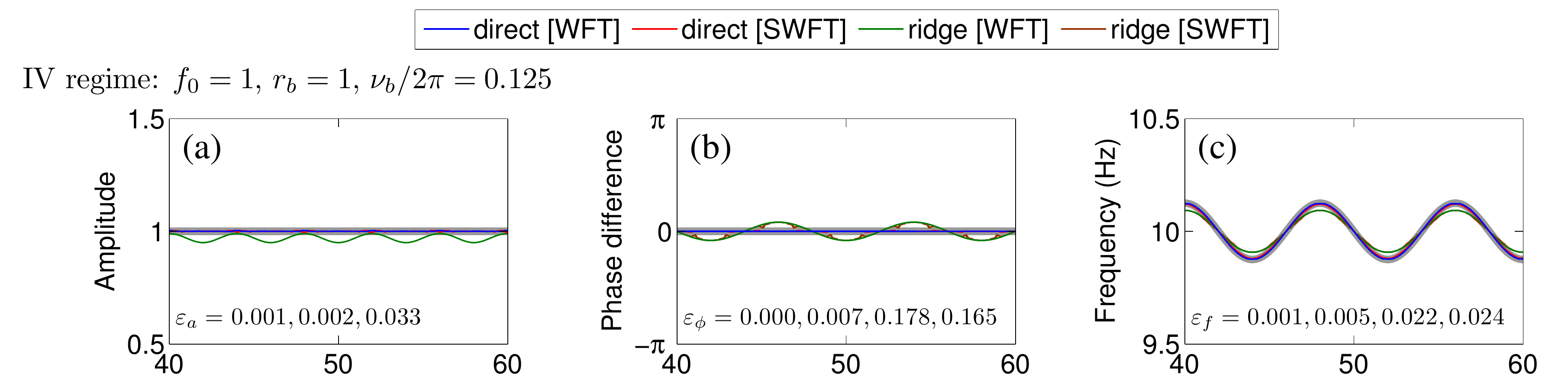}\\
\includegraphics[width=1.0\linewidth]{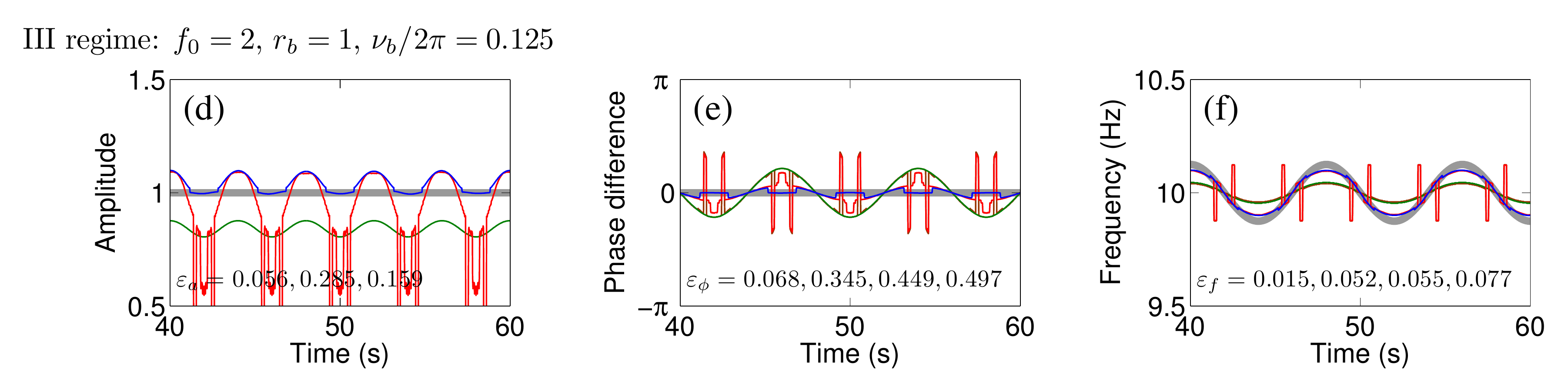}\\
\caption{Amplitude (a,d), phase (b,e) and frequency (c,f) of the FM component (\ref{FMs}) as reconstructed from its WFT and SWFT (colored lines), compared to the true values (thick gray lines). Values of $\varepsilon_{a,\phi,f}$ shown are in the same order as lines in legend, corresponding to ${\rm direct [WFT]}$ (blue), ${\rm direct [SWFT]}$ (red), ${\rm ridge [WFT]}$ (green), ${\rm ridge [SWFT]}$ (brown). In (a,d), ridge reconstruction from the SWFT is not shown as it is not appropriate for amplitude (see Part I). In (b,e), the difference between the reconstructed and true phase is shown. The signal (\ref{FMs}) was sampled at $50$ Hz for $100$ s, and it was simulated with $\varphi=\varphi_b=0$ and $\nu/2\pi=10$; all other parameters are indicated on the figure.}
\label{FMex}
\end{figure*}

Fig.\ \ref{FMex} shows examples of reconstructed amplitudes, phases and frequencies from (S)WFTs in two cases, corresponding to Regimes IV and III (for (S)WT all is similar). In Regime IV (a-c), when (S)WFT behavior is appropriate, all methods work quite well. However, in every case, the direct methods outperform the ridge ones, especially for amplitude estimation. This is because the direct estimates are by definition exact in this case, while ridge estimates possess errors dependent on the strength of the frequency modulation (see Sec.\ \ref{sec:optrec}). Methods based on the WFT and SWFT have same performance for Regime IV. When (S)WFT behavior changes to that of the III type (d-f), the performance of all methods become much worse, but direct methods still greatly outperform the ridge-based ones. However, due to the existence of a few peaks at certain times, the extracted time-frequency support no longer includes all the region wherein the whole FM component is contained, so that direct estimates are no longer exact. It can be seen that e.g.\ the frequency reconstructed by ${\rm direct [WFT]}$ deviates from the real one mainly when the frequency modulation is maximal ($\phi_b=0,\pi$) which, as we now know, is exactly where the additional peaks occur. Nevertheless, the most striking feature is the appearance of considerable differences between reconstruction by WFT or SWFT in Regime III, with the former now being much better than the latter for both direct and ridge methods. Thus, the SWFT-based estimates now have rapid jumps at certain times, corrupting the reconstructed values. This happens because there are time intervals within which SWFT power is redistributed to the side frequencies for some reason, so that the dominant peaks become located away from the true frequency; something similar also occurs for strong amplitude-modulation in Regime III (see discussion of Fig.\ \ref{AMrecWFT}).

\begin{figure*}[t!]
\includegraphics[width=1.0\linewidth]{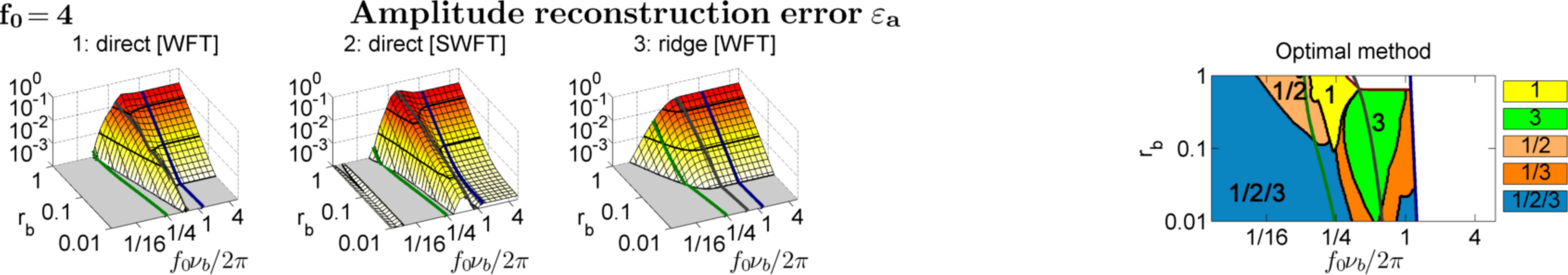}\\
\includegraphics[width=1.0\linewidth]{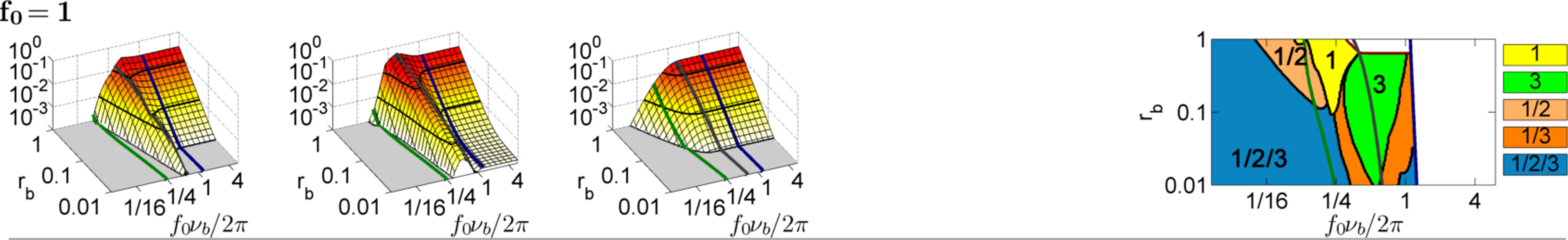}\\
\includegraphics[width=1.0\linewidth]{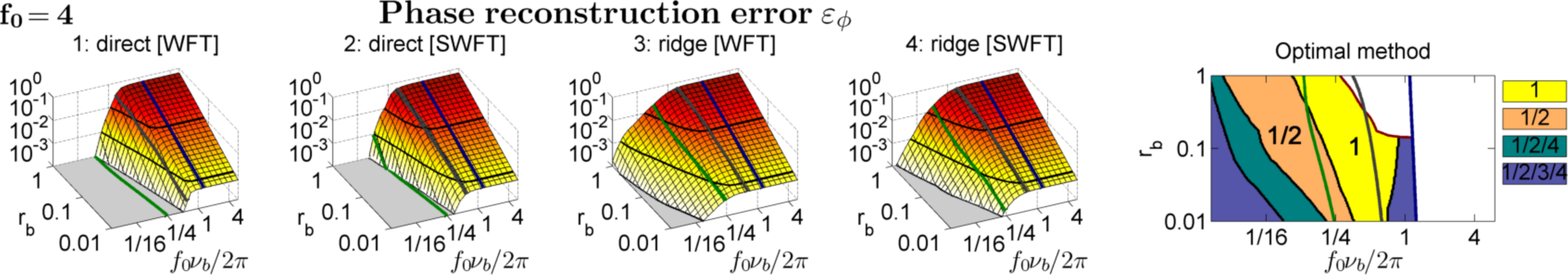}\\
\includegraphics[width=1.0\linewidth]{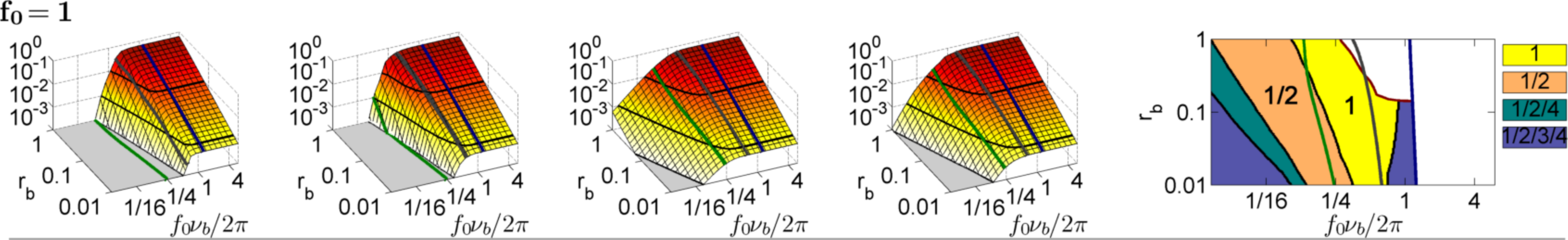}\\
\includegraphics[width=1.0\linewidth]{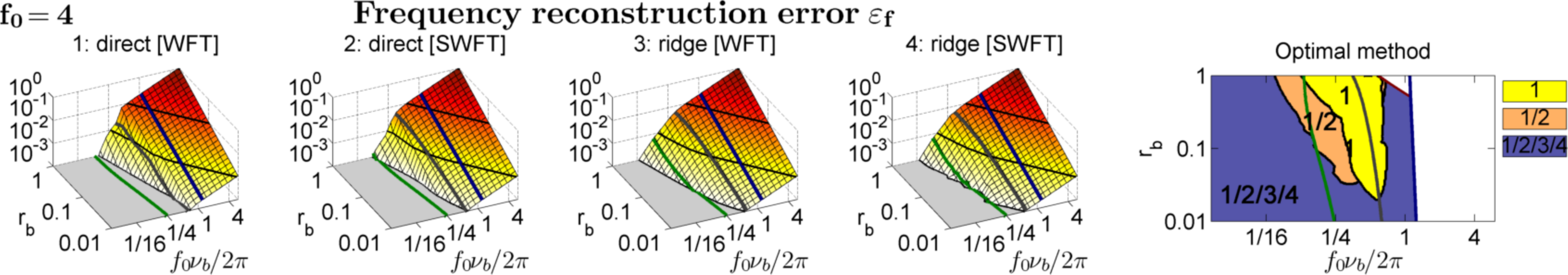}\\
\includegraphics[width=1.0\linewidth]{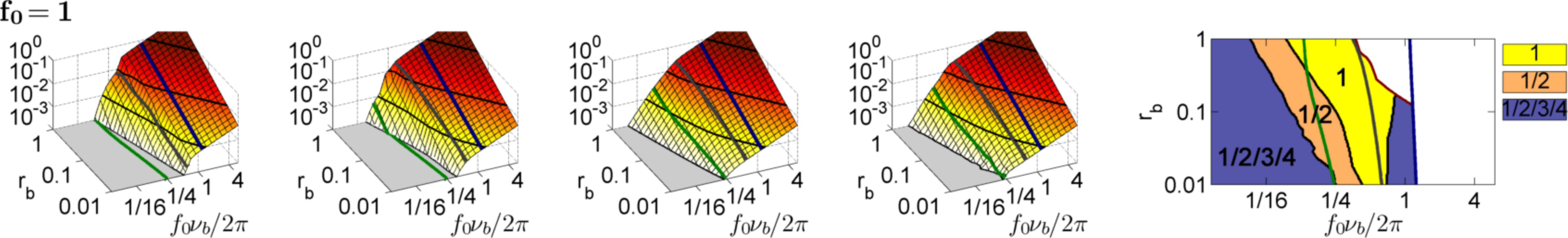}\\
\caption{Left four panels in each row show (S)WFT-based amplitude (1-2 rows), phase (3-4 rows) and frequency (5-6 rows) reconstruction errors for the FM component $s(t)=\cos(\nu t+r_b\sin\nu_a t)$ in dependence on parameters $r_b,\nu_b$ and Gaussian window resolution parameter $f_0$. We used $\nu/2\pi=10$, but results do not depend on its value, at least assuming $2\nu_b<\nu$ (so the analytic signal approximation remains valid). Estimation of amplitude from SWFT ridges is not appropriate, and so the corresponding results are discarded. Thin black lines indicate the levels of $0.001,0.01,0.1,1$, and the behavior of errors $<0.001$ for simplicity is not shown; thick green, gray and blue lines are the same as in Fig. \ref{FMbdWFT}, showing the borders of regions corresponding to behavior in Regimes I-IV. The right panels in each row show regions in parameter space where each method is optimal, i.e.\ gives the smallest estimation error; if the resultant errors of two methods differ on less than $0.001$, they are regarded as having similar performance (the corresponding regions are denoted as `'method1/method2`'). Comparison does not make much sense when WFT behavior is of I type or if all reconstruction errors are $>0.1$, i.e.\ all methods fail; the corresponding regions are white-filled in right panels. The signal was sampled at 50 Hz for 100 s.}
\label{FMrecWFT}
\end{figure*}

\begin{figure*}[t!]
\includegraphics[width=1.0\linewidth]{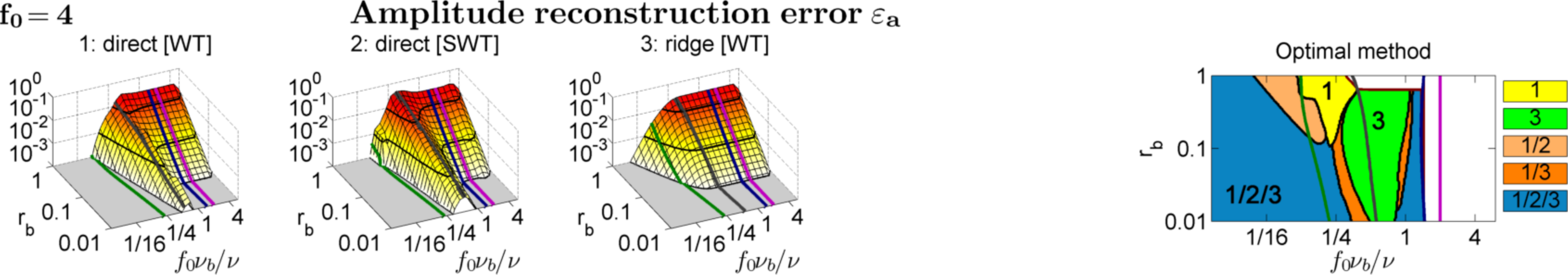}\\
\includegraphics[width=1.0\linewidth]{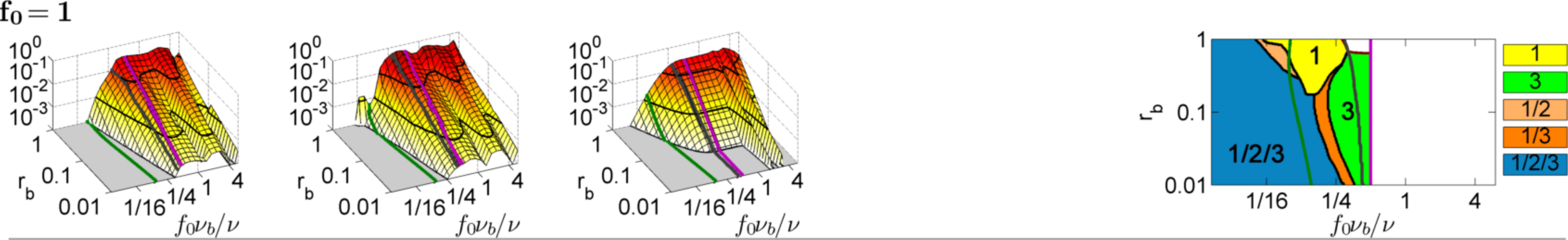}\\
\includegraphics[width=1.0\linewidth]{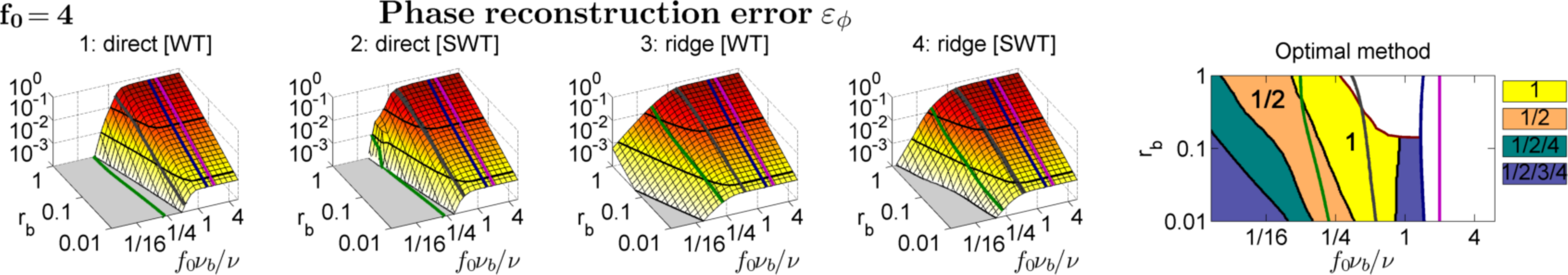}\\
\includegraphics[width=1.0\linewidth]{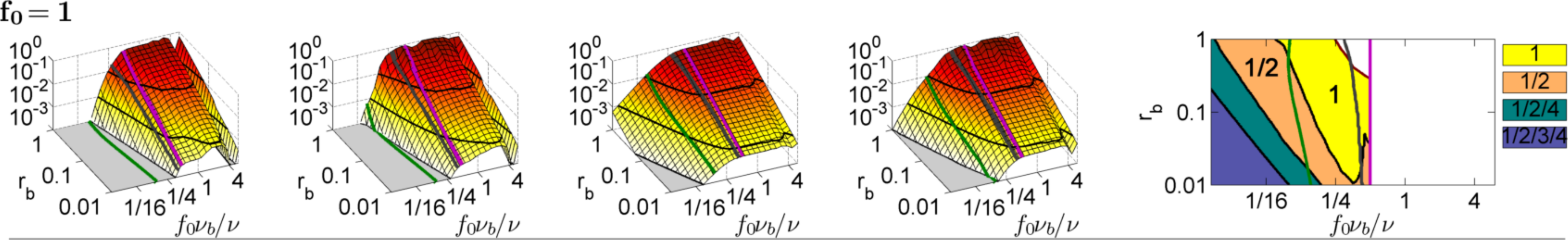}\\
\includegraphics[width=1.0\linewidth]{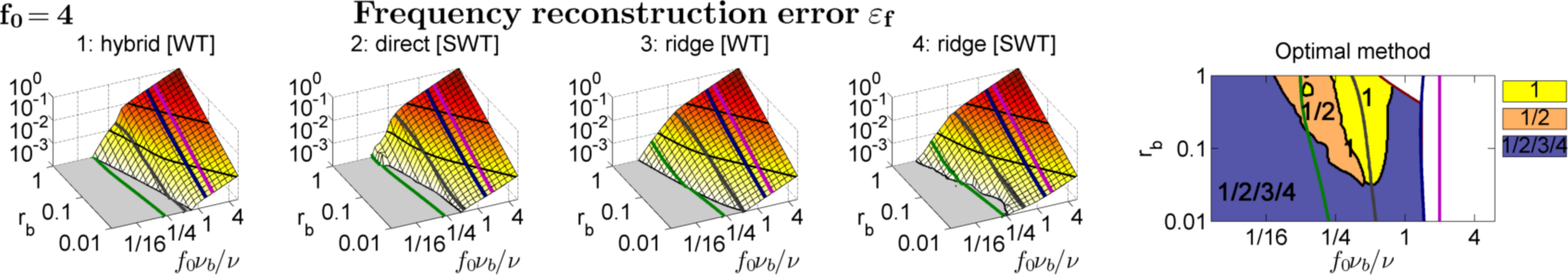}\\
\includegraphics[width=1.0\linewidth]{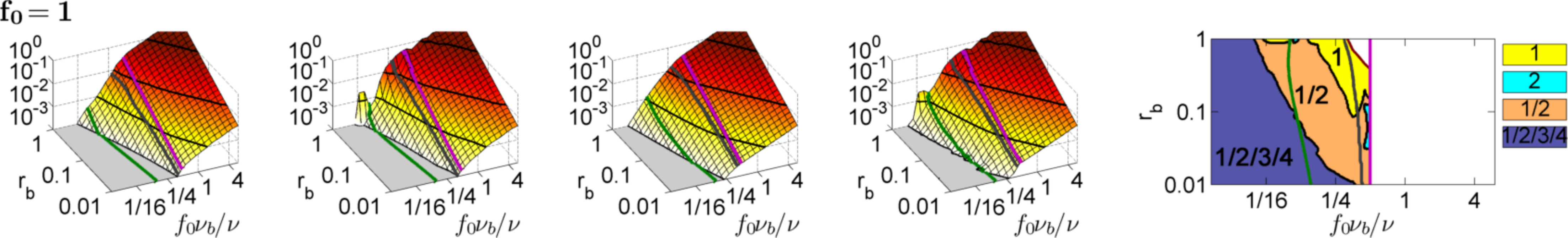}\\
\caption{Same as Fig.\ \ref{FMrecWFT}, but all parameters are reconstructed from the signal's WT/SWT. The direct frequency reconstruction for Morlet wavelet is not possible, so we use hybrid reconstruction (see Part I). Thick magenta lines separate the parameter regions where $\nu_b>\nu/2$, so that the analytic signal approximation is not valid there (implying additional theoretical error which cannot be reduced and thus making any comparisons not appropriate). Note, that the frequency reconstruction error for the WT (\ref{apfrel}) is defined in relation to the component's mean frequency, so $\varepsilon_f$ here is divided by $\nu/2\pi=10$ in comparison to $\varepsilon_f$ in Fig.\ \ref{FMrecWFT}; this makes the dependences depicted consistent for different frequencies $\nu$ in (\ref{FMs}).}
\label{FMrecWT}
\end{figure*}

The dependences on the signal and window/wavelet parameters for the reconstruction errors arising in different methods are shown in Fig.\ \ref{FMrecWFT} for the (S)WFT and in Fig.\ \ref{FMrecWT} for the (S)WT. As usual, for the (S)WFT the dependence of reconstruction errors on $r_b$ and $f_0\nu_b$ is the same for any value of resolution parameter, while for the (S)WT it changes for different $f_0$. The performance of all methods depends largely on the TFR behavior type, with the best estimates being obtained within Regime IV. Similarly to the case of the AM component (and in contrast to the case of two tones), for the FM component (\ref{FMs}) the direct methods outperform the ridge-based ones for the estimation of all three characteristics, giving almost exact estimates when the TFR behavior is of the Regime IV type (so that all tones in (\ref{FMs}) are attributed to the same TFS). Synchrosqueezing does not improve the performance of the direct methods (in fact, it can only worsen it) but, in contrast to the previous examples, it seems to reduce slightly the error of the ridge estimation of phase (but not those of amplitude or frequency) within Regime IV.

\subsection{The effect of noise}\label{sec:WN}

We now turn to the effect of noise on the different TFRs. Consider a signal
\begin{equation}\label{WNs}
s(t)=\cos\nu t+\frac{\sigma}{\sqrt{2}}\zeta(t)
\end{equation}
where $\zeta(t)$ denotes white Gaussian noise of unit variance, and $\sigma/\sqrt{2}$ is its standard deviation. The $1/\sqrt{2}$ multiplier is introduced to make $\sigma$ equivalent to noise-to-signal ratio (standard deviation of the noise divided by that of the signal).

\subsubsection{Representation}

In general, any noise has its power continuously distributed over the whole or large part of a Fourier domain, and so can be regarded as a superposition of infinitely many tones with random independent phase shifts (for white noise they all have the same amplitude). Thus, in some sense the case of noise is similar to the case of interfering components; but, for noise, we have (in theory) infinitely many tones with infinitely close frequencies. For the case of (\ref{WNs}) it is hard to devise a helpful classification of TFR behavior, so we restrict ourselves to qualitative and illustrative considerations.

\begin{figure*}[t!]
\includegraphics[width=1.0\linewidth]{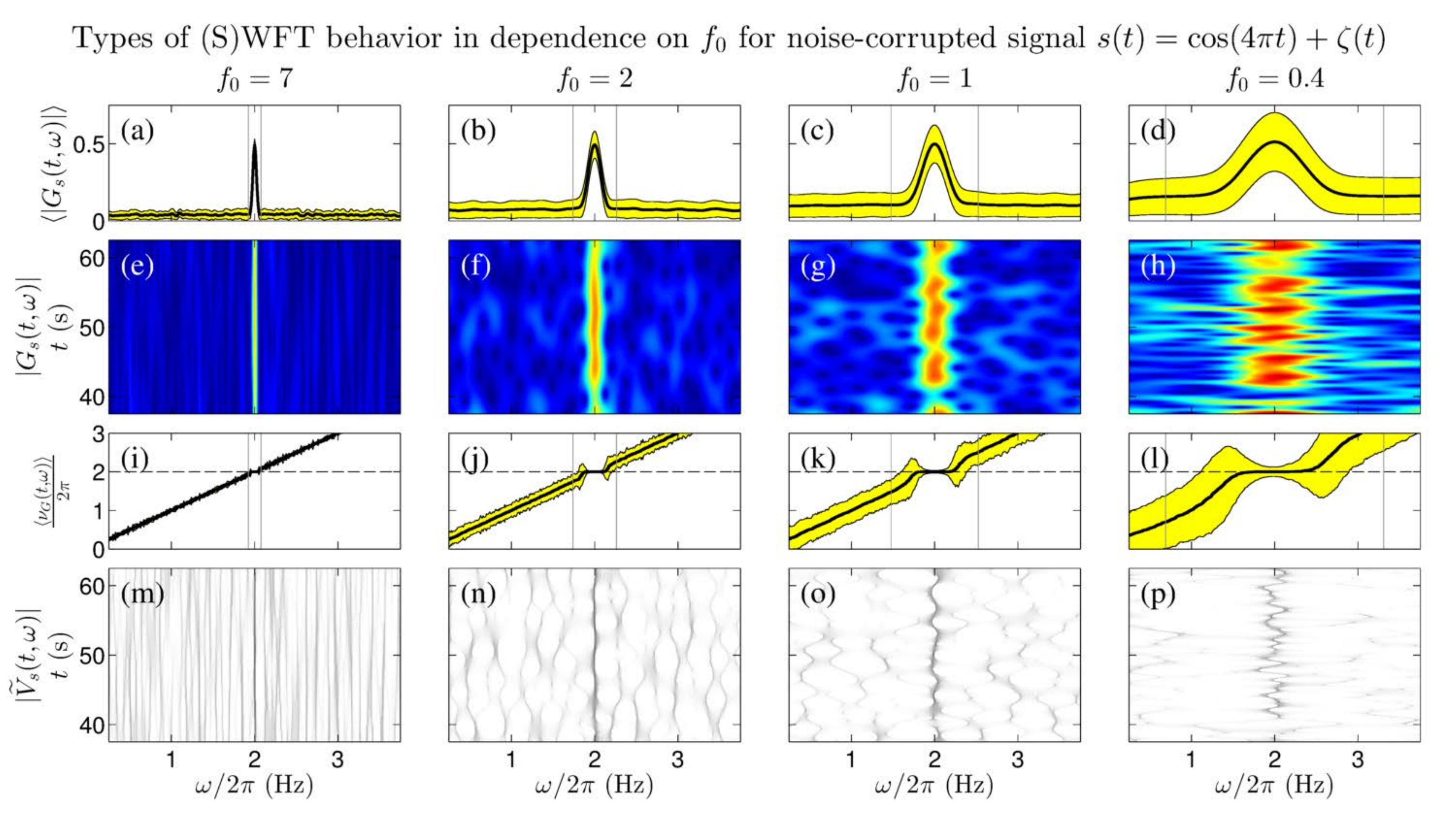}\\
\caption{Behavior of the WFT in dependence on $f_0$ for signal consisting of tone corrupted by white noise of unit deviation $s(t)=\cos(4\pi t)+\zeta(t)$, sampled at $50$ Hz for $100$ s. (a-d): Time-averaged WFT amplitudes. (e-h) WFT amplitudes in time-frequency domain. (i-l): Time-averaged WFT frequency $\nu_G(\omega,t)$; dashed line shows the tone's frequency $\nu/2\pi=2$. (m-p): SWFT amplitudes in time-frequency domain.}
\label{WNbt}
\end{figure*}

The (S)WFT of the signal (\ref{WNs}) is presented in Fig.\ \ref{WNbt} for different values of $f_0$. We observe that, the higher the $f_0$, the better we can distinguish genuine tone within the noise. This is because the more the spread of $\hat{g}(\xi)$ or $\hat{\psi}(\xi)$ is (which is inversely proportional to $f_0$), the more noise tones are picked up while calculating the TFR. Thus, one can note that even the mean noise amplitude in the WFT increases with decreasing $f_0$, indicating that the noise contribution to each frequency bin has become stronger. Note that, for this reason, in the WT amplitude $|W_s(\omega,t)|$ the noise intensity increases with $\omega$ (not shown) as a result of the logarithmic frequency resolution, while white noise has constant power on a linear frequency scale.

We conclude that, for the single tone corrupted by noise (\ref{WNs}) the higher $f_0$ is the better (in fact, the best one can use here is the usual Fourier transform estimate, which provides maximum possible frequency resolution). However, if the signal represents a noise-corrupted component with amplitude or frequency variation, then one needs to choose $f_0$ as a compromise between reducing the effect of noise and ensuring that the AM/FM component is still represented as a single entity in the TFR; for the previously considered cases (\ref{AMs}) and (\ref{FMs}) this will correspond to selecting $f_0$ at the border of Regime IV, determined by (\ref{AMt4}) and (\ref{FMt4}) with $\epsilon$ proportional to the noise power.

\subsubsection{Reconstruction}

We now investigate the performance of different reconstruction methods for a single tone corrupted by noise (\ref{WNs}). In all TFRs, we extract the time-frequency support around the maximum TFR amplitude at each time (maximum-based scheme). For the WT, however, such a scheme should be slightly modified since, as mentioned, the white noise WT amplitude increases with frequency so that, given a wide enough frequency range, the highest peak in the WT amplitude will always appear at the highest frequencies and will correspond to noise. To avoid such issues for the WT, we extract the TFS around the points corresponding to the highest peaks in $|W_s(\omega,t)/\sqrt{\omega}|$ (since the white noise level in the WT is $\sim\sqrt{\omega}$). Note, that the mentioned issues do not arise for the SWT. To study the general robustness to noise, which may be colored or non-Gaussian, we do not apply any noise filtering (e.g.\ WFT/WT hard or soft thresholding). Hence our investigation is not limited to white noise, and remains valid even if we replace $\zeta(t)$ in (\ref{WNs}) with any other (colored) noise having FT power at frequencies in the vicinity of $\nu$ equal to the power of unit-deviation white noise.

\begin{figure*}[t!]
\includegraphics[width=1.0\linewidth]{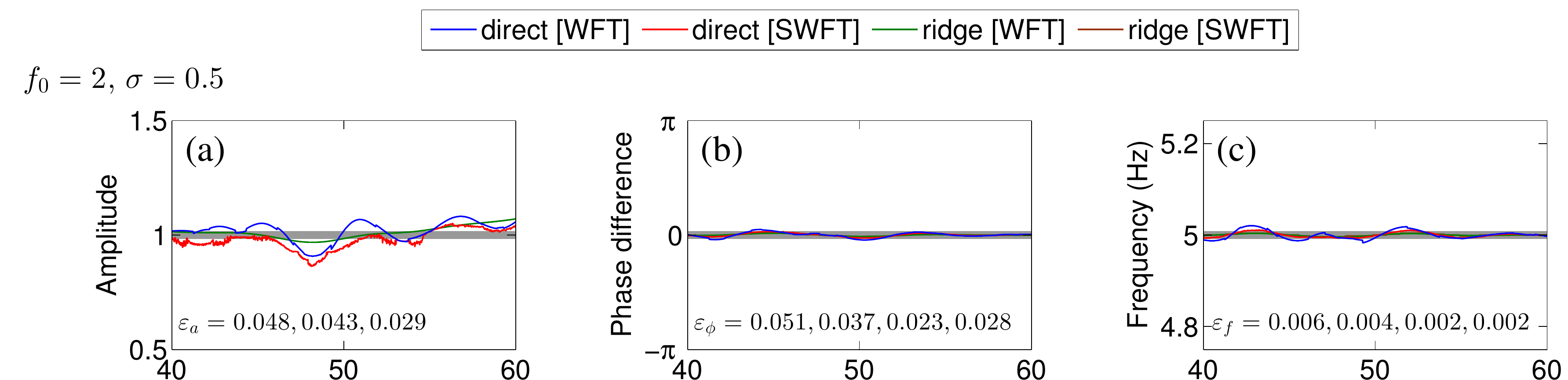}\\
\includegraphics[width=1.0\linewidth]{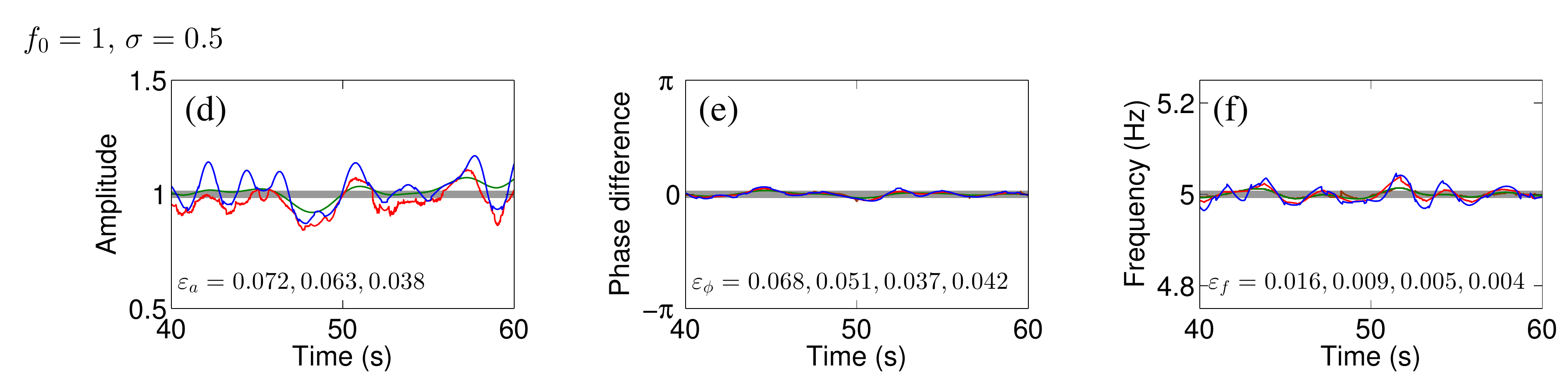}\\
\caption{Amplitude (a,d), phase (b,e) and frequency (c,f) of the noise-corrupted tone (\ref{WNs}) as reconstructed from its WFT and SWFT (colored lines), compared to the true values (thick gray lines). Values of $\varepsilon_{a,\phi,f}$ shown are in the same order as lines in legend, corresponding to ${\rm direct [WFT]}$ (blue), ${\rm direct [SWFT]}$ (red), ${\rm ridge [WFT]}$ (green), ${\rm ridge [SWFT]}$ (brown). In (a,d), ridge reconstruction from the SWFT is not shown as it is not appropriate for amplitude (see Part I). In (b,e), the difference between the reconstructed and true phase is shown. The signal was sampled at $50$ Hz for $100$ s.}
\label{WNex}
\end{figure*}

Fig.\ \ref{WNex} shows examples of amplitudes, phases and frequencies reconstructed from a Gaussian window (S)WFT at different $f_0$ (for (S)WT all is qualitatively the same). Comparing (a-c) and (d-f), one can see that the noise contribution increases with decreasing $f_0$, consistent with what was seen before in Fig.\ \ref{WNbt}. Clearly, the ridge methods are much more noise-robust than the direct methods, which is to be expected given their better performance for interfering tones (see Sec.\ \ref{sec:OV}). Thus, while ridge reconstruction accounts for noise contribution at only one frequency (corresponding to the TFR amplitude peak), in the direct methods one integrates over all frequencies which constitute the current TFS, thus picking up noise contributions from a much wider frequency band. Although all becomes more complicated in the case of SWFT/SWT, from the results presented it seems that the same reasoning applies. Note, that phase reconstruction by both methods is probably the most noise-robust among the three characteristics, at least as it appears visually (Fig.\ \ref{WNex}(b,e)).

\begin{figure*}[t!]
\includegraphics[width=1.0\linewidth]{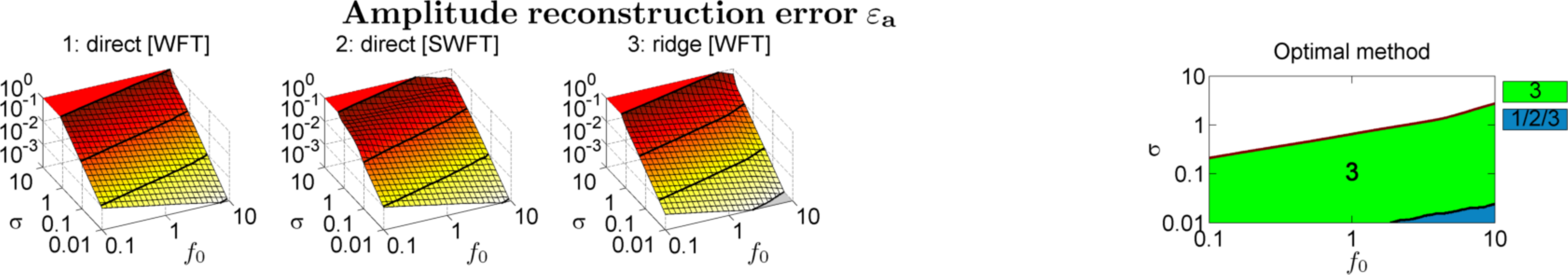}\\
\includegraphics[width=1.0\linewidth]{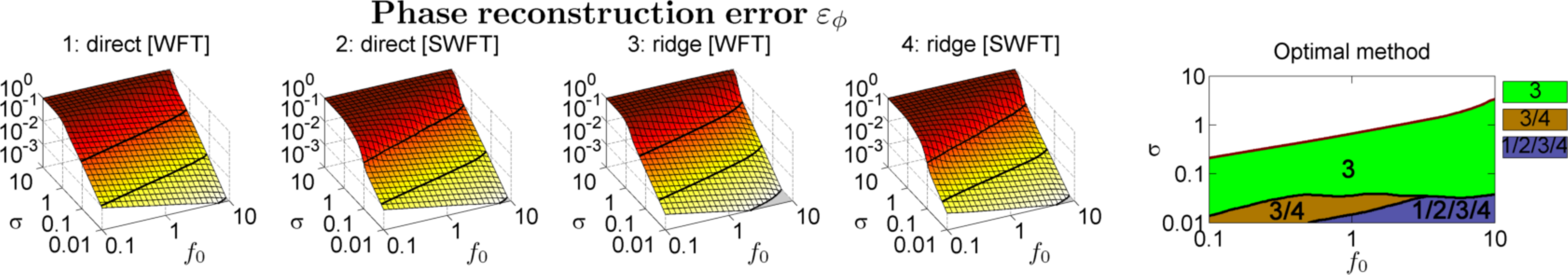}\\
\includegraphics[width=1.0\linewidth]{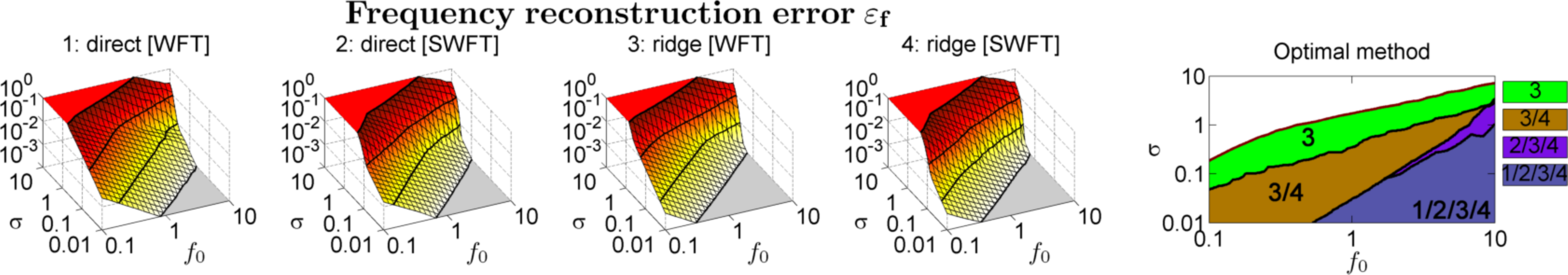}\\
\caption{Left four panels in each row show (S)WFT-based amplitude (1-2 rows), phase (3-4 rows) and frequency (5-6 rows) reconstruction errors for the noise-corrupted tone $s(t)=\cos\nu t+(\sigma/\sqrt{2})\zeta(t)$ in dependence on the noise-to-signal ratio $\sigma$ and Gaussian window resolution parameter $f_0$. For each $f_0$ and $\sigma$, the errors shown are averages over ten noise realizations. We used $\nu/2\pi=5$, but results do not depend on its value. Estimation of amplitude from SWFT ridges is not appropriate, and so the corresponding results are discarded. Thin black lines indicate the levels of $0.001,0.01,0.1,1$, and the behavior of errors $<0.001$ and $>1$ for simplicity is not shown. The right panels in each row show regions in parameter space where each method is optimal, i.e.\ gives the smallest estimation error; if the resultant errors of two methods differ on less than $0.001$, they are regarded as having similar performance (the corresponding regions are denoted as `'method1/method2`'). Comparison does not make much sense when all reconstruction errors are $>0.1$, i.e.\ all methods fail; the corresponding regions are white-filled in the right-hand panels. The signal was sampled at 50 Hz for 100 s.}
\label{WNrecWFT}
\end{figure*}

\begin{figure*}[t!]
\includegraphics[width=1.0\linewidth]{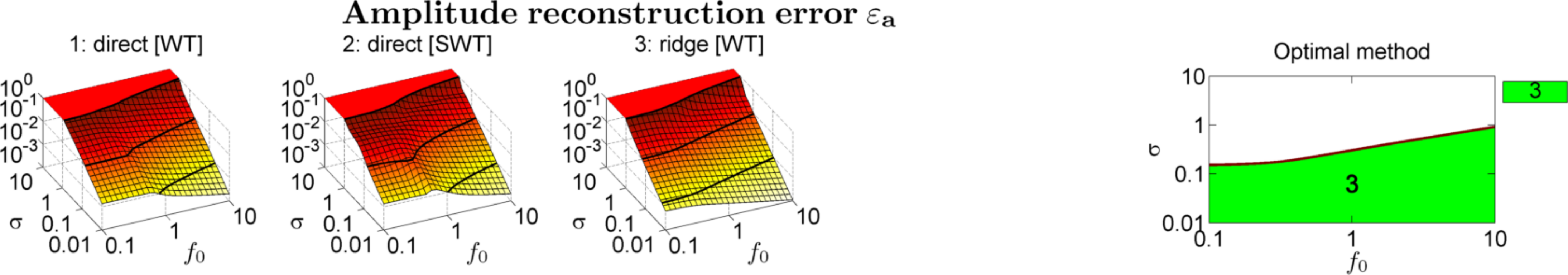}\\
\includegraphics[width=1.0\linewidth]{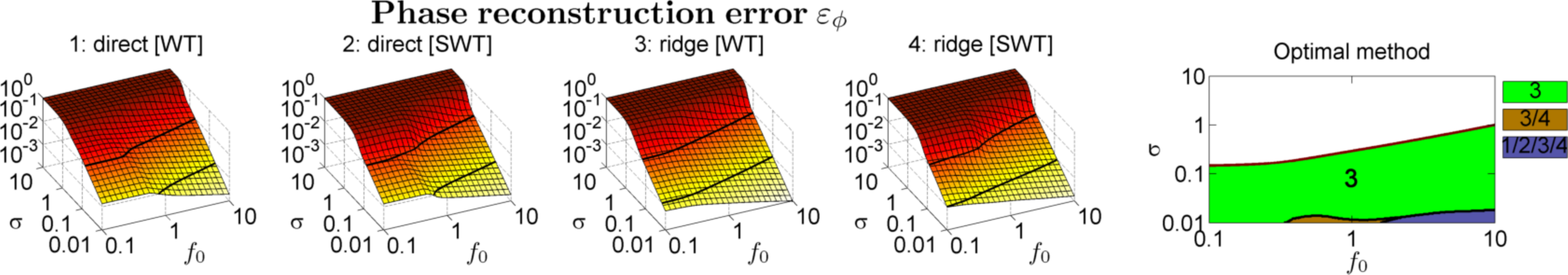}\\
\includegraphics[width=1.0\linewidth]{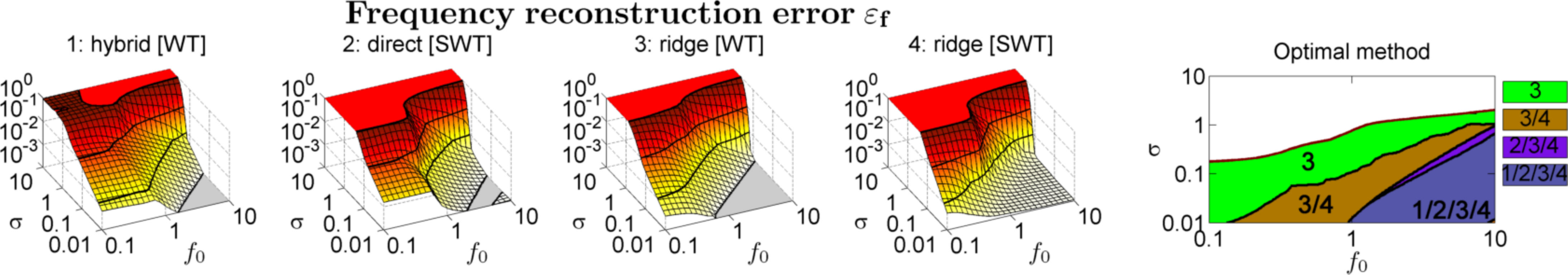}\\
\caption{Same as Fig.\ \ref{WNrecWFT}, but all parameters are reconstructed from the signal's WT/SWT. The direct frequency reconstruction for Morlet wavelet is not possible, so we use hybrid reconstruction (see Part I). Note that the frequency reconstruction error for the WT (\ref{apfrel}) is defined in relation to the component's mean frequency, so $\varepsilon_f$ here is divided on $\nu/2\pi=5$ in comparison to $\varepsilon_f$ in Fig.\ \ref{WNrecWFT}; this makes the dependences depicted consistent for different $\nu$.}
\label{WNrecWT}
\end{figure*}

Figures \ref{WNrecWFT} (WFT/SWFT) and \ref{WNrecWT} (WT/SWT) show the dependences of the reconstruction errors on noise intensity and the resolution parameter $f_0$ for the different methods. As can be seen, for both the (S)WFT- and (S)WT-based reconstructions, ridge methods are superior to direct ones. Interestingly, for (S)WFT with the direct method, $\varepsilon_a$ and $\varepsilon_f$ are $\approx\sigma$ for $f_0\approx0.1$ and $\approx\sigma/10$ for $f_0\approx10$. For any method and characteristic, synchrosqueezing does not provide any significant advantages in terms of reconstruction. Note that, for (S)WT-based reconstruction, the error rapidly and nonlinearly increases for small $f_0$ and then saturates, which is related to the fact that, below some $f_0$, the properties of the Morlet wavelet remain almost the same (see Part I); this does not happen for the lognormal wavelet.

\begin{remark}
The reconstruction errors in the present case will generally depend on the scheme used for the extraction of the tone's TFS. Here we employ the ``maximum-based'' scheme, selecting the time-frequency support around the maximum TFR amplitude peak over the whole frequency range. Such a TFS is closer to what one might encounter in a real situation, but might contain the spurious noise peaks located far from the actual frequency $\nu$ (thus giving rise to frequency discontinuities). This relates especially to the case $\sigma>1$, as can be deduced e.g.\ from $\varepsilon_f$ in Figs.\ \ref{WNrecWFT} and \ref{WNrecWT}. If we instead select the TFS around the TFR amplitude peak nearest to the actual tone frequency $\nu$ (``frequency-based`'' scheme), thus modelling the best possible curve extraction from the WFT/WT, this will reduce frequency reconstruction errors for $\sigma>1$, eliminating their rapid increase after this point and making the dependence of $\varepsilon_f$ on $f_0$ more linear (resembling that for $\varepsilon_a,\varepsilon_\phi$); the WFT/WT-based $\varepsilon_a$ and $\varepsilon_\phi$, on the other hand, will not change significantly. However, for SWFT/SWT frequency-based TFSs, selection will not model the best possible extraction since, due to the synchrosqueezed TFR's non-smoothness, it will lead to the picking up of small nearby noise flashes instead of the main curve, thus introducing some additional complex error into $\varepsilon_a,\varepsilon_\phi$, but reducing the frequency error $\varepsilon_f$ due to the bias of the extraction procedure. In any case, the conclusions remain the same for both schemes.
\end{remark}

\section{Optimal resolution}\label{sec:optres}

Summarizing the results of previous sections, the TFR behavior depends drastically on the choice of window/wavelet parameters, which in turn strongly affects the extraction and reconstruction of the components. Moreover, in the analysis of real data, inappropriate resolution can even lead to wrong conclusions, e.g.\ that there exist a number of oscillations while there is actually only a single oscillation but with amplitude and/or frequency modulation. Thus, the choice of appropriate window/wavelet parameters is crucial for both the interpretation and quantitative estimation of the signal structure.

Unfortunately, there is no universal choice, and the optimal resolution parameter depends on the signal composition. It should be chosen in such a way as to resolve the independent components (Regime I), e.g.\ tones, but at the same time to represent AM/FM components as single curves in the TFR (Regime IV). In other words, each component, with or without amplitude/frequency modulation, should at each time be mapped into a single individual TFS (i.e.\ an area where the TFR is single-peaked and unimodal).

However, usually this cannot be achieved for all components simultaneously, so the optimal choice requires compromise. Furthermore, in real situations, one often does not know for sure whether it is appropriate to represent the underlying process as a single AM/FM component, or as separate tones. This arbitrariness is unavoidable, but the TFR-based methods give the possibility of specifying a criterion by choosing appropriate window/wavelet parameters. Methods that can be used for this task are considered below. Note that, because synchrosqueezing does not change the resolution properties of the transform (see Sec.\ \ref{sec:Synchroneed}), the optimal parameters for the SWFT/SWT will be the same as for the WFT/WT.

\begin{remark}
Although in what follows we consider optimization of the (single) resolution parameter $f_0$, the same methods can be used quite generally for choosing any set of window/wavelet parameters. For example, for some signals it might be advantageous to use a chirped window in the WFT \cite{Stankovic:13,Tao:10,Almeida:94,Capus:03,Katkovnik:95}, e.g.\ $g(t)\sim e^{-t^2/2f_0^2}e^{i\alpha t^2}$, and selection of the optimal pair $\{f_0,\alpha\}$ in this particular case does not differ qualitatively from the selection of $f_0$ alone.
\end{remark}

\begin{remark}
Because the signal can be highly nonstationary, with characteristic amplitude/frequency modulation of the components changing in time, it might be advantageous to use a time-dependent $f_0(t)$ \cite{Jones:94,Emresoy:98,Zhong:10,Hon:12,Pei:12}; in addition, to achieve an accurate representations of the components occupying different frequency bands, one can allow the resolution parameter to depend also on frequency \cite{Jones:90,Jones:97}. The optimal evolution $f_0(t)$ or surface $f_0(\omega,t)$ can be found in the same way as the global $f_0$, but localizing the corresponding criterion to the neighborhood of each time or time-frequency coordinate. However, apart from the significantly higher computational cost of such approaches, selection of even a single resolution parameter for the whole signal represents a challenging problem that still remains unsolved for the general case, as will be seen below. Note also that the direct reconstruction of the component's parameters (see Part I) becomes problematic if $f_0$ depends on $\omega$. In what follows we consider the case of a global $f_0$ only.
\end{remark}

\subsection{``Monocomponent'' approaches}

Optimization of the TFR parameters for the representation of a single AM/FM component (possibly embedded in noise) was considered in \cite{Katkovnik:98a,Emresoy:98,Yin:13}. The main aim of these approaches is usually to estimate the component's instantaneous frequency $\nu(t)=\phi'(t)$ (see \cite{Boashash:03,Boashash:92a,Boashash:92b} for an overview of related concepts and algorithms), which is reconstructed from the optimized TFR. Real signals, however, are rarely monocomponent, and such an assumption is often too restrictive. In fact, these are the multicomponent signals for which time-frequency analysis is most useful (see Part I) and, as will be shown below, the main difficulty of adapting $f_0$ lies in estimating the number of components. We will therefore not consider the monocomponent case but will proceed with a more general approach.


\subsection{Functional approaches}

A particular class of methods for selecting an appropriate $f_0$ based on the signal's structure was considered in \cite{Baraniuk:01,Stankovic:01,Jones:95,Jones:94,Jones:90}. According to these works, the optimal window/wavelet parameters can be selected as being those that minimize a suitably chosen functional of the signal's TFR. In mathematical terms, one chooses an optimal resolution parameter $f_0$, matching the underlying structure of the signal, as
\begin{equation}\label{mfunc}
f_0={\rm argmin}\Big(F[H_s[f_0](\omega,t)]\Big)
\end{equation}
where $F[\cdot]$ is the specified functional of the WFT/WT calculated with the chosen $f_0$.

In practice, having chosen the functional (\ref{mfunc}), an optimal $f_0$ can be found by first calculating $F[...]$ for logarithmically sampled $f_0=f_0^{(\min)},2^{1/\tilde{n}_v} f_0^{(\min)}, 2^{2/\tilde{n}_v} f_0^{(\min)},...,f_0^{(\max)}$ (we use $\tilde{n}_v=2$ unless otherwise specified) and finding the $f_0$ corresponding to the minimum. Then, from that starting point, one can approach the exact minimum iteratively. In the following tests, however, for reasons of computational speed we will not perform this last step, restricting ourselves to the approximate minimum.

Obviously, one should search for optimal resolution parameter in some appropriate region: $f_0\in[f_0^{(\min)},f_0^{(\max)}]$. Indeed, it does not make much sense to consider parameters for which the time length $[0,T]$ (frequency range $[-f_s/2,f_s/2]$) of the signal can contain only limited part (e.g.\ 95\%) of the window/wavelet in time (frequency) domain. This gives the minimal (maximal) values of $f_0$ as
\begin{equation}\label{f0minmax}
\begin{aligned}
\mbox{\textbf{WFT: }}
&f_0^{(\min)}: \xi_2(0.05)[f_0^{(\min)}]-\xi_1(0.05)[f_0^{(\min)}]=2\pi f_s,\\
&f_0^{(\max)}: \tau_2(0.05)[f_0^{(\min)}]-\tau_1(0.05)[f_0^{(\min)}]=T,\\
\mbox{\textbf{WT: }}
&f_0^{(\min)}: \xi_2(0.05)[f_0^{(\min)}]-\xi_1(0.05)[f_0^{(\min)}]=\frac{\omega_\psi}{\omega_{\max}}2\pi f_s,\\
&f_0^{(\max)}: \tau_2(0.05)[f_0^{(\min)}]-\tau_1(0.05)[f_0^{(\min)}]=\frac{\omega_{\min}}{\omega_\psi}T,\\
\end{aligned}
\end{equation}
where $[\omega_{\min},\omega_{\max}]$ is the frequency range for which the TFR is calculated, while $\xi_{1,2}(\epsilon)$ ($\tau_{1,2}(\epsilon)$) denote the $\epsilon$-supports of the window/wavelet in frequency (time), which depend on $f_0$ (see Part I).

Minimization/maximization of almost all functionals proposed so far \cite{Baraniuk:01,Stankovic:01,Jones:95,Jones:94,Jones:90} can be reduced to the minimization of
\begin{equation}\label{fpq}
F_{p,q}[H_s(\omega,t)]=
\log\frac{\left(\int |H_s(\omega,t)|^pd\mu(\omega) dt\right)^{q/p}}{\int |H_s(\omega,t)|^qd\mu(\omega) dt},
\quad q>p>0.
\end{equation}
For example, the Renyi entropy of order $\alpha$, which was thoroughly investigated in \cite{Baraniuk:01,Sucic:11,Saulig:11,Malarvili:07}, corresponds to $p=2,q=2\alpha$ in (\ref{fpq}); the measures proposed in \cite{Stankovic:01} correspond to $p=1,q=2$; while maximization of the ratio of the fourth power of $L_4$-norm to the squared $L_2$-norm proposed in \cite{Jones:95,Jones:94,Jones:90} corresponds to the minimization of $F_{2,4}$. Note that, generally, $p$ and $q$ are not restricted to integers.

In (\ref{fpq}), the numerator power $q/p$ and inequality $q>p$ have clear meanings. Thus, to provide meaningful results, the minimizing functional should not depend on window/wavelet normalization, that can depend on the parameter to be optimized. This is achieved by the $q/p$ power of the numerator in (\ref{fpq}), which therefore cannot be changed to some other, independent value. Next, for the case of a single tone signal $s(t)=\cos(\nu t+\varphi)$ (delta-peak $s(t)=\delta(t-t_0)$) the optimal resolution parameter should minimize the spread of window/wavelet in frequency (time), corresponding to $f_0\rightarrow\infty$ ($f_0\rightarrow0$), so that the tone (delta-peak) is maximally localized in the TFR. Therefore, the minimum of $F_{p,q}$ should appear at the corresponding asymptotic values for these cases. One can show that, for a Gaussian window, the functional is $F_{p,q}\sim f_0^{1-q/p}$ for the WFT of a single tone signal, and (approximately) $F_{p,q}\sim f_0^{q/p-1}$ for the WFT of a delta-peak. Hence, to establish the required minima at $f_0\rightarrow\infty$ and $f_0\rightarrow0$ one should have $q/p>1$, i.e.\ $F_{p,q}$ should represent the ratio of the lower-order norm to the higher-order norm. Although hard to prove, it seems to be a general rule, valid for WFT/WT with any window/wavelet functions.

\begin{figure*}[t!]
\includegraphics[width=0.7\linewidth]{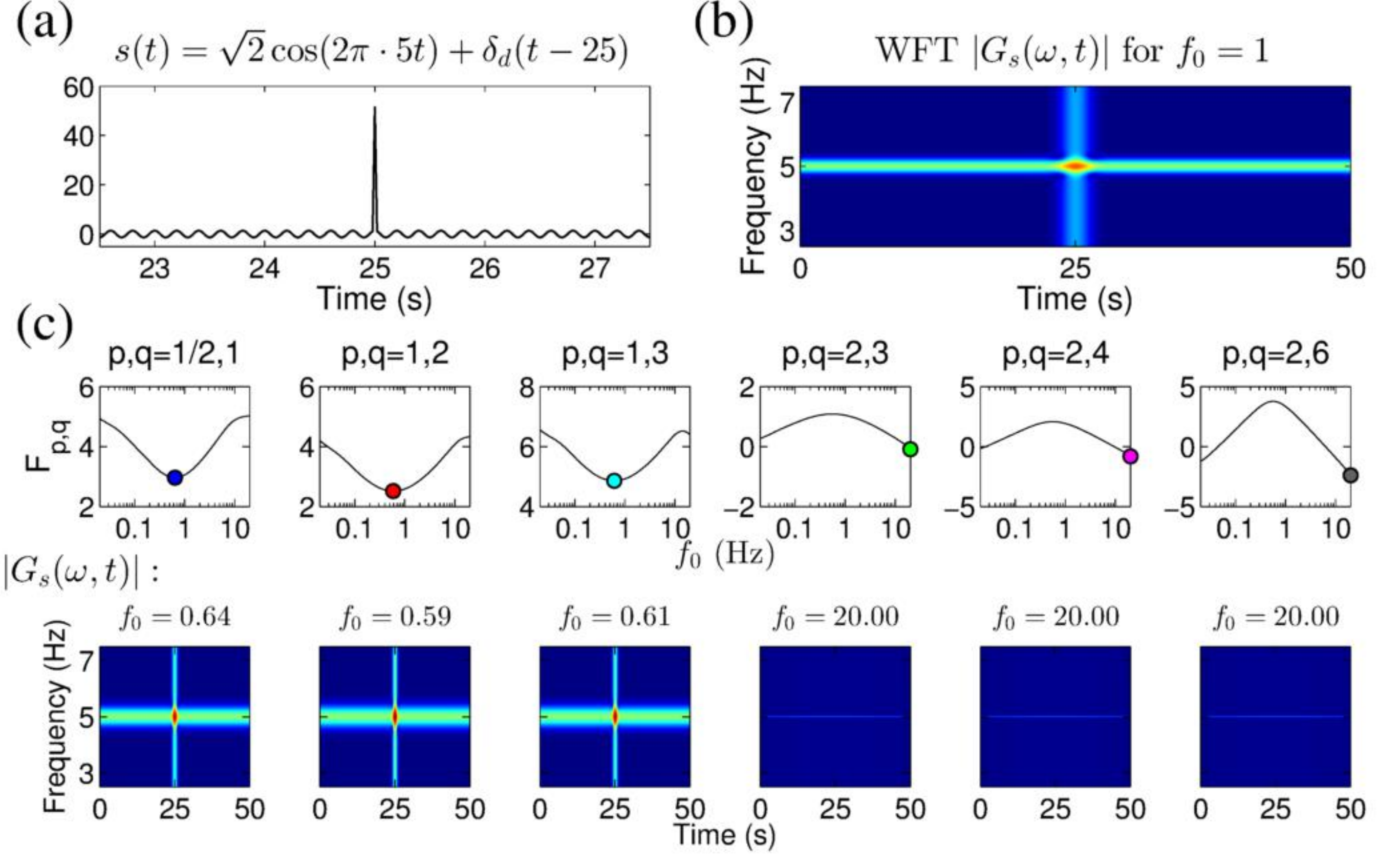}
\includegraphics[width=0.3\linewidth]{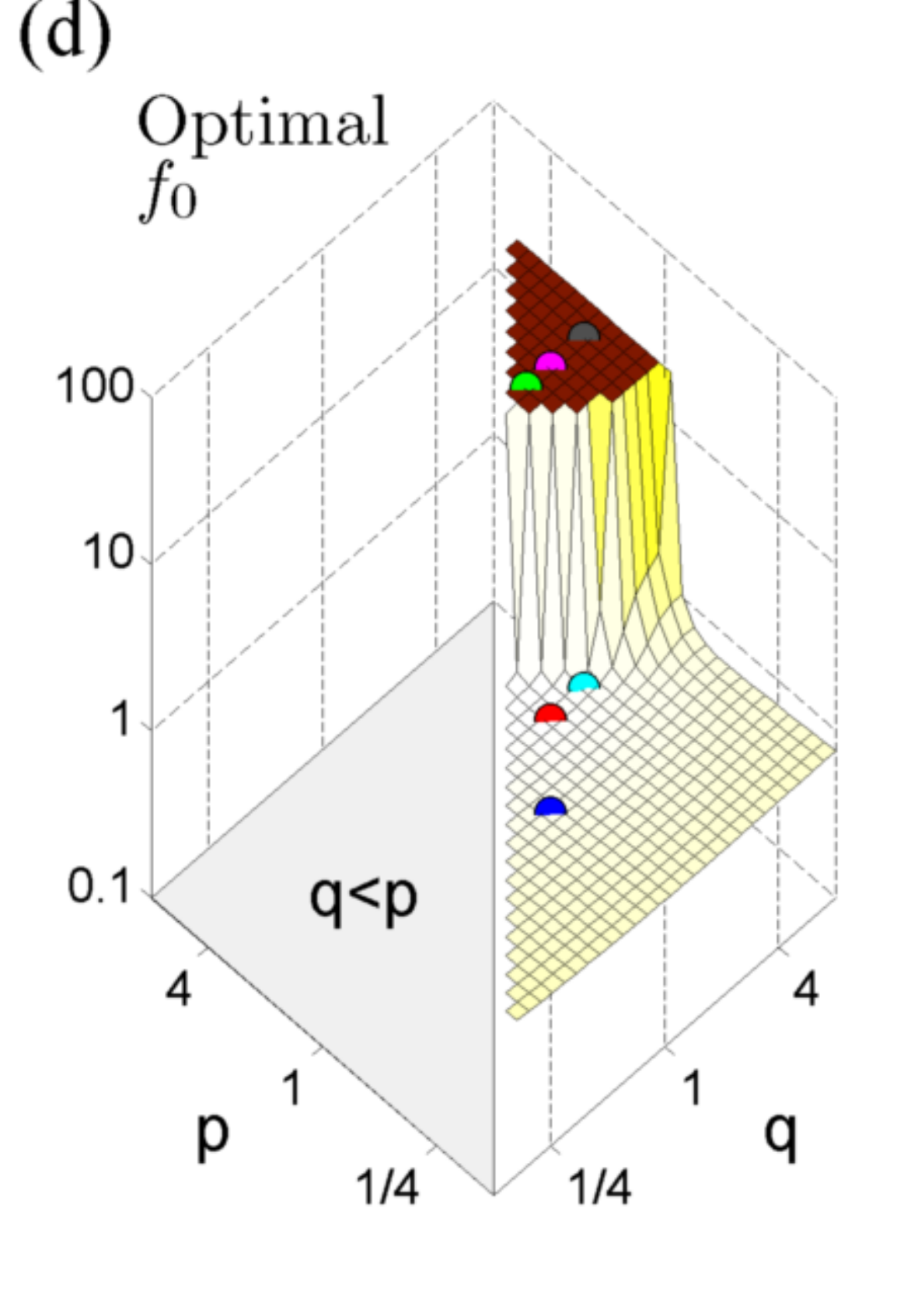}\\
\caption{Adapting the WFT according to (\ref{fpq}) for the signal $s(t)=\sqrt{2}\cos(2\pi \cdot 5t)+\delta_d(t-25)$, where $\delta_d(t-t_0)$ is the discrete analog of $\delta$-function, constructed by adding to $s(t_0)$ the square root of the total signal length $\sqrt{N}$; with such normalization, $\delta_d(t-t_0)$ has the same energy $\int |s(t)|^2dt$ as $\sqrt{2}\cos\nu t$. (a): Signal in the time domain. (b): Example of the signal's WFT calculated for $f_0=1$. (c): Examples of the dependencies of $F_{p,q}$ (\ref{fpq}) on $f_0$ for different pairs of $p,q$ (upper row); optimal $f_0$ minimizing $F_{p,q}$ are shown by filled circles, with the WFTs calculated for these values being presented at the bottom. (d): Optimal $f_0$ in dependence on $p$ and $q$ in (\ref{fpq}). The signal was sampled at $50$ Hz for $50$ s. The WFT was always calculated in the frequency range $[0,25]$ Hz, and the optimal $f_0$ was searched in a range $[0.05,20]$.}
\label{TFadapt}
\end{figure*}

To understand the influence of $p$ and $q$ used in (\ref{fpq}) on the results, it is useful to study a signal consisting of a sinusoid and a delta-pulse. As discussed above, the maximum TFR concentration for the former is achieved at $f_0\rightarrow\infty$, while for the latter at $f_0\rightarrow0$; when both are present, $f_0$ must represent some compromise. The results are presented in Fig.\ \ref{TFadapt} for the WFT with a Gaussian window (qualitatively, one observes the same picture for the WFT with other windows and for the WT). Evidently, there are two possibilities, determined mainly by the value of $p$: for $p\lesssim 2$, the minimum of $F_{p,q}$ (\ref{fpq}) occurs in the ``middle'', providing some compromise between the time and frequency resolutions and trying to represent reliably both a sinusoid and a delta-pulse while, for $p\gtrsim 2$, minima occur at the limiting values on both sides ($f_0\rightarrow 0,\infty$), thus achieving the maximally compressed representation of only one among two features (no compromise). Due to these issues, which were also noted in \cite{Stankovic:01}, the measure of \cite{Stankovic:01} ($p=1,q=2$) seems to be more appropriate in general than that of \cite{Jones:95,Jones:94,Jones:90} ($p=2,q=4$) or the Renyi entropy \cite{Baraniuk:01,Sucic:11,Saulig:11,Malarvili:07} ($p=2,q=6$).

\begin{figure*}[t!]
\includegraphics[width=0.7\linewidth]{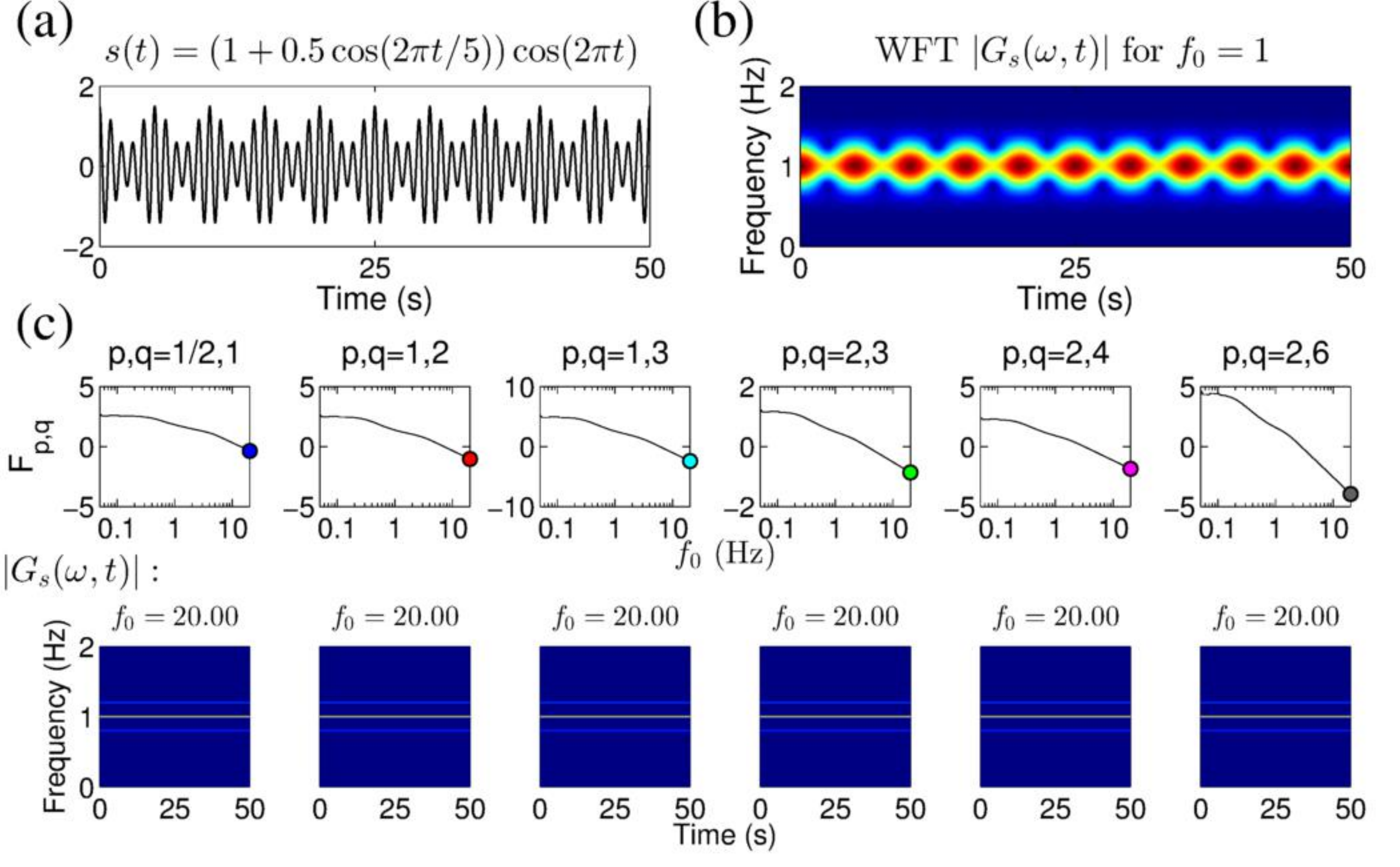}
\includegraphics[width=0.3\linewidth]{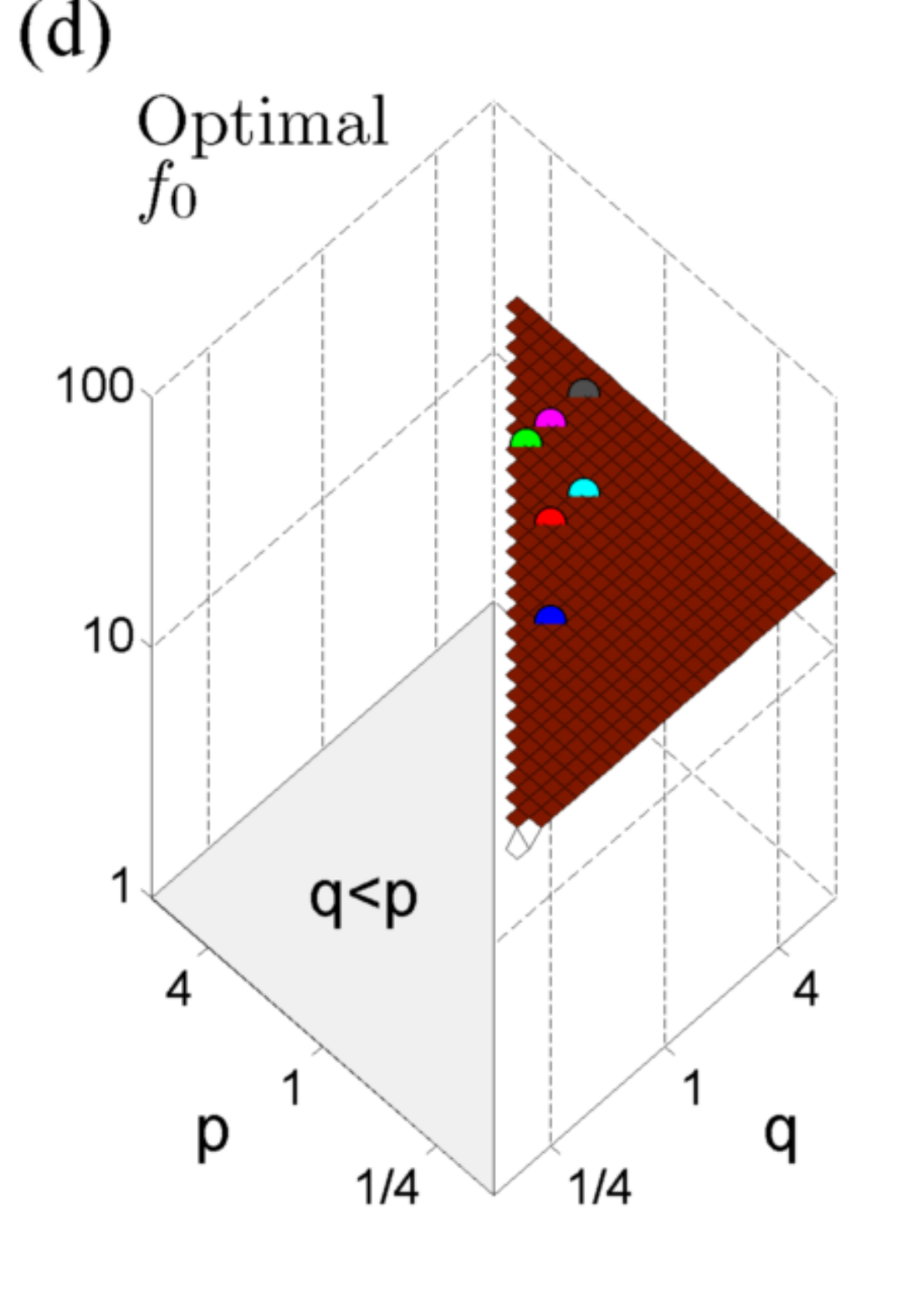}\\
\caption{Same as Fig.\ \ref{TFadapt}, but for the AM component $s(t)=(1+0.5\cos(2\pi t/5))\cos(2\pi t)$, sampled at $50$ Hz for $50$ s. The WFT was always calculated in the frequency range $[0,2]$ Hz, and the optimal $f_0$ was searched in a range $[0.05,20]$.}
\label{AMadapt}
\end{figure*}

Nevertheless, although seemingly promising, the functional (\ref{fpq}) is inappropriate in the general case. It was tested mainly on signals consisting of Gaussian pulses, delta-pulses, sinusoids and chirps (or at least components with frequency modulation containing strong linear term $\sim t$) \cite{Baraniuk:01,Stankovic:01,Jones:94,Jones:90}. But in real application one often deals with AM/FM components occupying some particular band, in which case the functional (\ref{fpq}) often fails to give an appropriate TFR. This is illustrated in Fig.\ \ref{AMadapt} for the example of single AM component with sinusoidal amplitude modulation. As previously discussed (see Sec.\ \ref{sec:AM}), such an AM component can be represented as a sum of tones (\ref{AMs}) which, due to the particular relationships of amplitudes, phases and frequencies, give rise to the amplitude modulation via mutual interference. At the same time, minimizing $F_{p,q}$ minimizes any interference, without differentiating between the desirable and undesirable one; this leads to a separation of AM/FM components into maximally concentrated tones.

\begin{figure*}[t!]
\includegraphics[width=0.7\linewidth]{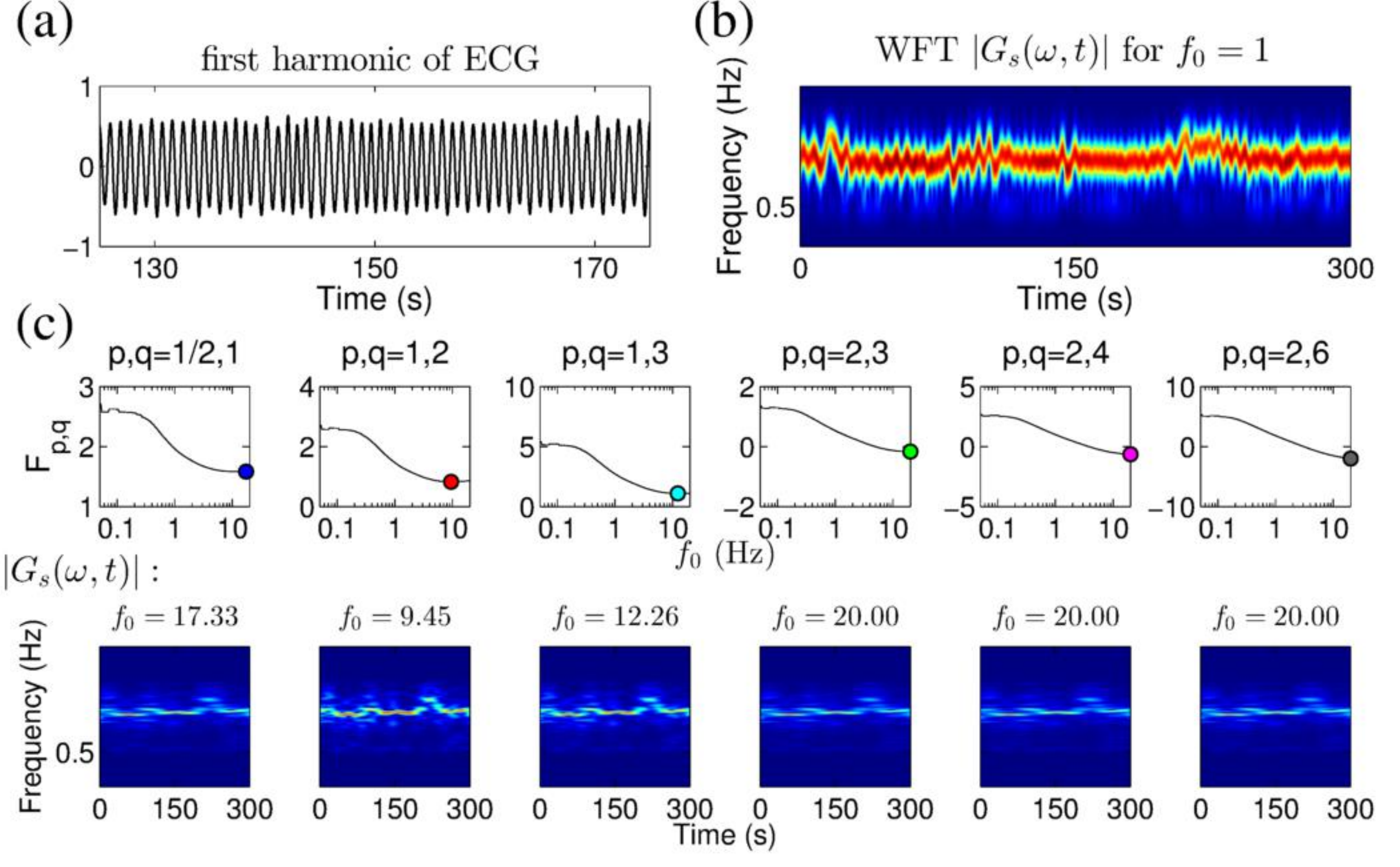}
\includegraphics[width=0.3\linewidth]{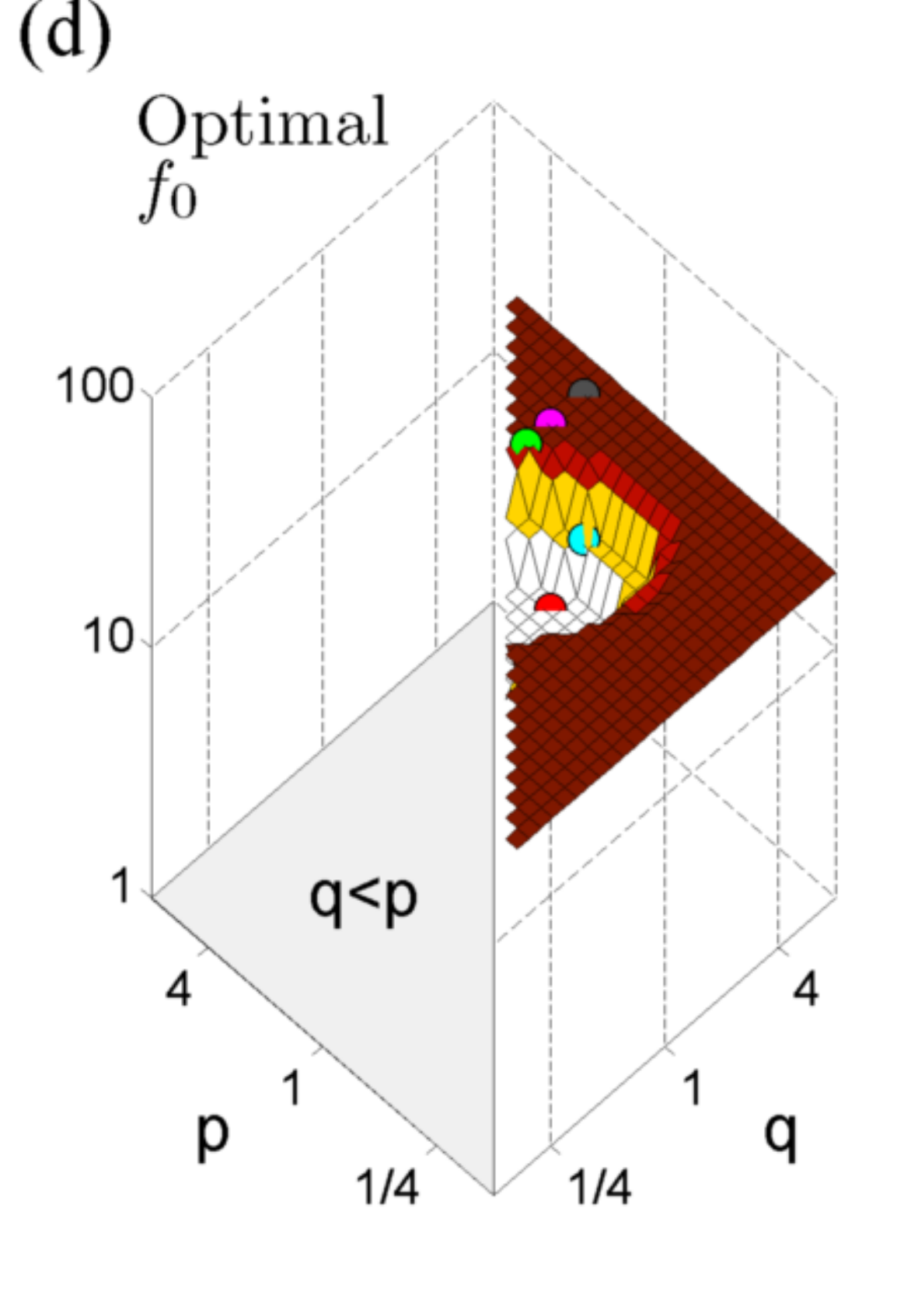}\\
\caption{Same as Fig.\ \ref{TFadapt}, but for the first harmonic of the ECG signal obtained by bandpass filtering it in the region $[0.5,1.5]$ Hz. The signal was sampled at 40 Hz for 30 min, but only central 5 min part is taken for analysis, while other is used for padding to eliminate boundary effects. The WFT was always calculated in the frequency range $[0,2]$ Hz, and the optimal $f_0$ was searched in a range $[0.05,20]$.}
\label{ECGadapt}
\end{figure*}

This issue is not restricted to simple amplitude/frequency modulation, but occurs quite generally for persistent AM/FM components occupying well-defined frequency bands (i.e.\ excluding chirps). This is illustrated in Fig.\ \ref{ECGadapt} for a real ECG signal. As can be seen, although the heart rate modulation is generally quite complex, minimizing $F_{p,q}$ still tries to separate it into tones, giving too large $f_0$, similar to what was seen for a component with simple sinusoidal amplitude modulation. This leads to selection of an inappropriate $f_0$, for which the TFR becomes unsuitable for the analysis, e.g.\ one cannot extract the instantaneous heart frequency from it.

In conclusion, the functionals of \cite{Baraniuk:01,Stankovic:01,Jones:95,Jones:94,Jones:90} are useful for particular types of signals (containing mainly chirps, tones, Gaussian pulses and delta-peaks), but not in general. However, the idea of adapting a TFR by minimizing a suitably chosen functional is very powerful, and possibly some functional that is more universal than \ref{fpq} can be designed.

\subsection{Other approaches}

In addition to the methods already considered, there exist many other ways of selecting optimal window/wavelet parameters. For example, it has been proposed \cite{Hon:12,Zhong:10,Kawahara:99} to construct an adaptive TFR based on certain characteristics (e.g.\ the ridge frequencies) estimated from some initial TFR calculated for particular $f_0$. It is easy to see, however, that these approaches are susceptible to the choice of the initial $f_0$ whereas, if the latter is chosen adaptively based on some variation of the functional (\ref{fpq}), as in \cite{Pei:12}, then the above-mentioned drawbacks of the functional approach will apply.

Another idea is to optimize the WFT/WT based on its local moments in the time-frequency plane \cite{Jones:97}. This method is very expensive computationally, however, being O($N^3$); though the cost might possibly be reduced if estimating the global resolution parameter and not its time-frequency varying version $f_0(\omega,t)$, as originally. Some ``monocomponent'' methods can also be generalized to the case of multicomponent signals by introducing frequency dependence into the resolution parameter \cite{Stankovic:99}, but this will at the same time increase the computational complexity and give rise to additional issues.

Finally, there are various methods \cite{Jones:95,Stankovic:00,Stankovic:01,Baraniuk:01,Baraniuk:93a,Baraniuk:93b,Baraniuk:94,Sang:95,Stankovic:99,Stankovic:98,Katkovnik:98b} for optimizing types of TFR other than the WFT/WT considered here. Most of them represent a modification of one of the approaches already mentioned, but tailored for a particular representation, e.g.\ the Wigner-Ville distribution.

The majority of the existing approaches, however, are only suitable for a particular class of signals, usually those necessarily containing chirps or Gaussian pulses in addition to other components. Thus, to the best of the authors' knowledge, there are at present no universal methods for selecting the optimal window/wavelet parameters. The condition for universality in the present case is that the approach tries to resolve any independent components, e.g.\ tones, but at the same time favors representation of the AM/FM components as single entities. In other words, the interference between independent components should be minimized, while the interference between AM/FM-induced tones should be maximized.

Given the duality of the representation of AM/FM components as a single entity and as a sum of tones, it is questionable whether such a universal approach can in principle be developed. However, what suggests that this may in fact be feasible is that an adaptive signal decomposition method known as basis pursuit \cite{Chen:01} seem to have all the desirable properties, albeit at a computational cost of O($N^3$). Therefore, it might be possible that, for example, the approach (\ref{mfunc}) can be made more universal by choosing instead of (\ref{fpq}) a better functional, which will distinguish between independent components and AM/FM-induced ones by implicitly taking into account the specific relationships between their amplitudes, phases and frequencies. This remains an open problem.

\section{Optimal reconstruction}\label{sec:optrec}

In previous sections we have investigated two possible methods of reconstructing components from their time-frequency supports in the TFR. Ridge reconstruction appears to be more robust to interference and noise, while direct reconstruction performs better in the case of considerable amplitude/frequency modulation. This is quite understandable, since, as we saw, the AM/FM component can be represented as a sum of tones with particular amplitude, phase and frequency relationships (recall (\ref{AMs}) and (\ref{FMs})), so the amplitude and frequency modulation can be viewed as arising due to interference between these AM/FM-induced tones. The more the method is susceptible to interference, therefore, the better it should pick amplitude/frequency modulation, and {\it vice versa}. The remaining questions are how best to choose the method to use for a particular signal, and to decide whether there is any possibility of combining the advantages while simultaneously eliminating the drawbacks of both procedures.

Let us study analytically the errors of the direct and ridge estimates. We will investigate only reconstruction from WFT and WT bearing in mind that, as we saw in previous sections, the SWFT- and SWT-based estimates usually have similar accuracy. Consider a signal consisting of $M+1$ AM/FM components:
\begin{equation}\label{irs}
\begin{aligned}
&s(t)=s_0(t)+\sum_{m=1}^Ms_{m}(t)\\
&s_{0,m}(t)=A_{0,m}(t)\cos \phi_{0,m}(t),\\
&\nu_{0,m}(t)\equiv\phi_{0,m}'(t),\\
\end{aligned}
\end{equation}
for which we want to estimate the parameters of $s_0(t)$, i.e.\ $A_0(t)$, $\phi_0(t)$, $\nu_0(t)$. Evidently, the signal (\ref{irs}) is what one usually deals with in real cases, i.e.\ there is some component of interest surrounded by other ones; often noise is also present, but it can be viewed as large number of tones densely distributed in frequency (see Sec.\ \ref{sec:WN}). The WFT/WT of the summary signal (\ref{irs}) can be represented as
\begin{equation}\label{comptfr}
H_s(\omega,t)=H_{s_0}(\omega,t)+\sum_{m=1}^M H_{s_m}(\omega,t),
\end{equation}
where $H_{s_0}(\omega,t)$ and $H_{s_m}(\omega,t)$ denote the WFT/WT of the component of interest ($s_0(t)$) and side components ($s_m(t)$) in (\ref{irs}), respectively. To avoid dealing with the effects of badly chosen window/wavelet parameters, for the remainder of this section we assume that the behavior of $H_{s_0}(\omega,t)$ is of the IV type, i.e.\ if the signal consisted of only the single $s_0(t)$ then it would be perfectly represented in the TFR.

Suppose that we have successfully extracted the TFS\newline $[\omega_-(t),\omega_+(t)]$ corresponding to the image of $s_0(t)$ in the current TFR, with amplitude peaks at $\omega_p(t)$. Then we can obtain direct/ridge estimates of component's amplitude, phase and frequency $A_0^{(direct|ridge)}$, $\phi_0^{(direct|ridge)}$, $\nu_0^{(direct|ridge)}$ as described in Part I. It is convenient to parametrize the errors of these estimates as
\begin{equation}\label{irTIerr}
\begin{aligned}
\Delta A^{(direct|ridge)}(t)\equiv& A_0^{(direct|ridge)}(t)-A_0(t)\\
\equiv& \Delta A_T^{(direct|ridge)}(t)+\Delta A_I^{(direct|ridge)}(t),\\
\Delta \phi^{(direct|ridge)}(t)\equiv& \phi_0^{(direct|ridge)}(t)-\phi_0(t)\\
\equiv& \Delta \phi_T^{(direct|ridge)}(t)+\Delta \phi_I^{(direct|ridge)}(t),\\
\Delta \nu^{(direct|ridge)}(t)\equiv& \nu_0^{(direct|ridge)}(t)-\nu_0(t)\\
\equiv& \Delta \nu_T^{(direct|ridge)}(t)+\Delta \nu_I^{(direct|ridge)}(t).
\end{aligned}
\end{equation}
In (\ref{irTIerr}), $\Delta A^{(direct|ridge)}_T(t)$ is the theoretical, inherent inaccuracy of amplitude reconstruction, i.e.\ the error it would have if there were only one component in the signal ($s(t)=s_0(t)$) and its extracted TFS $[\omega_-(t),\omega_+(t)]$ contained all of its power. The second term, $\Delta A^{(direct|ridge)}_I(t)$, represents the error related to interference with the other components $s_m(t)$ present in the signal. Due to this interference the extracted TFS $[\omega_-(t),\omega_+(t)]$ will usually contain only some proportion of the $s_0(t)$, mixed with parts of other $s_m(t)$. The same classification applies to the phase and frequency reconstruction errors. Below we assume that both the theoretical and interference-related errors are small enough that the first order expansion over them is valid.

\subsection{Ridge reconstruction errors}\label{sec:optrecA}

Assuming for simplicity a continuous frequency scale (so there are no discretization errors), the WFT/WT-based ridge estimates can be considered in the form (see Part I):
\begin{equation}\label{estridge}
\nu^{(ridge)}(t)=\omega_p(t),\;
A^{(ridge)}(t)e^{i\phi^{(ridge)}(t)}=\frac{2H_s(\omega_p(t),t)}{\hat{h}_{\max}}.
\end{equation}
To the best of our knowledge, the latest theoretical estimates of ridge errors are given in \cite{Lilly:10}:
\begin{equation}\label{irTridge}
\begin{aligned}
&\Delta A_T^{(ridge)}(t)=\frac{1}{2}P^2\big(\nu_0(t)\big)A_0''(t)+O(\delta_{N_T}^2),\\
&\Delta\phi_T^{(ridge)}(t)=\frac{1}{2}P^2\big(\nu_0(t)\big)\nu_0'(t)+O(\delta_{N_T}^2),\\
&\Delta\nu_T^{(ridge)}(t)=P^2\big(\nu_0(t)\big)\left(\frac{1}{2}\nu_0''(t)
+\frac{A_0'(t)}{A_0(t)}\nu_0'(t)\right)+O(\delta_{N_T}^3),\\
&P^2(\omega)\equiv \left[\frac{\partial_\nu^2\hat{h}_\nu(\omega)}{\hat{h}_{\nu}(\omega)}\right]_{\nu=\omega}=
\left[\begin{array}{l}
-\frac{\hat{g}''(0)}{\hat{g}(0)}\mbox{ for the WFT},\\
-\frac{\omega_\psi^2}{\omega^2}\frac{\hat{\psi}''(\omega_\psi)}{\hat{\psi}(\omega_\psi)}\mbox{ for the WT},\\
\end{array}\right.
\end{aligned}
\end{equation}
where the value of $\delta_{N_T}$ is determined by the strength of the amplitude/frequency modulation in relation to the window/wavelet parameters, see \cite{Lilly:10}. The latter is relatively small when the behavior of the WFT/WT $H_{s_0}(\omega,t)$ is of the IV type, as assumed in this section, but might become non-negligible otherwise, requiring one to take account of higher order terms in (\ref{irTridge}). The magnitude of the error is additionally determined by the proportionality factor $P^2(\omega)$ in (\ref{irTridge}), which is $P^2(\omega)=f_0^2$ ($P^2(\omega)=(2\pi f_0)^2/\omega^2$) for Gaussian window (lognormal wavelet).

By numerical simulation, we have found that the theoretical prediction (\ref{irTridge}) matches extremely well with the actual ridge reconstruction errors for both the WFT and WT, at least in the case of Regime IV behavior. Importantly, the inaccuracy of ridge frequency estimation $\Delta\nu_T^{(ridge)}(t)$ (\ref{irTridge}) is of the second order in $\delta_{N_T}$ and so is its contribution to $\Delta A_T^{(ridge)}(t),\Delta \phi_T^{(ridge)}(t)$ \cite{Lilly:10}. Note that, for windows/wavelets with $\hat{g}''(0)=0$ and $\hat{\psi}''(\omega_\psi)=0$, one has $P(\omega)=0$, and all theoretical ridge errors (\ref{irTridge}) become of higher order in $\delta_{N_T}$; however, the reconstruction accuracy is to a large extent determined by the TFR behavior, so the time-frequency resolution of the window/wavelet is the most important thing.

\begin{remark}
Only the WT was considered in \cite{Lilly:10}, so that the expression for the WFT (\ref{irTridge}) was not given. Its rigorous derivation requires cumbersome calculations based on the machinery developed in \cite{Lilly:10}, but from simple logical considerations based on comparison of the linear/logarithmic frequency resolution of WFT/WT one can infer the expressions given in (\ref{irTridge}); their correctness was also confirmed numerically.
\end{remark}

\begin{remark}
Interestingly, for the case of AM component (which can always be represented in the form $s(t)=A(t)\cos(\nu t+\varphi)=(a_0+\sum_ka_k\cos(\nu_k t+\varphi_k))\cos(\nu t+\varphi)$, where we assume all $\nu_k<\nu$), one can derive an exact theoretical error for WFT-based ridge reconstruction if $\hat{g}(\xi)$ is symmetric. As discussed previously (see Sec.\ \ref{sec:AM}), there will be no errors of phase and frequency estimation: $\Delta\phi_T^{(ridge)}=\Delta\nu_T^{(ridge)}=0$. Next, the peak in the WFT amplitude will always occur at $\omega_p(t)=\nu$, with the WFT there being $H_s(\nu,t)=(A/2)e^{i(\nu t+\varphi)}[a_0\hat{g}(0)+\sum_k a_k \hat{g}(\nu_k)\cos(\nu_k t+\varphi_k)]$, which gives the ridge amplitude estimate\newline $A^{(ridge)}(t)=\frac{1}{2\pi}\int \hat{g}(\xi)\hat{A}(\xi)e^{i\xi t}d\xi$.
\end{remark}

Determination of the interference-related inaccuracy is a highly non-trivial task. To approach it, let us make the very rough assumption that the components of the TFRs (\ref{comptfr}) at the ridge points $\omega_p(t)$ are
\begin{equation}\label{ridgeapp}
\begin{aligned}
H_{s_0}(\omega_p(t),t)\approx& H_{s_0}(\tilde{\omega}_p(t),t)\\
=& \frac{A_0(t)+\Delta A_T^{ridge}(t)}{2}\hat{h}_{\nu_0(t)}(\nu_0(t))e^{i(\phi_0(t)+\Delta\phi_T^{ridge}(t))},\\
H_{s_m}(\omega_p(t),t)\approx& \frac{A_m(t)}{2}\hat{h}_{\nu_m(t)}(\nu_0(t)),
\end{aligned}
\end{equation}
where $\tilde{\omega}_p(t)=\nu_0(t)+\Delta\nu^{(ridge)}_T$ denote the positions of the ridge points in $H_{s_0}(\omega,t)$.

The first approximation in (\ref{ridgeapp}) can be shown to be accurate to the second order over the reconstruction errors. Thus, denoting the TFR phase as $\phi_H(\omega,t)\equiv {\rm arg}[H_{s_0}(\omega,t)]$, and taking into account that $\Delta\nu^{(ridge)}_I(t)=\omega_p(t)-\tilde{\omega}_p(t)$, one has
\begin{equation}\label{rapp}
\begin{aligned}
H_{s_0}(\omega_p(t),t)
=&H_{s_0}(\tilde{\omega}_p(t),t)+\Big[|\partial_\omega H_{s_0}(\tilde{\omega}_p(t),t)|\\
&+i|H_{s_0}(\tilde{\omega}_p(t),t)|\partial_\omega\phi_H(\tilde{\omega}_p(t),t)\Big]
e^{i\phi_H(\tilde{\omega}_p(t),t)}\Delta\nu^{(ridge)}_I(t)\\
&+O\big([\Delta\nu^{(ridge)}_I(t)]^2\big).
\end{aligned}
\end{equation}
By definition $|\partial_\omega H_s(\tilde{\omega}_p(t),t)|=0$, while for the assumed Regime IV behavior of $H_{s_0}(\omega,t)$ the frequency-derivative of the TFR phase at the peak $\partial_\omega\phi_H(\tilde{\omega}_p(t),t)$ is of the first order over the theoretical ridge errors (\ref{irTridge}) \cite{Lilly:10} (while for tones it is exactly zero at all frequencies).

The quality of the second approximation in (\ref{ridgeapp}) is in general harder to estimate. However, e.g.\ for tones it can be easily seen that $2H_{s_m}(\omega_p(t),t)=A_m\hat{h}_{\nu_m(t)}(\omega_p(t))=A_m\hat{h}_{\nu_m(t)}(\nu_0(t))+O\big(\hat{h}_{\nu_m(t)}(\nu_0(t))\big) \Delta\nu^{(ridge)}(t)$. Therefore, since $H_{s_m}(\omega_p(t),t)$ is by itself proportional to the interference-related error (so that $\hat{h}_{\nu_m(t)}(\nu_0(t))$ can be assumed small), the expression for $H_{s_m}(\omega_p(t),t)$ (\ref{ridgeapp}) is valid up to the second order over the reconstruction errors. So it follows that, when $s_m(t)$ have slow amplitude and frequency variations, the second approximation in (\ref{ridgeapp}) holds. Note, that for the following derivations it does not need to be valid for all $s_m(t)$, but only for those making some contribution to the interference errors, i.e.\ having non-negligible $|H_{s_m}(\omega_p(t),t)/H_{s_0}(\omega_p(t),t)|$.

Based on (\ref{ridgeapp}) and the ridge reconstruction formulas, one can show the interference-related errors of the ridge method to be
\begin{equation}\label{irIridge}
\begin{aligned}
\Delta A_I^{(ridge)}(t)\approx& \left|A_0(t)+\sum_{m=1}^M A_{m}(t)\frac{\hat{h}_{\nu_{-1}(t)}(\nu_0(t))}{\hat{h}_{\max}}e^{i(\phi_{m}(t)-\phi_0(t))}\right|-A_0(t)\\
\approx&\sum_{m=1}^M A_m(t)\frac{\hat{h}_{\nu_m(t)}(\nu_0(t))}{\hat{h}_{\max}}\cos(\phi_m(t)-\phi_0(t)),\\
\Delta \phi_I^{(ridge)}(t)\approx& {\rm arg}\left[A_0(t)+\sum_{m=1}^MA_{m}(t)\frac{\hat{h}_{\nu_{m}(t)}(\nu_0(t))}{\hat{h}_{\max}}e^{i(\phi_{m}(t)-\phi_0(t))}\right]\\
\approx&\sum_{m=1}^M \frac{A_m(t)\hat{h}_{\nu_m(t)}(\nu_0(t))}{A_0(t)\hat{h}_{\max}}\sin(\phi_m(t)-\phi_0(t)),\\
\Delta \nu_I^{(ridge)}(t)=&\nu_H(\omega_p(t),t)-\nu_0(t)-\Delta\nu_T^{ridge}(t)\\
\approx& \sum_{m=1}^M\frac{A_m(t)\hat{h}_{\nu_m(t)}(\nu_0(t))}{A_0(t)\hat{h}_{\max}}
[\nu_m(t)-\nu_0(t)]\cos(\phi_m(t)-\phi_0(t)),\\
\end{aligned}
\end{equation}
where the expression for $\Delta \nu_I^{(ridge)}(t)$ was derived using (\ref{gformGEN}) with $\omega=\nu_0(t),\;\{a_n,\nu_n\}\rightarrow \{A_{n}(t),\nu_n(t)\}$ (the motivation behind this being the same as for (\ref{ridgeapp})), and the approximation $\omega_p(t)\approx \nu_H(\omega_p(t),t)\approx \nu_H(\nu_0(t),t)$, which is of second order over the ridge reconstruction errors \cite{Lilly:10}.

\subsection{Direct reconstruction errors}\label{sec:optrecB}

When there is only one AM/FM component and the TFR behavior is of the IV type, the direct estimates are by definition exact (up to the accuracy $\epsilon$ with which Regime IV is determined), so that there are no theoretical errors:
\begin{equation}\label{irTdir}
\Delta A_T^{(direct)}(t)=\Delta \phi_T^{(direct)}(t)=\Delta \nu_T^{(direct)}=0.
\end{equation}

The interference-related errors are more sophisticated, and to treat them we employ a rough simplification similar to (\ref{ridgeapp}) used for the ridge case. Thus, we assume that
\begin{equation}\label{dirapp}
\begin{aligned}
\int_{\mu(\omega_-(t))}^{\mu(\omega_+(t))} H_{s_0}(\omega,t)d\mu(\omega)
&\approx \int_{\mu(\omega_-(t))}^{\mu(\omega_+(t))} \frac{A_0(t)}{2}\hat{h}_{\nu_0(t)}(\omega)e^{i\phi_0(t)}d\mu(\omega),\\
\int_{\mu(\omega_-(t))}^{\mu(\omega_+(t))} H_{s_m}(\omega,t)d\mu(\omega)
&\approx \int_{\mu(\omega_-(t))}^{\mu(\omega_+(t))} \frac{A_m(t)}{2}\hat{h}_{\nu_m(t)}(\omega)e^{i\phi_m(t)}d\mu(\omega),\\
\end{aligned}
\end{equation}
This approximation is evidently consistent in the sense that, for $\mu(\omega_{\pm}(t))=\pm\infty$, it becomes exact, with the first line of (\ref{dirapp}) being equal to $s_0^a(t)$ (the analytic signal of the component considered), and the second one equal to $s_m^a(t)$. Moreover, (\ref{dirapp}) is also exact when all components are represented by tones. Therefore, the approximation (\ref{dirapp}) holds in the case when component have slowly-varying amplitudes and frequencies (in respect to window/wavelet time resolution). Note that, for the estimation of interference-related errors, the approximation (\ref{dirapp}) does not need to be valid for $s_m(t)$ with negligible $\big|\int_{\mu(\omega_-(t))}^{\mu(\omega_+(t))} H_{s_m}(\omega,t)d\mu(\omega)\big/\int_{\mu(\omega_-(t))}^{\mu(\omega_+(t))} H_{s_0}(\omega,t)d\mu(\omega)\big|$, which thus almost do not interfere with the component of interest.

Substituting (\ref{comptfr}) and (\ref{dirapp}) into the direct estimation formulas (see Part I), one obtains the interference-related errors as
\begin{equation}\label{irIdir}
\begin{aligned}
X(\omega,t)\equiv& A_0(t)\hat{h}_{\nu_0(t)}(\omega)+
\sum_{m=1}^MA_m(t)\hat{h}_{\nu_m(t)}(\omega)e^{i(\phi_m(t)-\phi_0(t))}\\
\Delta A_I^{(direct)}(t)\approx&\bigg|\frac{C_h^{-1}}{2}
\int_{\mu(\omega_{-}(t))}^{\mu(\omega_{+}(t))}X(\omega,t)d\mu(\omega)\bigg|-A_0(t)\\
\approx&-A_0(t)\Big[1-\widetilde{Q}_{\nu_0(t)}\big(\omega_-(t),\omega_{+}(t)\big)\Big]\\
&+\sum_{m=1}^MA_m(t)\widetilde{Q}_{\nu_{m}(t)}\big(\omega_-(t),\omega_{+}(t)\big)\cos(\phi_m(t)-\phi_0(t)),\\
\Delta \phi_I^{(direct)}(t)=&{\rm arg}\left[\frac{C_h^{-1}}{2}\int_{\mu(\omega_{-}(t))}^{\mu(\omega_{+}(t))}X(\omega,t)d\mu(\omega)\right]\\
=&\sum_{m=1}^M\frac{A_m(t)}{A_0(t)}\widetilde{Q}_{\nu_m(t)}(\omega_-(t),\omega_{+}(t))\sin(\phi_m(t)-\phi_0(t)),\\
\Delta \nu_I^{(direct)}(t)=&-\nu_0(t)-\overline{\omega}_h+{\rm Re}
\frac{\big(D_h^{-1}/2\big)\int_{\mu(\omega_{-}(t))}^{\mu(\omega_{+}(t))}X(\omega,t)
\omega d\mu(\omega)}{[A_0(t)+\Delta A^{(direct)}(t)]e^{i\Delta\phi^{(direct)}(t)}}\\
\approx& -(\nu_0(t)+\overline{\omega}_h)\frac{\Delta A_I^{(direct)}}{A_0(t)}
-\frac{D_h^{-1}}{2}\\
&\times\bigg[\int_{-\infty}^{\mu(\omega_-(t))}\hat{h}_{\nu_0(t)}(\omega)\omega d\mu(\omega)\\
&+\int_{\mu(\omega_+(t))}^{\infty}\hat{h}_{\nu_0(t)}(\omega)\omega d\mu(\omega)
+\sum_{m=1}^M \frac{A_m(t)}{A_0(t)}\\
&\times\cos(\phi_m(t)-\phi_0(t))\int_{\mu(\omega_-(t))}^{\mu(\omega_+(t))}\hat{h}_{\nu_0(t)}(\omega)\omega d\mu(\omega)\bigg],\\
\end{aligned}
\end{equation}
where we have denoted
\begin{equation}\label{nIdir}
\begin{aligned}
\overline{\omega}_h\equiv&
\left[\begin{array}{l}
\overline{\omega}_g\equiv \big(C_g^{-1}/2\big)\int \omega \hat{g}(\omega)d\omega\;\mbox{ for the WFT},\\
0\;\mbox{ for the WT},\\
\end{array}\right.
\\
D_h\equiv&
\left[\begin{array}{l}
C_g\;\mbox{ for the WFT},\\
D_\psi\equiv \frac{\omega_\psi}{2}\int_0^\infty \hat{\psi}^*(\omega)\frac{d\omega}{\omega^2}\;\mbox{ for the WT},\\
\end{array}\right.
\end{aligned}
\end{equation}
and the expressions for $\Delta \nu_I^{(direct)}(t)$ were derived assuming the slightly modified form of (\ref{dirapp}) with $d\mu(\omega)\rightarrow\omega d\mu(\omega)$, for which the same considerations apply.

When $D_\psi=\infty$, as for the Morlet wavelet, one is forced to reconstruct the frequency by the hybrid method (see Part I). Under the assumption that $\int_{\omega_-(t)}^{\omega_+(t)}\nu_H(\omega,t)H_s(\omega,t)d\omega\approx \int_{\omega_-(t)}^{\omega_+(t)}\tilde{\nu}_H(\omega,t)H_s(\omega,t)d\omega$, where $\tilde{\nu}_H(\omega,t)$ is given by $\nu_H(\omega,t)$ in (\ref{gformGEN}) with $\{a_n,\nu_n\}\rightarrow\{A_n(t),\nu_n(t)\}$, the hybrid frequency estimation errors are
\begin{equation}\label{irTIhybrid}
\begin{aligned}
\Delta\nu_T^{(hybrid)}(t)\approx&0,\\
\Delta\nu_I^{(hybrid)}(t)\approx&\sum_{m=1}^M\frac{A_m(t)}{A_0(t)}[\nu_m(t)-\nu_0(t)]
Q_{\nu_m(t)}\big(\omega_-(t),\omega_+(t)\big)\\
&\times\cos(\phi_m(t)-\phi_0(t)),
\end{aligned}
\end{equation}
which can be derived in a similar way to that used for the expressions in (\ref{irIdir}).

\subsection{Advantages and drawbacks of each method}\label{sec:optrecC}

Comparing (\ref{irTridge}) and (\ref{irTdir}), it is clear that direct method outperforms ridge method in terms of theoretical error. At the same time, as seen from (\ref{irIridge}) and (\ref{irIdir}), in the direct method one picks the contribution of the side components over all TFS, while ridge reconstruction accounts for interference only at the peak. As a result, the ridge estimates are superior to the direct ones in respect of the interference-related errors.

Therefore, the choice of the method depends very much on the particular signal and the representation of the component of interest in its TFR. When the noise is small and different components are well-separated in the TFR, then the direct method should be used; otherwise, if the noise and/or the interference with other components is strong, the ridge method is the better choice. Ridge reconstruction is also superior if the component of interest has no or very weak (in terms of the window/wavelet time-resolution) amplitude and frequency variations, implying a small theoretical error (\ref{irTridge}). Thus, the latter is exactly zero for tones, in which case ridge estimates are always the best. Recall also, that ridge reconstruction is less susceptible to boundary effects than the direct estimation (see Sec.\ \ref{sec:OV}).

Regarding the TFR and the form of the window/wavelet, the (S)WFT with a window function $\hat{g}(\xi)$ that is symmetric in frequency provides clear advantages in terms of the accuracy of the resultant estimates. For example, as mentioned above, there are no phase/frequency reconstruction errors for the AM components in this case, while otherwise they exist. This is because the AM/FM-induced tones appear symmetrically around the main tone (see Sec.\ \ref{sec:AM} and \ref{sec:FM}), so that a symmetric $\hat{g}(\xi)$ is the best form to match this structure; the (S)WFT with frequency-asymmetric windows, and the (S)WT (due to its logarithmic frequency scale), do not reflect such a symmetry. Note also, that $\hat{\psi}(\xi)$ that is symmetric on a logarithmic scale, such as the lognormal wavelet, offers slightly better reconstruction possibilities as compared to other wavelets; this is because at low $\Delta\nu/\nu$ one has $\log(1+\Delta\nu/\nu)\approx-\log(1-\Delta\nu/\nu)\approx\Delta\nu/\nu$, so that the symmetries of $\nu\pm\Delta\nu$ around $\nu$ on linear and logarithmic scales become nearly equivalent.

The most important characteristic of the window/wavelet, however, is its time-frequency resolution. Thus, for a single component one can usually adjust parameters to represent and reconstruct it perfectly using any window/wavelet. But in real cases, when the signal consists of a number of components, there is usually no choice of parameters for which all components can be recovered perfectly, and one needs to make a compromise. The time-frequency resolution determines how good such compromise might be in principle, i.e.\ the best accuracy with which all components can be reconstructed. The choice between the (S)WFT and (S)WT, on the other hand, depends on the signal properties, as discussed in Part I: (S)WT is to be preferred when the AM/FM components at lower frequencies are closer to each other and less time-varying than those at higher frequencies, while the (S)WFT is more suitable otherwise.

\subsection{Adaptive choice of the method}\label{sec:optrecD}

To choose the best method automatically, one can devise an empirical criterion as follows. Suppose we have calculated the TFR of a signal and extracted from it the ridge curve $\omega_p(t)$ and TFS $[\omega_-(t),\omega_+(t)]$ corresponding to some component. Its associated parameters can then be reconstructed by both the direct and ridge methods; the resultant estimates will be denoted as $A^{(d,r)}(t)$, $\phi^{(d,r)}(t)$ and $\nu^{({d,r})}(t)$, where ``d'' and ``r'' stand for ``direct'' and ``ridge'', respectively. To understand which reconstruction method is more accurate, we calculate the TFR (using the same window/wavelet as originally) of the signal $s^{(d)}(t)=A^{(d)}(t)\cos\phi^{(d)}(t)$, extract the ridge curve and TFS from it (taking simple maxima $\omega_p(t)=\operatorname{argmax}_\omega|H_s(\omega,t)|$ is sufficient here), and reconstruct by the direct method the ``refined'' parameters $\tilde{A}^{(d)}(t),\tilde{\phi}^{(d)}(t),\tilde{\nu}^{(d)}(t)$. The same procedure is performed for the ``ridge'' signal $s^{(r)}(t)=A^{(r)}(t)\cos\phi^{(r)}(t)$, now using the ridge method to reconstruct the refined estimates.

Obviously, if e.g.\ the direct estimates are accurate, one should have $\{\tilde{A}^{(d)}(t),\tilde{\phi}^{(d)}(t),\tilde{\nu}^{(d)}(t)\}\approx\{A^{(d)}(t),\phi^{(d)}(t),\nu^{(d)}(t)\}$. Therefore, one can assess which method is better on the basis of the \emph{discrepancies} between the original and refined estimates, which can be quantified using the corresponding relative errors (\ref{apfrel}) as
\begin{equation}\label{apfrel0}
\begin{aligned}
&\tilde{\varepsilon}_a^{(d,r)}\equiv \kappa_a^{(d,r)}\frac{\sqrt{\langle (\tilde{A}^{(d,r)}(t)-A^{(d,r)}(t))^2 \rangle}}{\langle[A^{(d,r)}(t)]^2\rangle},\\
&\tilde{\varepsilon}_\phi^{(d,r)}\equiv \kappa_\phi^{(d,r)}\sqrt{1-|\langle e^{i(\tilde{\phi}^{(d,r)}(t)-\tilde{\phi}^{(d,r)}(t))} \rangle|^2},\\
&\tilde{\varepsilon}_\nu^{(d,r)}\equiv \kappa_\nu^{(d,r)}\frac{\sqrt{\langle (\tilde{\nu}^{(d,r)}(t)-\nu^{(d,r)}(t))^2\rangle}}{2\pi},\\
\end{aligned}
\end{equation}
where $\kappa_{a,\phi,\nu}^{(d,r)}$ are the coefficients that can be used to tune the performance of the approach (they were found empirically to be $\kappa_{a,\phi,\nu}^{(d)}=\{3,4,2\}$, $\kappa_{a,\phi,\nu}^{(r)}=1$). For each parameter, the choice between its direct and ridge estimate is then made based on the corresponding discrepancy (\ref{apfrel0}): the smaller it is, the more accurate the reconstructed parameter is expected to be.

\begin{figure*}[t!]
\includegraphics[width=1.0\linewidth]{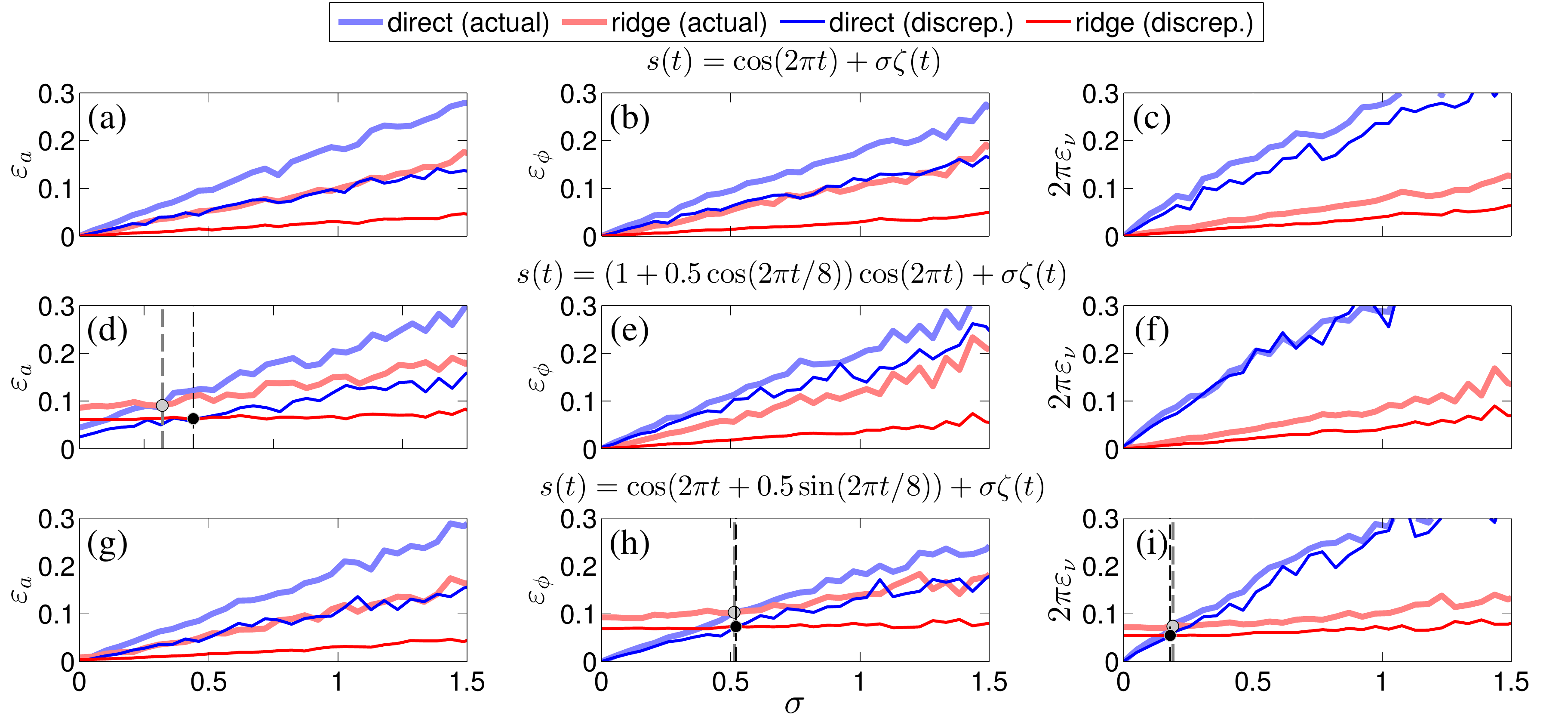}\\
\caption{The actual reconstruction errors (\ref{apfrel}) of the direct and ridge methods (light-blue and light-red lines, respectively) and the corresponding discrepancies (\ref{apfrel0}) (direct -- blue, ridge -- red) in their dependence on the noise level $\sigma$. (a,d,g): Amplitude reconstruction errors. (b,e,h): Phase reconstruction errors. (c,f,i): Frequency reconstruction errors. The signals associated with each row are given by the equations above the central panels (b,e,h), with $\zeta(t)$ denoting unit-deviation Gaussian white noise; each signal was sampled at 50 Hz for 200 s. Where present, the gray (or black) points with the corresponding dashed lines indicate the intersections between the true errors (\ref{apfrel}) (or the discrepancies (\ref{apfrel0})) of the direct and ridge methods.}
\label{fig:choicerec}
\end{figure*}

Despite being empirical, the approach outlined above works very well in practice, selecting the best estimates in the majority of cases. This is illustrated in Fig.\ \ref{fig:choicerec}, where the discrepancies (\ref{apfrel0}) are shown together with the actual reconstruction errors (\ref{apfrel}) for each method. As can be seen, the values of $\tilde{\varepsilon}_{a,\phi,\nu}^{(d,r)}$ are proportional to the true errors and allow one to judge reliably about the relative performance of the two reconstruction methods. Thus, as discussed previously, for a single tone signal embedded in noise the ridge estimates are always preferred, and the criterion based on (\ref{apfrel0}) correctly reflects this fact (see Fig. \ref{fig:choicerec}(a-c)). Next, when amplitude/frequency modulation is present, at low noise levels the direct estimates are preferred, but with increasing noise strength their inaccuracy grows faster than in the case of ridge reconstruction. Therefore, beyond some threshold noise level (indicated by gray vertical dashed lines in Fig.\ \ref{fig:choicerec}) ridge estimates become the more accurate; this threshold and the optimal method in each case can be well recovered from the behavior of the discrepancies (\ref{apfrel0}), as is clear from Fig.\ \ref{fig:choicerec}(d,h,i).

\section{Concentration or resolution? Do we really need synchrosqueezing?}\label{sec:Synchroneed}

As we have seen, the behavior of the WFT/WT projects onto the SWFT/SWT. Thus, if e.g.\ two tones are not well separated in the WFT/WT, then they will be not well separated in the SWFT/SWT as well. Furthermore, in all examples studied so far synchrosqueezing did not significantly improve (but often worsened) the accuracy of parameters' reconstruction by both the direct and ridge methods. Taken together, this indicates that synchrosqueezing does not increase time or frequency resolution, as might have seemed the case at the first glance: it only improves the ``readability'' of the TFR \cite{Auger:95}, providing a more visually appealing picture.

While not providing considerable advantages, the SWFT/SWT has a few drawbacks in comparison to the WFT/WT, namely:
\begin{enumerate}
\item The SWFT/SWT amplitude depends on the discretization of the frequency scale, making ridge amplitude reconstruction ill-defined. At the same time, for the usual WFT/WT one can estimate the amplitude by both the direct and ridge methods.
\item The behavior of the SWFT/SWT is more complicated than that of the WFT/WT, being harder to study both analytically and practically. Thus, even in the case when all components are well represented and there is no noise, synchrosqueezed TFRs might still have many side TFSs containing a small amount of power. Additionally, components with fast amplitude or frequency modulation might be represented in the SWFT/SWT in a quite weird way (see e.g.\ Fig.\ \ref{FMex} and the related discussion).
\item For the WFT/WT, both ridge and direct estimates of the component's frequency are relatively unaffected by frequency discretization effects, while accurate estimation of instantaneous frequency from the SWFT/SWT by any method requires very small frequency bins, thus increasing its computational cost (see Part I).
\end{enumerate}
Hence, the usefulness of synchrosqueezing is questionable, because in terms of components reconstruction it only introduces additional complications, while not providing significant advantages. Even in terms of ridge curve extraction, i.e.\ tracing the components in the time-frequency plane, the SWFT/SWT also does not seem to be more suitable than the WFT/WT. Note, however, that in mathematical terms synchrosqueezing does not bring much disadvantage either, i.e.\ the results obtained from the SWFT/SWT in an appropriate way are correct and will be qualitatively, and to a large extent quantitatively, the same as the corresponding results obtained from the WFT/WT.

\begin{remark}
As discussed in Sec.\ \ref{sec:assumptions}, in this (second) part of the work we consider $\hat{g}(\xi)$ and $\hat{g}(\xi>0)$ to be at least approximately unimodal. For multimodal windows/wavelets, on the other hand, synchrosqueezing has the advantageous property of joining together the component's power contained in all sidelobes into the one TFS in the SWFT/SWT, hence making the latter more interpretable than the underlying WFT/WT. Thus, in the process of synchrosqueezing one utilizes the relationships between the instantaneous frequencies $\nu_{H}(\omega,t)$, implicitly determining and differentiating between the sidelobes corresponding to independent components and those corresponding to the same one. However, in practice it appears that such a property is greatly affected even by small interference between components (or by considerable amplitude/frequency modulation), in which case the behavior of the SWFT/SWT becomes very complex, with the power of the component often being distributed over few TFSs. Generally, the advantages of synchrosqueezing for windows/wavelets which are multimodal in frequency is a separate topic. In any case, multimodal $\hat{g}(\xi)$ and $\hat{\psi}(\xi>0)$ are rarely used because of being inconvenient in terms of the resultant representation, as well as usually having poor time-frequency resolution.
\end{remark}

The fact that synchrosqueezing increases the TFR concentration, but at the same time does not give better results in terms of resolving components in frequency or representing the time-variability of their parameters, leads to reconsideration of a more general question: does the concentration of the TFR alone represent the main measure of its performance, as is often believed? Our results argue against such a view.

In general, the ``ideal'' representation of a signal $s(t)=$\\
$\sum_k A_k(t)\cos\phi_k(t)$ can be regarded as being $I(\omega,t)\sim \sum_k A_k(t)\delta(\omega-\phi_k'(t))$. Hence, the inverse of the (somehow defined) ``distance'' between the perfect representation $I(\omega,t)$ and the calculated TFR can be considered as a measure of its performance. What one aims to achieve, therefore, is not just to increase the TFR concentration, but to increase it \emph{around the instantaneous frequencies $\phi_k'(t)$} and/or to improve the representation of the amplitude variations. For example, a TFR having peaks at $\omega=\phi_k'(t)$, but not being too concentrated, is obviously to be preferred to an extremely concentrated TFR with peaks distant from the true instantaneous frequencies. In other words, the main goal is to represent appropriately all the components present in the signal, so that their parameters can accurately be recovered. The most important characteristics of the TFR are therefore its resolution properties and their conformity with the signal, and not simply the concentration.

\begin{figure*}[t!]
\includegraphics[width=1.0\linewidth]{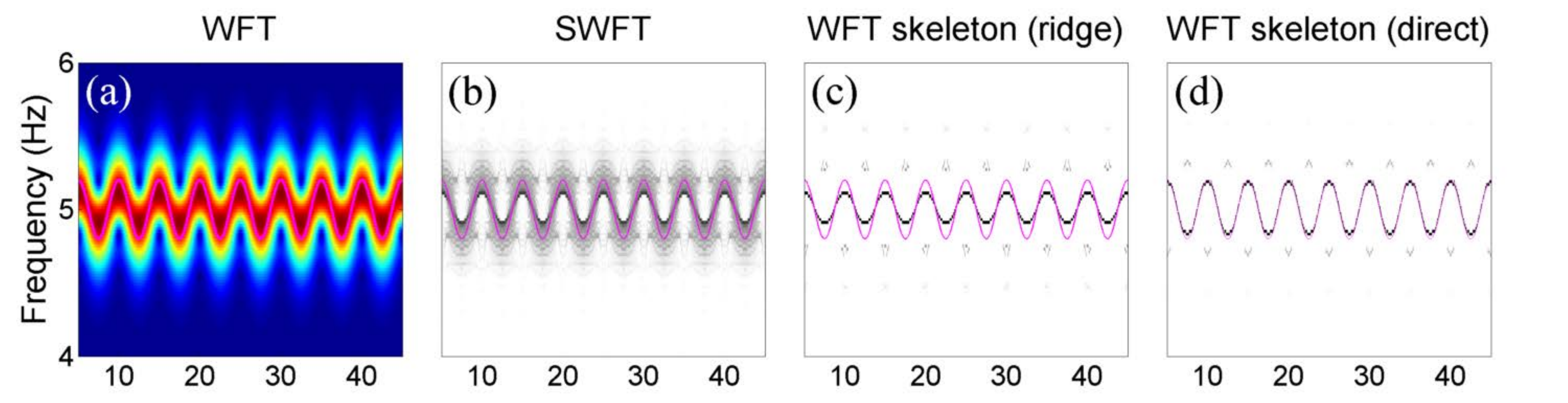}
\includegraphics[width=1.0\linewidth]{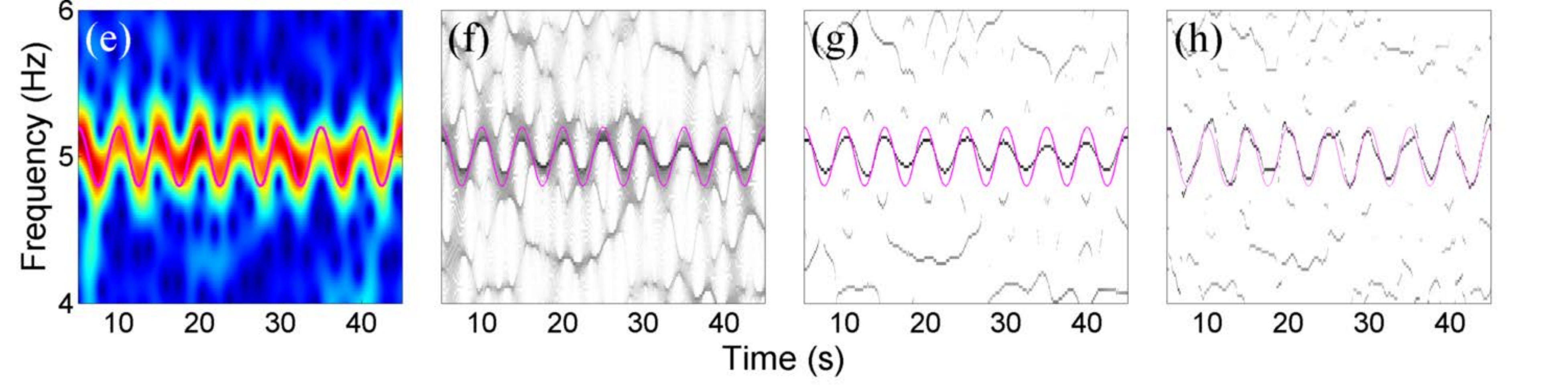}
\caption{Comparison of the WFT, SWFT and WFT skeletons based on ridge and direct reconstruction for: (a-d) the FM component $s(t)=\cos(10\pi t+\sin(2\pi t/5))$; (e-h) the same component additionally corrupted by white noise of $\sqrt{2}$ standard deviation. Magenta lines show the true frequency of the component. The small intermittent components appearing at both sides of the main frequency in (c,d) are due to III type of TFR behavior. The signal was sampled at 100 Hz for 50 s.}
\label{skeletons}
\end{figure*}

Considering synchrosqueezing, from the previous sections it is clear that the positions of the SWFT/SWT ridges are no closer to the actual frequencies than the WFT/WT ridges. Thus, ``curves'' in synchrosqueezed TFRs, although being more concentrated, are not located around the actual instantaneous frequencies of the components (though the latter can be fully recovered from the full TFS). This is illustrated in Fig. \ref{skeletons}, where the SWFT is compared with the ridge/direct WFT skeletons. For each time, the latter are constructed by partitioning the WFT into regions of unimodal amplitude (time-frequency supports $[\omega_-^{(m)}(t),\omega_+^{(m)}(t)]$), reconstructing from them the amplitudes $A^{(m)}(t)$, phases $\phi^{(m)}(t)$ and frequencies $\nu^{(m)}(t)$ using the chosen method, and then assigning $A^{(m)}(t)e^{i\phi^{(m)}(t)}$ to the frequency bin where the estimated frequency $\nu^{(m)}(t)$ lies; the WT skeletons can be constructed in the same way. For example, up to frequency discretization effects, the ridge-based WFT/WT skeleton is simply the WFT/WT with only peaks left (and multiplied by $2/\hat{g}(0)$), while other coefficients are set to zero. The MatLab codes for calculating TFR skeletons can be downloaded from \cite{freecodes} together with the other codes used in this work.

As can be seen from pairwise comparison of (b,f) and (c,g) in Fig. \ref{skeletons}, the SWFT is very similar to a simple ridge-based WFT skeleton (which is additionally more concentrated and easy to interpret): in both former and latter cases, the ``curves'' are located not around the true component frequency, but have similar deviations from it. On the other hand, the direct WFT skeleton in the noiseless case provides almost perfect representation (Fig. \ref{skeletons}(d)), being clearly superior to the SWFT or ridge skeleton, though the picture becomes more complicated when the noise is present (Fig.\ \ref{skeletons}(e-h)). However, both skeletons are constructed from the original WFT and obviously do not improve neither time, nor frequency, nor joint time-frequency resolution (as the accuracy of the parameters' estimates remains the same), providing advantages mainly in terms of visual appearance, similarly to the case of the SWFT.

\section{Conclusions}\label{sec:conclusions}

The results of this work can be summarized as follows:

\begin{enumerate}
\item The appropriate choice of window/wavelet resolution parameter $f_0$ is of crucial importance in time-frequency analysis. It determines the tradeoff between time and frequency resolutions of the TFR, with different choices leading to different TFR behaviors, and therefore different quantitative and qualitative results. If the frequency resolution is too high, the AM/FM component might be represented in the TFR as a number of independent tones, whereas if it is too low, then two interfering tones can be merged together and appear as a single component. We have considered and illustrated this issue on a numerous examples, and provided the conditions for each type of TFR behavior (see Tables \ref{tab:OVc}, \ref{tab:AMc} and \ref{tab:FMc}).

\item The optimal $f_0$ depends on the signal. Several adaptation schemes have been reviewed, but none of them is fully universal. The question of how best to select the appropriate window/wavelet parameters for a given signal remains open.

\item In the absence of an adaptation scheme, one can choose $f_0$ based on the desired resolution properties of the TFR. For the Gaussian window WFT and lognormal wavelet WT, in order to resolve two tones at frequencies $\nu_{1,2}$ with relative error $\epsilon$, one needs $f_0\geq \frac{2n_G(\epsilon)}{|\nu_2-\nu_1|}$ and $f_0\geq\frac{2n_G(\epsilon)}{2\pi|\log(\nu_2/\nu_1)|}$, respectively ($n_G(\epsilon)$ is the number of standard deviations within which the $1-\epsilon$ part of the normal distribution resides, e.g.\ $n_G(0.05)\approx2$, see Part I). Next, from (\ref{apfrel}),(\ref{irTridge}) it follows that to recover the AM/FM component with error $\lesssim\epsilon$ using the ridge method, one should choose $f_0\leq \sqrt{2\epsilon/\max{\Big[}\frac{\langle[A''(t)]^2\rangle^{1/2}}{\langle[A(t)]^2\rangle^{1/2}}, \langle2[\nu'(t)]^2\rangle^{1/2}{\Big]}}$ and $f_0\leq (2\pi)^{-1}\times$ \newline $\sqrt{2\epsilon/\max{\Big[}\frac{\langle[A''(t)/\nu^2(t)]^2\rangle^{1/2}}{\langle[A(t)]^2\rangle^{1/2}}, \langle2[\nu'(t)/\nu^2(t)]^2\rangle^{1/2}{\Big]}}$ for a Gaussian window and lognormal wavelet, respectively; provided $\epsilon$ is small enough (e.g.\ $\epsilon=0.05$), this will guarantee that the component is represented reliably in the TFR. Note that, in respect of different kinds of signals (having different characteristic frequency bands), the choice of $f_0$ for the WT seems to be slightly more universal than for the WFT (with the most widespread being $f_0=1$).

\item The relative performance of the direct and ridge reconstruction methods depends on the signal and the conformity of the TFR resolution properties with its structure. Direct estimates are exact when the component is reliably represented in the TFR and there is no noise or interference; ridge estimates are more robust to such complications (as well as to boundary distortions), but have inherent errors related to amplitude/frequency modulation. Hence, direct methods are to be preferred in the case of relatively clean signals with frequency components that are well-separated (as they appear in the TFR), while ridge reconstruction is superior for signals considerably corrupted by noise or with highly interfering components. We have suggested a simple automatic procedure for selection of the optimal reconstruction method in Sec.\ \ref{sec:optrecD}.

\item Synchrosqueezing does not provide significant advantages in terms of components' reconstruction, at least for the windows/wavelets which are unimodal in frequency, but it introduces additional complications. Thus, although being more concentrated, the SWFT/SWT actually has the same time and frequency resolutions as the WFT/WT from which it is constructed, and therefore does not offer the possibility of better tracking of parameters' time-variations or of the resolution of components that lie closer in frequency.

\end{enumerate}

\section*{References}

\bibliographystyle{elsarticle-num}
\bibliography{resolutionbib}

\end{document}